\documentclass[10pt]{amsart}

\usepackage{latexsym}
\usepackage{amsmath}
\usepackage{amssymb,amsthm}
\usepackage[english]{babel}
\usepackage{graphicx}
\usepackage{epic, eepic}
\usepackage{symbolesmagEng}
\usepackage{comment}
\usepackage{dessm6a}
\usepackage{dessm6cEng}
\usepackage{dessmagintroEng}
\usepackage{dessCS}
\usepackage{dessCNEng}

\title{Stabilization of the linear system of magnetoelasticity}
\author{Thomas Duyckaerts}
\date{\today}
\begin{document}
\begin{abstract}
We give a necessary and sufficient condition, of geometrical type, for the
uniform decay of energy of solutions of the linear system of
magnetoelasticity in a bounded domain with smooth boundary. A
Dirichlet-type boundary condition is assumed. When the
geometrical condition is not fulfilled, we show polynomial decay of
the energy, for smooth initial conditions. Our strategy is to use
micro-local defect measures to show suitable observability inequalities on
high-frequency solutions of the Lam\'e system. 
\end{abstract}
\maketitle



\section{Introduction}

\subsection{The system of magnetoelasticity}
Let $\Omega$ be a bounded, simply connected domain of $\RR^3$, with a
smooth boundary. Let us consider the following system, modelling the
displacement of a elastic solid in a magnetic field: 
\begin{equation}
\label{magneto}
\begin{aligned}
\left. \begin{gathered}
\partial _{t}^2 v - \mu \Delta v - (\lambda+ \mu ) \nabla \text{div} v
- \kappa \text{ rot}h \wedge {\Bbf}=0\\
\beta\partial _{ t}h + \text{ rot} \text{ rot} h - \beta \text{ rot }(
\partial _{t}v \wedge {\Bbf})=0\\
\div h=0 
\end{gathered}\right\} & (t,y)\in (0,\infty)\times\Omega\\
v=0,\quad 
h. \bfn=0,\quad 
\rot h \wedge \bfn=0  \quad&(t,y)\in(0,\infty)\times\partial\Omega, 
\end{aligned}
\end{equation}
where $v= (v_{1}, v_{2}, v_{3})$ is the displacement vector of the solid,
and $h=(h_{1}, h_{2}, h_{3})$ the magnetic field. The system is located in a
constant exterior magnetic field $\Bbf=
(B,0,0)$. We have denoted by $\Delta$, $\nabla$, div, curl respectively the Laplace
operator, gradient, divergence and curl operators according to the space
variable $y$, in the euclidian metric of $\RR^3$. The positive constants
$\kappa$ and $\beta$ are coupling constants, and $\bfn$ is the external
normal vector to the boundary of $\Omega$. The real Lam\'e constants $\lambda$
and $\mu$ are such that: $\lambda+2\mu>0$, $\mu>0$ and $\lambda+\mu \neq 0$.\par
The system (\ref{magneto}) has a natural time-decreasing energy:
\begin{gather*} 
E(t) =\frac 12 \int _{\Omega}| \partial_{t}v | ^2 + \mu | \nabla v| ^2 + (\lambda+\mu)|\div v|^2 + \kappa | h| ^2dy.
\end{gather*}
When $\Omega$ is simply connected, G.~Perla Menzala and E.~Zuazua have
showed that this energy tends to zero as time tends to infinity, which is a
simple consequence, using La Salle invariance principle, of the
non-existence of stationnary solution for (\ref{magneto}). The goal of
this paper is to give estimates on the speed of this convergence.\par
The system (\ref{magneto}) may be seen as a coupling between the Lam\'e system:
\begin{equation}
\label{lame.intro}
\partial_t^2 u-\mu \Delta u-(\lambda+\mu)\nabla \div u =0, 
\end{equation}
with Dirichlet boundary conditions, which is a conservative system, and the
following heat equation:
$$ \beta \partial_t  g -\Delta g=0.$$
The decay of energy is produced by this strongly dissipative equation. From the point of view of $v$, the
dissipation is caused by the coupling term: $R(v)\egaldef \rot(\partial_t v \wedge \Bbf)$.\par
Let us first consider the uniform decay with respect to initial condition
of the energy:
\begin{equation}
\label{stab.unif}
E(t)\leq f(t) E(0),\quad f(t) \underset{t\rightarrow +\infty}{\longrightarrow}0,
\end{equation}
where $f$ if independent of the initial condition.
In this case it is easy to show, using the semi-group property of the
equation (\ref{magneto}), that $f$ maybe taken as a negative exponential
function.\par
In paragraph \ref{intro.dec.unif} we state, with a technical hypothesis on $\Omega$, a necessary and
sufficient condition on the geometry of the problem for (\ref{stab.unif})
to hold. When this condition is not
fulfilled, there exist rays on $\Omega$, named $\Bbf$-resistant
rays, along which the energy of some solutions of (\ref{magneto})
concentrates, and the dissipative term $R(v)$ is very small. Indeed, when
such a ray exists, there is a sequence of solutions of (\ref{magneto})
concentrating on the ray and which is in first approximation parallel to $\Bbf$.\par
When there is no uniform stabilization we show (with the same technical
property on $\Omega$ than before), that solutions of (\ref{magneto}) decay
with polynomial speed for smooth enough initial data (cf paragraph
\ref{intro.dec.poly}). The speed of decay still depends on the geometry of
$\Omega$. In this case, the possible existence of boundary $\Bbf$-resistant
rays (i.e. living only in the boundary of $\Omega$) of infinite life-length
is the main obstacle to the decay.\par
Before giving more explicit results, let us mention some earlier works on
related subjects. As it was already stated, the convergence to $0$ for the
energy of magnetoelasticity in a bounded, simply connected domain was shown
by G.~Perla Menzana and E.~Zuazua in \cite{PeZu98}, but their method does not
give any information on the rate of convergence.  By energy methods, Mu\~
noz Rivera and Racke \cite{MuRa01}, Mu\~noz Rivera and de Lima Santos 
\cite{LiMu03} have shown the rate of convergence to be at least polynomial, in dimension $2$
or $3$, but only for some precise types of domains. Andreou and Dassios
\cite{AnDa97} have examined the same system on the entire space $\RR^3$,
showing again polynomial decay for some initial conditions.\par
The linear system of thermoelasticity has been more precisely
understood. In this system, the Lam\'e equations are coupled with a scalar
heat equation. The dissipation is caused by the longitudinal part of the
Lam\'e equation (the curl-free part of $v$). In \cite{LeZu98} and
\cite{BuLe99}, the authors give (under a spectral assumption) a necessary and sufficient condition on
$\Omega$, of geometrical nature, for the uniform decay in dimension $2$ or
$3$. Namely, this decay is equivalent to  
the non-existence of rays, called ``transversal polarization rays'', carrying the transversal component of $v$ (the
divergence-free component), which resists to the dissipation. In
\cite{LeZu98}, the authors also prove the polynomial decay in dimension $2$,
under the same spectral assumption, which is namely that the operator associated to the equation does not
admit any real eigenvalue. As shown in \cite{PeZu98}, this spectral condition is always fulfilled for
the system of magnetoelasticity in a bounded, simply-connected
domain.\par
The comparison of the two systems of thermo and magnetoelasticity show
that thermoelasticity is slightly less dissipative (the coupling of the
Lam\'e system with the heat equation is weaker), and more difficult to
describe, because of the non-trivial polarization of transversal waves. 
\subsection{Uniform decay}
\label{intro.dec.unif}
Assume that $\partial \Omega$ has no contact of infinite order with its
tangents. Thus, the hamiltonian flow of the symbol of a d'Alembertian
$\partial_t^2-c^2\Delta$, which is defined locally in
$S^*(\RR\times\Omega)$ (the spherical cotangent bundle of $\Omega$), maybe
extended until the boundary of this bundle to a global flow, the generalized bicharacteristic flow, wich may be seen as a continuous flow on
the spherical compressed cotangent bundle $S^*_b (\RR\times
\overline{\Omega})$ (cf \cite[chap. 24.3]{Hor}). We shall call
bicharacteristic rays or just rays the
characteristic curves of this flow. Such a curve $\gamma$
will be said  parallel to $\Bbf$ if its direction of propagation is always 
parallel to $\Bbf$ and orthogonal to $\Bbf$ if its direction of 
propagation is always orthogonal to $\Bbf$. We refer to section
\ref{chap.mesures} for the exact definitions of
$S^*_b(\RR\times\overline{\Omega})$ and of the generalized bicharacteristic
flow.\par
The Lam\'e system (\ref{lame.intro}) may be written as the sum of two wave
equations known as the longitudinal and transversal wave equations, of
respective speed  $c_L \egaldef \sqrt{\lambda+2\mu}$ and
$c_T \egaldef \sqrt{\mu}$ (cf paragraph \ref{par.Lame}). The assumption
$\lambda+\mu\neq 0$ is equivalent to $c_L\neq c_T$.
\begin{Def}
One calls {\bf longitudinal ray} (respectively {\bf transversal ray})
any bicharacteristic ray for the operator  $\partial_t^2-c_L^2\Delta$
(respectively $\partial_t^2-c_T^2\Delta$). One calls {\bf $\Bbf$-resistant ray}
any {\bf continuous} application $\gamma$ from an open interval $I=(s_0,s_n)$ to
$S_b^*(\RR\times \overline{\Omega})$ such that there exists a finite number
of reals $s_0<s_1<...<s_n$ such that:
\begin{itemize}
\item  on $(s_{j-1},s_j)$, $j\in \{1,...,n\}$, $\gamma$ is a longitudinal
  ray parallel to $\Bbf$, or a transversal ray orthogonal to $\Bbf$;
\item if $j\in \{1,...,n-1\}$, $\gamma(s_j)$ is an hyperbolic point for the
  longitudinal and transversal waves (cf
  paragraph \ref{var.car}) and one of the following assertions is true:
\begin{itemize}
\item ($L\rightarrow T$) case: $\gamma$ is a longitudinal ray on
  $]s_{j-1},s_j[$, and a transversal ray on $]s_j,s_{j+1}[$;
\item ($T\rightarrow L$) case: $\gamma$ is a transversal ray on
  $]s_{j-1},s_j[$, and a longitudinal ray on $]s_j,s_{j+1}[$.
\end{itemize}
\end{itemize}
(cf figure \ref{rayons})
\end{Def}
Near $s_j$, $1\leq j\leq n-1$, the continuity imposed by the definition of
$\gamma$ gives a condition on the angles of incidence and refraction. In the
case ($L\rightarrow T$), if we denote by $\alpha_L$ the angle between the
longitudinal incoming ray and the tangent to $\partial \Omega$ in the plane
of incidence, and by $\beta_T$ the angle between the transversal outcoming
ray and this tangent (cf figure \ref{rayons},~c), we have:
$$ \tan \alpha_L=\frac{c_T}{c_L},\quad \tan
\beta_T=\frac{c_L}{c_T},$$
(which implies $\alpha_L+\beta_T=\pi/2$).
In the case $(T\rightarrow L)$, and with similar notations, we have:
$$ \tan \alpha_T=\frac{c_L}{c_T},\quad \tan
\beta_L=\frac{c_T}{c_L}.$$
\begin{rem}
The $\Bbf$-resistant rays of figure \ref{rayons} are all planar, but this is
not a general property. 
\end{rem}
\begin{figure}[htbp]
\begin{center}
\dessintro
\end{center}
\caption{Examples of $\Bbf$-resistant rays}
\label{rayons}
\end{figure}
\begin{maintheo}
\label{th1}
Let $\Omega$ be a bounded, simply connected domain of $\RR^3$, with a
smooth boundary, having no contact of infinite order with its tangents.\\
The energy of the system of magnetoelasticity in $\Omega$ decays uniformly
if and only if there exists an $L>0$ such that every $\Bbf$-resistant ray
on $\Omega$ is of length at most $L$.
\end{maintheo}
\begin{rem}
As it will be shown in the proof, the transversal rays carry the component
of $v$ which is orthogonal to the direction of propagation, and the
longitudinal rays the component of $v$ which is parallel to this
direction. A $\Bbf$-resistant ray, whose direction of propagation is
orthogonal to $\Bbf$ in the transversal case and parallel to $\Bbf$ in the
longitudinal case carries essentially the component of $v$ which is
parallel to $\Bbf$, thus cancelling the dissipative term:
$$ R(v)\egaldef\rot(\partial_t v\wedge \Bbf).$$
From this point of view, the theorem \ref{th1} is very natural. 
\end{rem}
\begin{rem}
It is essential to assume $c_L \neq c_T$. Otherwise, the first equation in
(\ref{magneto}) would be a wave equation with wave speed $c_l=c_T$. Every
solution of (\ref{magneto}) such that:
$$ v_{\restriction t=0} \bot \Bbf,\quad \partial_t v_{\restriction t=0} \bot
\Bbf,\quad h_{\restriction t=0}=0$$ 
would be of constant energy.
\end{rem}
\begin{rem}
If $\Omega$ is not simply connected, there exists a finite dimensionnal
space $E$ of stationary solutions of (\ref{magneto}), whose components along
$v$ are null. The study of the decay to zero of the solutions may be
replaced by the study of their convergence to the eigenfunctions
corresponding to the space $E$ (cf \cite[chap. 5]{PeZu98}). We won't
develop this aspect here.
\end{rem}
\begin{rem}
The condition of uniform decay is not fulfilled in simple cases, like the
one of a bowl, but is generic in the class of $C^{\infty}$ open sets.
\end{rem}

\subsection{Polynomial decay} 
\label{intro.dec.poly}
Now we state a result of polynomial decay for initial data which are
sufficiently smooth.
The existence of a boundary $\Bbf$-resistant ray of infinite life-length is equivalent to the existence of a smooth closed curve
$\Gamma$ of $\partial \Omega$, included in a plane $\PPP$ normal to $\Bbf$,
boundary of a convex set of $\PPP$, and such that on $\Gamma$, ${\mathbf
  n}$ is normal to $\Bbf$ (cf figure \ref{rayons},~d). On such a curve,
$\Bbf$ stays tangential to the boundary. Let $\aleph(\Gamma)$ be the minimal
order
of contact of $\partial \Omega$ with a tangent parallel to
$\Bbf$. If such a curve $\Gamma$ exists, and if the boundary has no contact
of infinite order with its tangents, then:
$$ 2\leq \aleph(\Gamma)<\infty.$$
If $(v,h)$ is a sufficiently smooth solution of (\ref{magneto}), we shall
denote by   $E^{(j)}(t)$ the energy of order $j$ of $(v,h)$:
$$ E^{(j)}(t)\egaldef\frac{d^j}{dt^j} E(t).$$
Let $X_j$ be the subspace of $X$ of all initial data of
(\ref{magneto}) such that $E^{(j)}(0)$ is finite. It is exactly the domain
of $\AAA^j$, where $\AAA$ is the linear operator of magnetoelasticity
defined in section \ref{chap.A.F}.
\begin{maintheo}
\label{thdecpoly}
Let $\Omega$ be a bounded, simply connected domain with smooth boundary
having no contact of infinite order with its tangents.\\
 {\bf a)} Assume there is no boundary $\Bbf$-resistant ray on $\Omega$ of
 infinite life-length. Then:
$$ \exists C>0, \quad \forall V_0\in X_1,\quad \forall t\geq 0,\quad E(t) \leq \frac{C}{t+1} E^{(1)}(0).$$
{\bf b)} Assume on the contrary that such rays exist. Le
$\Gamma_1,...,\Gamma_M$ be the support of this rays, and:
$$ K\egaldef\sup_{m=1..M} \aleph(\Gamma_m).$$ 
Then:
$$ \exists C>0,\quad \forall V_0 \in X_{K},\quad \forall t\geq 0,\quad E(t) \leq
\frac{C}{t+1} E^{(K)}(0).$$
\end{maintheo} 
\begin{rem}
By an easy interpolation argument, one may deduce from {\bf b} the polynomial decay of any
solution with initial condition in $X_1$:
$$ \forall V_0 \in X_{1},\quad \forall t\geq 0,\quad E(t) \leq
\frac{C}{(t+1)^{\frac{1}{K}}} E^{(1)}(0).$$
\end{rem}
\begin{rem}
Theorem \ref{thdecpoly} completes the works of J.~E.~Mu\~ noz Rivera et
M.~De Lima Santos \cite{LiMu03} which show, for some types of domains of
$\RR^3$, a decay in $1/t$ for initial data in $E^{(7)}$. 
Note that the domains considered in their work (all of which have contacts
of infinite order with their tangents) do not fall within the scope of our article.
\end{rem}
The remainder of the paper is organized as follows. In section \ref{chap.A.F},
we reduce theorems \ref{th1} and \ref{thdecpoly} to high-frequency observability
inequalities on the Lam\'e system (\ref{lame.intro}). This is based on two
arguments: the setting aside of low frequencies, which is a consequence of
the non-existence of stationnary solution for the equation (\ref{magneto})
shown in \cite{PeZu98},
 and the decoupling, by simple calculations, of the two equations (the Lam\'e
 system and the heat equation) which compose (\ref{magneto}). In section
 \ref{chap.mesures}, we introduce micro-local defect measures (an object
 due to  P.~G\'erard \cite{PG91} and L.~Tatar \cite{Tar90}, and 
 in this particular setting to N.~Burq and G.~Lebeau \cite{BuLe99}), in
 order to study the lack of compactness of a sequence of
 high-frequency solutions of the Lam\'e system. The main result of this
 section (apart from the existence of the measures), is a propagation
 theorem which was stated and shown in \cite{BuLe99}. In section
 \ref{chap.CS}, we prove the observability inequality on solutions of the
 Lam\'e system (\ref{lame.intro}) which implies theorem \ref{th1}. The
 method of proof is to introduce, in a contradiction argument, a sequence of high
 frequency solutions of (\ref{lame.intro}) which contredicts this
 inequality, and to use propagation arguments on the defect measures of
 this sequence. Section \ref{chap.CN} is devoted to the necessary condition
 of theorem \ref{th1}, and is inspired by \cite{NBPG}: defect
 measures are used to construct
 a sequence of solutions of (\ref{lame.intro}) concentrating on a
 $\Bbf$-resistant ray and contradicting an
 observability inequality. Finally, in section \ref{chap.poly}, we prove by
 similar arguments than those of section \ref{chap.CS} an observability 
 inequality with loss of derivatives which implies the polynomial
 decay.\par
The author would very like to thank his thesis advisor, Nicolas Burq for his
invaluable help.

\section{Observability inequality for the Lam\'e system} 
\label{chap.A.F}
\subsection{Notations and preliminary results}
In this subsection are gathered a few basic facts about equations
(\ref{magneto}) and (\ref{lame.intro}), as well as some notations. The main
results of section \ref{chap.A.F} are stated in the next subsection.\par
If $U$ is an open set of $\RR^3$ or $\RR^4$ we set:
$$\Hv^s(U):=H^s(U,\CC^3), \quad \Lv(U):=L^2(U,\CC^3).$$
\subsubsection{Magnetoelasticity}
Consider the following spaces:
\begin{gather*}
H:=\big{\{} g \in \Lv(\Omega), \; \div g=0 \text{ in }\Omega,\; g.n=0 \text{ in } \partial \Omega\big\}  \\
\Hv_0^1:=\{f \in \Hv^1(\Omega), \; f=0 \text{ in } \partial \Omega \}\\
H':=\big{\{}f\in H\cap\Hv^2(\Omega),\; \rot f\wedge n=0 \text{ in }\partial \Omega\big{\}}, 
\end{gather*}
and the following norms:
\begin{equation*} 
\|g\|^2_H:=\kappa \|g\|^2_{\Lv(\Omega)}, \quad \|f\|^2_{\Hv^1_0(\Omega)}:=(\lambda+\mu)\|\div f\|^2_{L^2(\Omega)}+\mu\|\nabla f\|^2_{\left(\Lv(\Omega)\right)^3}.
\end{equation*}
Let $\AAA$ be the unbounded operator on
$X:=\Hv_0^1(\Omega)\times\Lv(\Omega)\times H$, with domain $D(\AAA)$,
defined by:
\begin{gather*} 
\AAA(V_0)\egaldef \left( \begin{array}{c} -v_1  \\ - \Delta_e v_0 -\kappa
    (\rot h_0) \wedge \Bbf \\ -\rot(v_1 \wedge \Bbf) + \frac{1}{\beta} \rot
    \rot h_0  \end{array}\right)  \\  D(\AAA)\egaldef (\Hv^2\cap\Hv^1_0)\times\Hv^1_0\times H'. 
\end{gather*}
where $V_0=(v_0,v_1,h_0)$ denotes an element of $X$.
Equation (\ref{magneto}) may be rewritten: 
\begin{equation}
\label{magneto2}
\partial_t V+\AAA V=0,\quad V=(v,\partial_t v,h).
\end{equation}
The following proposition is due to G.~Perla Menzala and E.~Zuazua \cite{PeZu98}:
\begin{prop}
a) The operator $\AAA$ is maximal accretive. For any initial data $V_0\in
X$, there exists an unique weak
solution $V(t)=(v(t),\partial_t v(t), h(t))\in C^0([0,+\infty[;X)$
of \eqref{magneto2} such that $V(0)=V_0$. Functions $v$ and $h$ are
solutions in the distributional sense of the three first lines of system
\eqref{magneto}.\\
b) The energy $E(t)$:
\begin{equation}
\label{defEu}
E(t) := \frac 12 \|V(t)\|^2_X=\frac 12 \int _{\Omega}| \partial_{t}v | ^2 + \mu | \nabla v| ^2 + (\lambda+\mu)|\div v|^2 + \kappa | h| ^2dy
\end{equation}
is decreasing. More precisely:
\begin{equation}
\label{dEdt}
\forall t\geq 0,\quad E(t)-E(0)= -\frac{\kappa}{\beta} \int_0^t \int _{\Omega} | \text{rot} h| ^2 dy.
\end{equation}
c) If $\Omega$ is simply connected:
$$ \forall V_0 \in X,\quad E(t) \underset{t\rightarrow +\infty}{\longrightarrow } 0. $$
\end{prop}
The assertions a) and b) are straightforward applications of the semi-group
theory for the oprator $\AAA$. The assertion c) is a consequence of the
non-existence of stationnary solutions for the system.
\subsubsection{Lam\'e system}
Let us now consider the Lam\'e system with Dirichlet boundary conditions:
\begin{equation}
\label{Lame1}
\begin{gathered}
\partial_t^2 u-\Delta_e u=0 \text{ in } \RR\times \Omega\\ 
u_{\restriction \partial \Omega}=0\\
(u_{\restriction t=0},\partial_t u_{\restriction t=0})=(u_0,u_1).
\end{gathered}
\end{equation}
Let $X_e$ be the space $\Hv^1_0\times \Lv$ and $\LLL$ the unbounded
operator on $X_e$ defined by:
\begin{equation}
\LLL\egaldef \left[\begin{array}{cc} 0 &-\text{Id} \\ -\Delta_e & 0
  \end{array}\right] \qquad D(\LLL)\egaldef \Hv^2\cap\Hv_0^1\times\Lv.
\end{equation}
Taking $(u_0,u_1)$ in the energy space $X_e$, the equation (\ref{Lame1}) may be written:
\begin{equation}
\label{Lame2}
\partial_t U+\LLL U=0, \; U(t)=(u,\partial_t u).
 \end{equation}
\begin{prop}
The operator $\LLL$ is maximal and unitary. For any initial data
$U_0=(u_0,u_1)\in X_e$, the system (\ref{Lame2}) has an unique weak
solution $U\in C^0(\RR,X_e)$. Furthermore, the function $u$ is a solution
of (\ref{Lame1}) in the distributional sense. At last, the energy: 
$$\frac 12 \left(\|u(t)\|^2_{\Hv^1}+\|\partial_t
  u\|^2_{\Lv}\right)=\frac 12 \|U\|_{X_e}^2$$ 
of this solution is constant.
\end{prop}  
\subsubsection{Two useful lemma}
The two following standard lemma will be of great help in all this
paper. The first one is due to the fact that $\Omega$ is simply connected 
(cf \cite[Appendix I, lemma 1.6]{Te79}):
\begin{lem}
 \label{divrot}
The $\Hv^1$ norm on $H\cap \Hv^1$ is equivalent to the norm: $\|u\|\egaldef
\|\rot u\|_{\Lv}$.
\end{lem}
The second lemma is a elementary energy estimate on solutions of the
non-homogeneous Lam\'e system. If $w(t)$ is a function with values in some
Hilbert space, we set: $(w_0,w_1)\egaldef (w,\partial_t w)_{\restriction_{t=0}}$.
\begin{lem}
\label{In En Lame}
Let $T>0$, $W\in C^0((0,T),X_e)$ and $F \in L^2((0,T),X_e)$ such that:
\begin{equation}
\label{Hyp In En Lame}
\partial_t W+\LLL W=F, \quad t\in (0,T).
\end{equation}
Then:
\begin{equation*}
\int_0^T \| W(t)\|^2_{X_e}\, dt\leq C\left\{\|W(0)\|^2_{X_e}+\int_0^T \|F(t)\|^2_{X_e} dt \right\},
\end{equation*}
where $C$ only depends on $T$.
In particular, if:
\begin{gather*}
 w\in C^1((0,T),\Lv(\Omega))\cap C^0((0,T),\Hv^1_0(\Omega)),\quad f\in \Lv((0,T)\times\Omega)\\
\partial_t^2 w-\Delta_e w= f.
\end{gather*}
Then:
$$ \|w\|^2_{\Hv^1((0,T)\times\Omega)} \leq C \left(\|w_0\|_{\Hv^1(\Omega)}^2+\|w_1\|_{\Lv(\Omega)}^2 + \|f\|^2_{\Lv((0,T)\times\Omega)} \right).$$
\end{lem}
\begin{proof}
To prove the first inequality, we may suppose $W_0\in D(\LLL)$. The $X_e$
scalar product of (\ref{Hyp In En Lame}) with $W$ gives
\begin{align*}
{\rm Re } (\partial_t W,W)_{X_e}+\underbrace{{\rm Re} (\LLL W,W)_{X_e}}_{0}
& =(W,F)_{X_e}\\
\frac 12 \frac{d}{dt} \|W(t)\|^2_{X_e} & \leq \|W(t)\|_{X_e}\left\|F(t)\right\|_{X_e}\\
\frac{d}{dt} \|W(t)\|_{X_e} & \leq \|F(t)\|_{X_e}\\
\|W(s)\|_{X_e} & \leq \|W_0\|_{X_e}+ \int_0^s \|F(t)\|_{X_e} \,dt,\quad s\in (0,T)\\
\|W(s)\|^2_{X_e} & \leq C\left\{ \|W_0\|^2_{X_e}+s\int_0^s\|F(t)\|^2_{X_e}\,dt\right\}.
\end{align*}
(the last line is a consequence of Cauchy-Schwarz inequality).
Next, we bound, in the right member of the inequality the integral from $0$
to $s$ by the same integral from $0$ to $T$, and we integrate with respect
to $s$ between $0$ and $T$, which yields the first part of the lemma. The
second part is an easy consequence of it.
\end{proof}
\subsection{Results}
\label{A.F.result}
\begin{prop}[Uniform decay]
\label{CNSineg}
Let $\Omega$ be a smooth, simply connected, bounded domain of $\RR^3$.
\begin{enumerate}
\item[{\bf a)}] 
Assume that there exist $T>0$ and $C>0$ such that for any solution $U$ of (\ref{Lame2}):
\begin{equation}
\label{inegCS}
\|u_0\|^2_{\Hv^1_0}+\|u_1\|^2_{\Lv}\leq C \left( \|\rot(\partial_t u\wedge \Bbf)\|^2_{\Hv^{-1}((0,T)\times \Omega)} +\|u_0\|^2_{\Lv}+\|u_1\|^2_{\Hv^{-1}}\right)
\end{equation}
then the energy of solutions of the system of magnetoelasticity
(\ref{magneto}) decays uniformly with respect to initial data.
\item[{\bf b)}] 
Conversely, if the energy of solutions of (\ref{magneto}) decays uniformly,
then there exist $T>0$ and $C>0$ such that for any solution of
(\ref{Lame1}) of finite energy:
\begin{equation}
\label{inegCN}
\|u_0\|^2_{\Hv_0^2}+\|u_1\|^2_{\Lv} \leq C\int_0^T \|\partial_t u \wedge \Bbf \|^2_{\Lv(\Omega)}.
\end{equation}
\end{enumerate}
\end{prop}
\begin{rem}
The two inequalities (\ref{inegCS}) and (\ref{inegCN}) are indeed equivalent
(by theorem \ref{th1}).
\end{rem}
Let $U=(u,\partial_t u)$ be a solution to (\ref{Lame2}) with initial data
$U_0\in D(\LLL^N)$. Set:
$$ \QQQ^N_T(u)\egaldef\sum_{l=0}^N \|\rot (\partial_t^{l+1} u \wedge
B)\|^2_{\Hv^{-1}\left((0,T)\times\Omega\right)}.$$
\begin{prop}[A sufficient condition of polynomial decay]
\label{decpoly}
Let $\Omega$ be a bounded, simply connected domain of $\RR^3$. Assume that
there exist $T>0$, $C>0$ and an integer $N\geq 1$ such that for every
solution of the Lam\'e system \eqref{Lame2} with initial data:
\begin{equation}
\label{regpoly1}
U_0=(u_0,u_1)\in D(\LLL^N),
\end{equation}
the following inequality holds:
\begin{equation}
\label{inegpoly}
\|u_0\|^2_{\Hv_0^1}+\|u_1\|^2_{\Lv} \leq C\left(\QQQ_T^N (u)+\|u_0\|^2_{\Lv}+\|u_1\|^2_{\Hv^{-1}}\right).
\end{equation}
Then there exists $C>0$ such that for every solution $V$ of
\eqref{magneto2} with initial date $V_0\in D(\AAA^N)$, and for all positive $t$,
$$ \|V(t)\|_X^2 \leq \frac{C}{t+1}\|V_0\|^2_{D(\AAA^N)} $$
\end{prop}
\begin{rem}
results such as propositions \ref{CNSineg} and \ref{decpoly} are fairly
classical in this setting. To prove them, we shall avoid the usual abstract
decoupling argument (see \cite{LeZu98}) but rather use simple energy
estimates on systems of magnetoelasticity and Lam\'e.
\end{rem}

\subsection{Uniform decay}
We prove here the proposition \ref{CNSineg}. We first write a necessary and
sufficient condition of uniform decay for solutions of a general
dissipative equation. The second step of the proof consists in applying this
condition to the system of magnetoelasticity, furthermor decoupling it in the
system of Lam\'e and an heat equation.
\subsubsection{Abstract framework}
\label{cadre abstrait}
Let $\PPP$ be a maximal, accretive operator on an Hilbert space $X$, with
dense domain $D(\PPP)$. Denote by $\|..\|$ the norm of $X$, $\|..\|_1$ the
natural norm of $D(\PPP^1)$ and $\|..\|_{-1}$ the norm of its dual space,
with respect to the pivot space $X$.
Assume the embedding:
 $$X \longrightarrow D(\PPP)' $$ 
is compact. For $z_0\in X$, we will denote by $z(t)$ the solution (obtained
for example by standard semi-group theory) of:
\begin{equation}
\label{dvdt+Pv}
\frac{d\,z}{dt}+\PPP z=0,\quad z_{t=0}=z_0 
\end{equation}
By accretivity of $\PPP$, the energy $\frac 12 \|z\|^2$ is
time-decreasing. The following uniqueness-compactness argument is by now
classical (cf \cite{BaLeRa92}):
\begin{lem}
\label{unifgen}
The two following assertions (i) and (ii) are equivalent:
\begin{align}
\tag{i}
& \exists C>0, \; \exists a>0, \;\forall z_0 \in X,\; \forall t>0, \;\|z(t)\|^2 \leq C \|z_0\|^2 e^{-at}\\ 
\notag
& \text{(the energy is uniformly decreasing)}\\
\tag{ii} &\text{a) } \exists T>0, \; \exists C>0,\; \forall z_0 \in X, \; \|z(T)\|^2 \leq C\left( \|z(0)\|^2-\|z(T)\|^2 + \|z(0)\|_{-1}^2 \right)\\
\notag
&\text{b) }\text{There is no non-zero solution of (\ref{dvdt+Pv}) of
  constant energy on $[0,+\infty[$.}
\end{align}
\end{lem}
\begin{corol}
\label{corolCNS}
The energy of  (\ref{magneto2}) is uniformly time-decreasing if and only if:
$$ \exists T>0,\;\exists C>0,\; \forall V_0 \in X,\; 
E(T) \leq C\left\{\int_0^T \int_\Omega |\rot h|^2\,dy\,dt+\|V_0\|^2_{D(\AAA)'}\right\}. $$
\end{corol}
Indeed, the non-existence of stationnary solution (the condition ii,b of
lemma \ref{unifgen} has been proved in \cite[p.356]{PeZu98}), which shows
the corollary.
\begin{proof}[Proof of lemma \ref{unifgen}]
It is easy to see that (i) may be replaced by:
\begin{equation}
 \tag{i' } \exists T>0, \; \exists C>0,\; \forall z_0 \in X, \quad
 \|z(T)\|^2 \leq C\left( \|z(0)\|^2-\|z(T)\|^2\right) 
\end{equation}
Clearly (i') implies (ii).\par
Assume (ii). For some $T>0$, set:
$$ 
 q_T(z)\egaldef\|z(0)\|^2-\|z(T)\|^2, \quad
G^T\egaldef\{ z_0 \in X, \; q_T(z)=0 \},$$
which is the kernel of a positive, bounded, quadratic form on $X$, thus a
closed subspace of $X$.\par  
According to (ii), a), and the compactness of the embedding from $X$ to
$D(P)'$, $G(T)$ is locally compact thus of finite dimension, for large $T$.
By assumption b),
$$ \bigcap_{T\geq 0} G^T=\{0\}. $$
Consequently, ${\rm dim}\, G^T$ being a time decreasing function of $T$,
when $T$ is large enough:
\begin{equation}
\label{GTvide}
G^T=\{0\}
\end{equation}
Let's fix such a $T$. The quadratic form $q_T$ is positive definite so that
its square root $\sqrt{q_T}$ is a pre-hilbertian norm on $X$, bounded from
above by the natural norm of $X$. Assume $(i')$ does not hold. Then there
exists a sequence $(z_0^k)$ of elements of $X$ such that:
\begin{equation}
\label{vkT=1}
 1=\|z^k(T)\|^2, \qquad \lim_{k \CVF +\infty} q_T(z^k)=0. 
\end{equation}
This implies that $\|z_0^k\|$ is bounded. Thus, we may extract form
$(z_0^k)$ a subsequence, which we will again denote by  $(z_0^k)$, such
that:
$$ z_0^k \tendkf z_0 \in X, \text{ weakly in }X.$$ 
Let $\varphi_T$ be the hermitian product given by $q_T$. We have:
$$ \lim_{k\CVF +\infty} \varphi_T(z^k,z) \tendk q_T(z),$$
which implies, with (\ref{vkT=1}), that $q_T(z)=0$ and thus, using
(\ref{GTvide}) that $z=0$. The compactness of the embedding of $X$ in
$D(P)'$ yields:
$$ \lim_{k\CVF +\infty} \|z_0^k\|_{-1} = 0. $$
Using $a)$ and (\ref{vkT=1}) we obtain the following contradictory assertion:
$$ 1\leq C\left(q_T(z^k)+\|z_0^k\|^2_{-1}\right)=o(1) \text{ quand } k \CVF +\infty. $$
\end{proof}

\subsubsection{Proof of proposition \ref{CNSineg}}
Assume the uniform time-decay of the energy of solutions of
\eqref{magneto2}. Then, by (\ref{dEdt}), there exist $T>0$ and $C>0$ such
that the following estimates hold for any solution $v$ of (\ref{magneto2}):
\begin{equation}
\label{inegdemCN}
\|v_0\|^2_X \leq C\int_0^T \|\rot h(t)\|_{\Lv(\Omega)}^2\,dt.  
\end{equation}
Let $U$ be a solution of the Lam\'e system with initial data
$U_0=(u_0,u_1)\in D(\LLL)$ and $V$ the solution of the system of
magnetoelasticity with initial data:
$$ V_0=(v,\partial_t v,h)_{\restriction t=0}=(u_0,u_1,0). $$
Set: $W(t)\egaldef V(t)-\big(u(t),\partial_t u(t),0\big)$. Then:
\begin{equation*}
 \partial_t W+\AAA W=\big(0,0,-\rot (\partial_t u\wedge \Bbf) \big) 
\end{equation*}
Take the scalar product in $X$ with $W$ of the two side of this equality,
then inegrate the real part with respect to time between $0$ and $T$. Using: 
\begin{gather*}
 \Re (\AAA W,W)_X=\frac{\kappa}{\beta} \| \rot h\|^2_{\Lv(\Omega)}\\
\big(\rot(\partial_t u\wedge \Bbf),h\big)_{\Lv(\Omega)}=\left(\partial_t u
\wedge \Bbf,\rot h\right)_{\Lv(\Omega)},
\end{gather*}
the fact that $W_{\restriction t=0}=0$, and in the second line,
the inequality (\ref{inegdemCN}), we get:
\begin{gather*}
\|W(T)\|^2_{X}+ \int_0^T \|\rot h(t)\|^2_{\Lv(\Omega)}\,dt \leq C\int_0^T \|\partial_t u(t) \wedge \Bbf\|^2_{\Lv(\Omega)}\,dt\\
\|u_0\|^2_{\Hv^1_0(\Omega)}+\|u_1\|^2_{\Lv(\Omega)}\leq C\int_0^T
\|\partial_t u(t) \wedge \Bbf\|^2_{\Lv(\Omega)}\, dt  
\end{gather*}
This shows point {b)}. To prove $a)$, assume that inequality
(\ref{inegCS}) holds. Consider a solution  $V=(v,\partial_t
v,h)$ of (\ref{magneto2}) with initial data $V_0=(v_0,v_1,h_0)$, and the
solution $u$ of Lam\'e system with initial data:
$$(u,\partial_t u)_{\restriction t=0}=(v_0,v_1). $$
Thus, by (\ref{inegCS}):
\begin{equation}
\label{interm1}
\|v_0\|^2_{\Hv^1_0}+\|v_1\|^2_{\Lv}\leq C \left( \|\rot(\partial_t u\wedge \Bbf)\|^2_{\Hv^{-1}((0,T)\times \Omega)} +\|v_0\|^2_{\Lv}+\|v_1\|^2_{\Hv^{-1}}\right).
\end{equation}
Furthermore, the energy inequality on the non-homogeneous Lam\'e system (lemma
\ref{In En Lame}) yields:
\begin{gather*}
\|u-v\|^2_{\Hv^1((0,T)\times\Omega)}\leq C\int_0^T \|\rot h(t) \|^2_{\Lv}\,dt\\
\|\rot(\partial_t u\wedge \Bbf)\|^2_{\Hv^{-1}((0,T)\times \Omega)}\leq C \left\{\int_0^T \|\rot h(t) \|^2_{\Lv}\,dt+\|\rot(\partial_t v\wedge \Bbf)\|^2_{\Hv^{-1}((0,T)\times \Omega)} \right\},
\end{gather*}
which implies, using (\ref{interm1}), lemma \ref{divrot}) and the following
equation:
\begin{gather}
\notag
 \beta\partial _{ t}h + \rot  \rot h= \beta \text{ rot }(
\partial _{t}v \wedge \Bbf),\\
\label{interm2}
\|v_0\|^2_{\Hv^1_0}+\|v_1\|^2_{\Lv}\leq C \left\{\int_0^T \|\rot h(t) \|^2_{\Lv}dt
+\|v_0\|^2_{\Lv}+\|v_1\|^2_{\Hv^{-1}}\right\}.
\end{gather}
In order to use corollary \ref{corolCNS}, we need to add to the left side of
inequality (\ref{interm2}) the $\Lv$-norm of $h(T)$. We may do so by taking
a larger $T$. Indeed, consider $s\in [0,T]$ such that: 
$$ \|h(s)\|^2_{\Lv}=\min_{t\in[0,T]}\|h(t)\|^2_{\Lv}.$$
Lemma \ref{divrot} gives an $\alpha>0$ such that:
\begin{gather*}
 g\in H\cap \Hv^1 \Rightarrow \|g\|^2_{\Lv(\Omega)} \leq \alpha \|\rot g\|^2_{\Lv(\Omega)} \\
\begin{aligned}
 \|h(s)\|^2_{\Lv(\Omega)} & \leq  \frac 1T \int_0^T \|h(t)\|^2_{\Lv(\Omega)} dt\\
& \leq  \frac{\alpha}{T} \int_0^T \| \rot h(t) \|^2_{\Lv(\Omega)} dt.
\end{aligned}
\end{gather*}
Inequality (\ref{interm2}) taken with initial time $t=s$ yields:
\begin{multline*}
E(s)=\frac 12 \left( \|u(s)\|^2_{\Hv_0^1(\Omega)}+\|\partial_t
  u(s)\|^2_{\Lv(\Omega)}+\kappa \|h(s)\|^2_{\Lv(\Omega)}\right) \\
\leq C\left\{ \int_0^{2T} \|\rot h(t)\|^2_{\Lv(\Omega)}dt + \|u(s)\|^2_{\Lv(\Omega)}+\|\partial_t u(s)\|^2_{\Hv^{-1}(\Omega)} \right\}.
\end{multline*}
The energy $E$ being time-decreasing, this implies the inequality of
corollary \ref{corolCNS}, and so the uniform decay of solutions of
(\ref{magneto}). The proof of b) is complete.

\subsection{Polynomial decay}

\subsubsection{Abstract framework}
We shall use here the notations of paragraph \ref{cadre abstrait} Let $N$
be a positive integer, and $Q_T$ the quadratic form defined by:
$$ Q_T(z)\egaldef\sum_{l=0}^N\left(\|\partial_t^l z(0)\|^2-\|\partial_t^l z(T)\|^2\right)=\sum_{l=0}^N q_T(\partial_t^l z),\quad D(Q_T)=D(\PPP^N).$$
The function $\partial_t^l z$ being a solution of (\ref{dvdt+Pv}), its
energy decays with time, so that $Q_T$ is positive. Recall the definition:
(cf \cite{RS})
\begin{Def}
\label{fermable}
The quadratic form $Q_T$ is said to be {\bf closable} when the closure
$X_Q^T$ of $D(Q_T)$ in $X$ for the norm:
$$ \|z_0\|_{Q_T}=\sqrt{\|z_0\|^2+Q_T(z)} $$
is complete for this norm.
\end{Def}
\begin{rem}
\label{critere de fermabilite}
This is equivalent to the fact that for all Cauchy sequence in $D(Q_T)$ for the norm
$\|.\|_{Q_T}$, $(z^k)$, converging to $0$ in $X$, we have:
$$ \lim_{k\rightarrow +\infty} Q_T(z^k)=0. $$ 
\end{rem}
We shall again assume the compactness of the embedding: $ X \longrightarrow D(\PPP)'. $.\par
The following classical argument goes back to Russel \cite{Ru75}.
\begin{lem}
\label{genpoly}
Under the following assumptions:
\begin{enumerate}
\item[{\bf a)}]
there exist $T,C>0$ such that:
\begin{equation}
\label{CSpoly1}
\forall z_0 \in D(\PPP^N),\quad \|z(T)\|^2\leq
C\left(\|z_0\|^2_{-1}+Q_T(z) \right);
\end{equation}
\item[{\bf b)}]
system (\ref{dvdt+Pv}) have no non-zero solution of constant energy on $[0,+infty[$;
\item[{\bf c)}] 
the quadratic form $Q_T$ is closable.
\end{enumerate}
There exists $C>0$ such that:
\begin{equation}
\label{decpoly1}
\forall z_0 \in D(\PPP^N), \; \forall t\geq 0, \quad \|z(t)\|^2 \leq \frac{C}{t+1} \|z_0\|^2_N.
\end{equation}
\end{lem}
\begin{proof}
We shall first use a compactness argument similar to that of proposition
\ref{unifgen}. Let $T$ be a large positive real number, such that
(\ref{CSpoly1}) holds. Consider $X_Q^T$, the subspace of $X$ introduced in
definition (\ref{fermable}). Extending $Q_T$ to $X_Q^T$ by continuity, we
can still write inequality (\ref{CSpoly1}) for $v_0\in X_Q^T$.
Consider the following closed subspace of $X_Q^T$: 
$$ J^T\egaldef \{ z_0 \in X_Q^T, \; Q_T(z)=0 \}. $$
By (\ref{CSpoly1}) we have, for any $z_0 \in J^T$:
$$ \|z_0\|_{Q_T} \leq C \|z_0\|_{-1}.$$
Using assumption $b)$ as in the proof of proposition \ref{unifgen}, we
obtain that for $T$ large enough:
\begin{equation}
\label{JTvide}
J^T=\{ 0\}.
\end{equation}
From now on, $T$ will be taken such that (\ref{JTvide}) holds. The same
process as in the proof of proposition \ref{unifgen} yields:
\begin{equation}
\label{CSpoly2}
\forall z_0 \in X_Q^T,\qquad \|z(T)\|^2\leq C\,Q_T(z).
\end{equation}
Elsewhere, there would exist a sequence $(z_0^k)$ of elements $X_Q^T$ such that:
\begin{equation}
\label{Q_T0}
 \|z^k(T)\|^2=1,\quad
 \lim_{k\CVF +\infty} Q_T(z^k)= 0. 
\end{equation}
Up to a subsequence, we may assume that $(z_0^k)$ converges weakly to $0$
in $X_Q^T$. By (\ref{JTvide}) and (\ref{Q_T0}), $z_0=0$. By compactness:  
$$ \|z_0^k\|_{-1} \tendk 0.$$
In view of (\ref{Q_T0}), this contredicts (\ref{CSpoly1}).
\par 
By the triangle inequality, we have:
$$\sqrt{q_T(x+y)}\leq \sqrt{q_T(x)}+\sqrt{q_T(y)}.$$ 
Noting that $\left(\Id+\PPP\right)^l$ is an isomorphism from 
$D\left(\PPP^l\right)$ to $X$, it is easy to show:
\begin{align*}
Q_T(z) &\leq C\sum_{l=0}^N q_T\left((\Id-\partial_t)^l z\right) \\
    &\leq C\sum_{l=0}^N \left(\|z(0)\|^2_l-\|z(T)\|^2_l\right).
\end{align*}
From this and (\ref{CSpoly2}), we deduce:
\begin{align}
\notag
z_0\in D(\PPP^N) &\Rightarrow \|z(T)\|^2\leq C \sum_{l=0}^N (\|z(0)\|^2_l-\|z(T)\|^2_l)\\
\notag
 z_0 \in X &\Rightarrow  \|z(T)\|_{-N}^2\leq C \sum_{l=0}^N
(\|z(0)\|^2_{l-N}-\|z(T)\|^2_{l-N})\\
\label{inegpoly2}
 z_0 \in X &\Rightarrow \|z(T)\|_{0}^4\leq C \|z(T)\|^2_N \sum_{l=0}^N
\left(\|z(0)\|^2_{l-N}-\|z(T)\|^2_{l-N}\right)
\end{align}
The last inequality follows from the trivial bound:
$\|f\|^2\leq\|f\|_N\|f\|_{-N}$.\\
Set: $\beta_n \egaldef \sum_{j=0}^N\|z(nT)\|^2_{-j}$. We have:
\begin{equation*}
\beta_n\leq (N+1) \|z(nT)\|^2.
\end{equation*}
Thus, using (\ref{inegpoly2}) with the initial data 
$t=nT$ instead of $t=0$, there exists a constant $C>0$ such that:
$$ \beta_{n+1}^2 \leq C\|z(0)\|^2_N(\beta_n-\beta_{n+1}). $$
The following lemma (standard in this setting, see \cite[lemma 2.1]{Ru75} completes the proof of lemma \ref{genpoly}:
\begin{lem}
\label{lem.suites}
Let $(\beta_n)_{n\in\NN}$ be a sequence of positive real numbers and
$M_0>0$ such that
$$\beta_{n+1}^2 \leq M_0 (\beta_{n}-\beta_{n+1}).$$
Then:
$$\forall n\geq 1,\; \beta_n \leq \frac {2M_0}{n}. $$
\end{lem}
Inequality (\ref{inegpoly2}) implies, with lemma \ref{lem.suites}:
$$ \|z(nT)\|^2\leq \frac{C}{T} \|z(0)\|^2_N,$$
from which we deduce, taking into account the decay of the energy of $z$ that
(\ref{decpoly1}) holds.
\end{proof}
\begin{proof}[Proof of lemma \ref{lem.suites}]
Let $\alpha_n\egaldef n\beta_n$. Then:
$$ \alpha_n-\alpha_{n+1} \geq \frac{\alpha_{n+1}}{n+1}\left(\frac{n}{M_0(n+1)}\alpha_{n+1}-1\right). $$
In particular:
$$ \alpha_{n+1} > 2M_0 \Rightarrow \alpha_n > \alpha_{n+1} $$
Assume there exists at least one integer $n$ such that
$\alpha_{n+1}>2M_0$. Let $N$ be the smallest of these integers. Then:
$$\alpha_n>\alpha_{n+1}>2M_0,$$ 
which contredicts the minimality of $N$ when $N\geq 1$ or the fact that $a_0$
is null when $N=0$.
\end{proof}
\subsubsection{Proof of proposition (\ref{decpoly})}
Let $V$ be a solution of (\ref{magneto2}), 
with initial data in $D(\AAA^N)$. The first step of the proof is to
approach $V$ by a solution $U$ of the Lam\'e system.\par
It is easy to see that $(\partial_t^N v,\partial_t^{N+1}v)_{\restriction
  t=0}$ is in $\Hv_0^1\times \Lv$. The operator $\Delta_e$ being an isomorphism
from $\Hv_0^1\cap \Hv^2$ to $\Lv$, the operator $\LLL$ is an isomorphism
from $D(\LLL)$ to $\Hv_0^1\times \Lv$. As a consequence, we may choose
$(u_0,u_1)$ such that: 
\begin{equation*}
(u_0,u_1)\in D(\LLL^N),\quad(-\LLL)^N(u_0,u_1)=(\partial_t^N v, \partial_t^{N+1} v)_{\restriction t=0}.
\end{equation*}
The corresponding solution of the Lam\'e system $U=(u,\partial_t u)$ satisfies:
\begin{equation*}
\left\{ 
\begin{array}{l}
\partial_t^{N+1} u_{\restriction t=0}=\partial_t^{N+1} v_{\restriction t=0}\\
\partial_t^N u_{\restriction t=0}=\partial_t^N v_{\restriction t=0}. 
\end{array} \right.
\end{equation*}
Set $w=v-u$. We will first show:
\begin{equation}
\label{inegU}
\sum_{l=0}^{N} \|\partial_t^l w\|^2_{\Hv((0,T)\times\Omega)} \leq C
\sum_{j=0}^{N} \| \partial_t^j \rot h \|^2_{\Lv((0,T)\times\Omega)}.
\end{equation}
We have:
\begin{equation}
\label{eqpolyU1}
\partial^2_t w-\Delta_e w=\kappa (\rot h)\wedge \Bbf,\quad
\left( \partial_t^N w,\partial_t^{N+1} w\right)_{\restriction t=0}=(0,0).
\end{equation}
Equation (\ref{eqpolyU1}) implies, for $l\in \{0,..,N-1\}$:
\begin{equation}
\label{ineg w/h}
 \|\Delta_e \partial_t^l w_{\restriction t=0}\|_{\Lv}\leq \|\partial_t^{l+2} w_{\restriction t=0}\|_{\Lv}+\|\rot \partial_t^l h_{\restriction t=0}\|_{\Lv}. 
\end{equation}
But any $g\in \Hv^2(\Omega)\cap \Hv_0^1(\Omega)$ satisfies: 
$$ \|g\|_{\Hv^1_0}\leq \|g\|_{\Hv^2}\leq C\|\Delta_e g\|_{\Lv},$$
which yields, with (\ref{ineg w/h}):
\begin{equation}
\label{ineg.w/h'}
\|\partial_t^l w_{\restriction t=0}\|_{\Hv_0^1}\leq
  C\left(\|\partial_t^{l+2} w_{\restriction t=0}\|_{\Lv}+\|\rot \partial_t^l h_{\restriction
      t=0}\|_{\Lv} \right).
\end{equation}
Since $(\partial_t^N w,\partial_t^{N+1} w)_{\restriction t=0}$ is null, we
deduce from (\ref{ineg.w/h'}):
\begin{align}
\notag
\forall l=0,...,N,\quad\|\partial_t^{l+1} w_{\restriction t=0}\|^2_{\Lv(\Omega)}+\|\partial_t^{l} w_{\restriction t=0}\|^2_{\Hv_0^1(\Omega)} &\leq C\sum_{j=0}^{N-1} \|\partial_t^j \rot h_{\restriction t=0}\|_{\Lv(\Omega)} \\
\label{inegU0}
\|\partial_t^{l+1} w_{\restriction t=0}\|^2_{\Lv(\Omega)}+\|\partial_t^{l}
w_{\restriction t=0}\|^2_{\Hv_0^1(\Omega)} &\leq  C \sum_{j=0}^{N} \|\partial_t^j \rot h\|_{\Lv((0,T)\times\Omega)}.
\end{align}
(the second ligne is a consequence of the standard trace theorem with respect
to the time variable).
With the energy estimates of lemma (\ref{In En Lame}), applied to
(\ref{eqpolyU1}), we get, for any $0\leq l\leq N$: 
\begin{equation}
\notag
\|\partial_t^{l} w\|^2_{\Hv^1((0,T)\times\Omega)}\leq C\left(\|\partial_t^{l+1} w_{\restriction t=0}\|^2_{\Lv(\Omega)}+\|\partial_t^{l} w_{\restriction t=0}\|^2_{\Hv_0^1(\Omega)}+\|\rot \partial_t^l h\|^2_{\Lv((0,T)\times\Omega)} \right),
\end{equation}
which yields  exactly (\ref{inegU}).\par
On the other side, assumption (\ref{inegpoly}) implies 
\begin{gather*}
\|u_0\|^2_{\Hv^1_0(\Omega)}+\|u_1\|^2_{\Lv(\Omega)} \leq C\left(\|u_0\|^2_{\Lv(\Omega)}+\|u_1\|^2_{\Hv^{-1}(\Omega)} +\QQQ^N_T(u) \right).
\end{gather*}
Hence:
\begin{align}
\notag
\|v_0\|^2_{\Hv_0^1}+\|v_1\|^2_{\Lv} \leq& C\left(\|u_0\|^2_{\Hv_0^1}+ \|u_1\|^2_{\Lv}+\|w_0\|^2_{\Hv_0^1}+ \|w_1\|^2_{\Lv} \right)\\
\notag
\leq& C\left(\|u_0\|^2_{\Lv}+
  \|u_1\|^2_{\Hv^{-1}}+\QQQ^N_T(v)+\QQQ^N_T(w)+\|w_0\|_{\Hv_0^1}^2+\|w_1\|^2_{\Lv} \right)\\
\label{ineg.v0.v1}
\leq &C\left(\|v_0\|^2_{\Lv}+
  \|v_1\|^2_{\Hv^{-1}}+\QQQ^N_T(v)+\QQQ^N_T(w)+\|w_0\|_{\Hv_0^1}^2+\|w_1\|^2_{\Lv} \right)
\end{align}
With (\ref{inegU}) we get: 
\begin{align*}
\QQQ^N_T(w)\leq& C \sum_{l=0}^{N} \|\partial_t^{l+1}w\|^2_{\Lv((0,T)\times\Omega)}\\
\leq& C\sum_{j=0}^N \|\partial_t^j \rot h \|^2_{\Lv((0,T)\times \Omega)}.
\end{align*}
Taking into account (\ref{inegU0}) (with $l=0$) and the equation: $\rot(\partial_t
v\wedge \Bbf)=\partial_t h+\frac {1}{\beta} \rot \rot h$, which implies
(with lemma \ref{divrot}): 
$$\QQQ_T^N(v) \leq \sum_{j=0}^{N} \|\partial_t^j \rot
h\|^2_{\Lv((0,T)\times \Omega)},$$ 
we deduce from (\ref{ineg.v0.v1}):
\begin{equation*}
\|v_0\|^2_{\Hv^1_0(\Omega)}+\|v_1\|^2_{\Lv(\Omega)}\leq C \left\{\|v_0\|^2_{\Lv
(\Omega)}+\|v_1\|^2_{\Hv^{-1}(\Omega)}+\sum_{j=0}^N \|\partial_t^j \rot h \|^2_{\Lv((0,T)\times \Omega)}\right\}.
\end{equation*}
Now from lemma \ref{divrot}:
$ \|h_0\|^2_{\Lv} \leq C\left(\|\rot h\|^2_{\Lv((0,T)\times \Omega)} +\|\partial_t \rot h\|^2_{\Lv((0,T)\times \Omega)} \right)$.\\
This gives the following inequality on solutions of (\ref{magneto2}) with
initial data in $D(\AAA)$:
$$ \|V_0\|^2_X \leq C \left\{ \|V_0\|^2_{D(\AAA)'}+\sum_{j=0}^N \left( \|\partial_t^j V(0)\|^2_X-\|\partial_t^j V(T)\|^2_X\right) \right\}.$$
It is easy to check, with the criterum given by remark \ref{critere de
  fermabilite}, that the quadratic form:
$$ Q_T(V)=\sum_{l=0}^N \int_{\Omega} |\rot \partial_t^l h|^2 \,dy,$$
with domain $D(\AAA^N)$, is closable. All the assumptions of lemma
\ref{genpoly}, with $\PPP=\AAA$ hold which completes the proof of
proposition \ref{decpoly}.
\section{Defect measures}
\label{chap.mesures}
Let $N$ be an integer. For an open subset $U$ of an euclidian space, we
set:
$$ \Lv(U)\egaldef L^2(U,\CC^N),\quad \Hv^s(U)\egaldef H^s(U,\CC^N).$$
We consider an open subset $\Omega$ of $\RR^n$, $n\geq 1$, and a sequence
$(u^k)$ of functions on $\RR_t\times\Omega_y$ such that:
\begin{equation}
\label{CVfaible}
 u^k \underset{k\rightarrow +\infty}{\longCVf} 0 \dans \Hv^1_{\loc}(\RR\times\overline{\Omega}), 
\end{equation}
(in the sense that
$(\varphi u^k)$ converges weakly to $0$ in
$\Hv^1(\RR\times\overline{\Omega})$ for all $\varphi\in
C^{\infty}_0(\RR\times\overline{\Omega})$). 
We assume that every $u^k$ is 
solution of a wave equation in $\Omega$: 
\begin{equation}
\label{onde} (\nu^2\partial_t^2-\Delta) u^k=0, \text{ in } 
\RR\times\Omega.
\end{equation} 
We shall introduce in this section a measure describing, from a micro-local
point of view, the defect of compactness in $\Hv^1$ of the sequence $(u^k)$.
This description is of fundamental importance to show the observability inequalities of
the preceding section, for the Lam\'e system may be decomposed in two waves
equation (see paragraph
\ref{ondes.T.L}). Micro-local defect measures have been independently
introduced by
 P.~G\'erard and L.~Tatar \cite{PG91,Tar90}. We shall follow the construction
 of N.~Burq and G.~Lebeau, which describes the defect of convergence up to
 the boundary of $\Omega$.\par
We assume, for the sake of simplicity that the functions $u^k$ are smooth,
so that their traces on the boundary are always defined. In the sequel we shall always reduce to this case.\par
In subsection \ref{mes.notations} we will give a few definitions and
notations. In subsection \ref{par.mes.existence} we will state an existence
theorem of micro-local defect measures and set out their first properties.
Subsection \ref{par.propa} is devoted to the propagation theorem of the
measure (proved in \cite{BuLe99}), and subsection \ref{mes.traces} to some
important properties of the traces of $u^k$ on the boundary. Finally, in
section \ref{par.Lame}, we shall apply the construction of the measure to
the case of a sequence of solutions of the Lam\'e system.

\subsection{Notations}
\label{mes.notations}
\subsubsection{Local coordinates}
Consider an open cover of $\Omega$: $\Omega= \bigcup_{j=0}^J \Omega_j,$
where $\overline{\Omega}_0\subset \Omega$ and, for all $j\geq 1$,
$\Omega_j$ 
is a small neighbourhood of a point of $\partial \Omega$, such that on
$\Omega_j$, there are geodesic normal coordinates: 
$$ z\in \Omega_j \mapsto (y',x_n) \in Y\egaldef Y' \times ]0,l[,$$
where $x_n$ is the distance to the boundary, and $Y'$ an open subset of
$\RR^{n-1}$. Most objects introduced here are global objectsm but we will
mainly use local coordinates. For a large part of this section we 
choose one of the open set $\Omega_j$, $j\geq 1$.\par
Set $\XXX\egaldef\RR\times Y$ and denote the elements of $X$ by: 
$$ x=(\underbrace{x_0,x_1,..,x_{n-1}}_{x'},x_n),\quad x_0=t,\quad 
y=(\underbrace{x_1,x_2,..,x_{n-1}}_{y'},x_n),$$
Let:
$$ \RR^{n+1}_+\egaldef\left\{(x',x_n) \in \RR^{n+1},\; x_n>0\right\},\quad
\overline{\RR}^{n+1}_+\egaldef \overline{\RR^{n+1}_+},\quad\partial \XXX=\XXX'\times\{0\},\quad \Xbar\egaldef\XXX'\times [0,l[.$$
The set $\overline{X}$ is  an open subset of $\overline{\RR}^n_+$. Let $g$
be the natural metric on $Y$, induced by the change of coordinates. 
In a geodesic system of coordinates, $g$ is of the form:
\begin{equation*}
g(y)=\left[ \begin{array}{cc} g'(y) & 0\\ 0 & 1 \end{array} \right],\quad
\detg=\det g.
\end{equation*}
\subsubsection{Bundles on $X$}
Let's consider $T^*X=X\times\RR^{n+1}$ the cotangent bundle of 
$X$ and $S^*X$ the spherical cotangent bundle, which is defined to be the quotient
$$ S^*X\egaldef(T^*X\backslash\{|\xi|=0\} )\big/ \RR_+^*,$$
by the action of $\RR^*_+:$ $(\lambda,\xi) \mapsto \lambda \xi$. 
The elements of those two bundles will be denoted by:
$$ \rho=(x,\xi),\quad x\in X,\quad \xi=(\xi',\xi_n) \in 
\RR^{n}\times \RR,\quad \xi=(\tau,\eta).$$
There is a natural euclidian norm for the $\eta$-component of $T^*X$:
$ \|\eta\|^2\egaldef\transp{\eta} g^{-1} \eta$.\par 
We will also consider $T^*\partial X\egaldef\partial X_{x'}\times
\RR^{n}_{\xi'}$ the boundary cotangent bundle and $S^*\partial X$ the
associated spherical bundle.
\subsubsection{Operators in  the interior of $\Omega$}
\label{a.l.interieur}
Le $S_i^m$ the set of matrix symboles of degree $m$ with compact support in
$X$, which are the functions:
$$a(x,\xi)\in C^{\infty}(X\times \RR^{n+1},\Matr_N(\CC)),$$
whose $x$-projection is of compact support in $X$, satisfying the following
estimates:
\stepcounter{equation}
\newcounter{est.cpt}
\setcounter{est.cpt}{\value{equation}}
\begin{equation}
\label{est.symbole}
\tag{$\arabic{equation}_m$}
\forall \alpha,\; \forall 
\beta,\;\exists C_{\alpha\beta},\;\left|\partial_x^{\alpha}\partial_{\xi}^{\beta}a(x,\xi)\right|\leq 
C_{\alpha\beta}\left( 1+|\xi|\right)^{m-|\beta|},
\end{equation}
and which have a principal symbol $a_m(x,\xi)$, homogeneous function of
degree $m$ in $\xi$, such that $a-a_m$ satisfies
($\arabic{est.cpt}_{m-1}$) for large $|\xi|$.  
The operator of symbol $a$, $A=a(x,D)$, is defined by:  
$$ A v(x)\egaldef\frac{1}{(2\pi)^{n+1}}\int a(x,\xi) \hat{v}(\xi) e^{ix.\xi} 
d \xi.$$
In order to act on functions which are only defined in $X$, it is
convenient to consider only the set $\AAA^m_i$ consisting of operators
$A$ which are of compact support in $X$, in the sense that $A=\varphi A
\varphi$ for a function $\varphi \in C^{\infty}_0(X)$. An operator in
$\AAA^m_i$ maps a distribution in $X$ to a 
compactly supported distribution in $X$. We shall denote by $\sigma_m(A)$
the principal symbol of an operator $A$ of degree $m$.
\subsubsection{Operators near the boundary}
\label{au.bord}
Let $S_b^m$ be the set of matricial \textbf{tangential} symbols of degree
$m$ with compact support in $X$, defined as the functions:
$$a(x,\xi')\in C^{\infty}(\Xbar\times \RR^{n},\Matr_N(\CC)),$$
whose $x$-projection has compact support in $\Xbar$, satisfying the estimations:
\stepcounter{equation}
\newcounter{est.cpt.b}
\setcounter{est.cpt.b}{\value{equation}}
\begin{equation}
\label{est.symbole.b}
\tag{$\arabic{equation}_m$}
\forall \alpha,\; \forall 
\beta,\;\exists C_{\alpha\beta},\;\left|\partial_x^{\alpha}\partial_{\xi'}^{\beta}a(x,\xi')\right|
\leq C_{\alpha\beta}\left( 1+|\xi'|\right)^{m-|\beta|},
\end{equation}
and which have a principal symbol $a_m(x,\xi')$, homogeneous of degree $m$
in $\xi'$ and such that $a-a_m$ satisfies the inegalities 
($\arabic{est.cpt.b}_{m-1}$) for large $|\xi'|$.
We define the operator of symbol $a$, $A=a(x,D')$, by:  
$$ A v(x)\egaldef\frac{1}{(2\pi)^n}\int a(x,\xi') \hat{v}(\xi',x_n) 
e^{ix'.\xi'} d 
\xi'.$$
Here, the Fourier transform of $v$ is only taken with respect to the
tangential variable $x'$. As in the interior case, we introduce the set
$\AAA^m_b$ of tangential operators $A$ with compact support in $\Xbar$, i.e such that
$A=\varphi A \varphi$ for a compactly supported function $\varphi\in C^{\infty}_0(\Xbar)$.\par
The set of all pseudo-differential operators of interest for us will be
denoted by:
$$ \AAA^m\egaldef \left\{ a=A_i+A_b,\; A_i \in \AAA^m_i,\;A_b \in \AAA^m_n\right\}. $$

\subsubsection{Sobolev spaces}
Let $s\in \RR$ and $\omega$ be an open set of $\RR^{n}$. As 
mentionned before, we denote by
$\Hv^s(\omega)$ the Sobolev space of vector-valued distributions (which may be defined as the set of
restrictions to $\omega$ of elements of $\Hv^s(\RR^n)$, endowed with the
quotient norm). We also consider the space $\Hv^s_{\loc}(\omega)$, the
space of vector-valued distributions such that: 
$$ \forall \varphi \in C^{\infty}_0(\omega), \quad \varphi u \in \Hv^s,$$
and $\Hv^s_{\comp}(\omega)$, the space of distributions in $\Hv^s(\omega)$
  compactly supported in $\omega$. The notation 
  $\Hv^s_{\loc}(Z)$,  will also be used when $Z$ is not open
  ($Z=\RR\times\overline{\Omega}$, or $Z=\overline{X}$, in the following
  natural sense:

$$ u^k \tendk u\; (\text{ or }=O(1)) \dans \Hv^s_{\loc}(Z)\iff \forall \varphi \in
C^{\infty}_0(Z), \quad \varphi u^k \tendk \varphi u\; (=O(1)) \dans \Hv^s,$$
where $C^{\infty}_0(Z)$ is the space of $C^{\infty}$, compactly supported
  functions in $Z$. We will also consider the following spaces,
  suitable for boundary-value problems: 
\begin{align*}
\Hv^{0,s}_{\loc}(\overline{X})&=L^2(0,l;\Hv^s_{\loc}(X'))\\
\Hv^{0,s}_{\comp}(\overline{X})&=\left \{u \in \Hv^{0,s}_{\loc}(\Xbar),\;
  \exists \varphi \in C^{\infty}_0(\Xbar),\; u=\varphi u\right\}.
\end{align*}
Note that the elements of $\AAA^m_i$
are continuous maps:
$$ \Hv^s_{\loc}(X)\longrightarrow 
\Hv^{s-m}_{\comp}(X),$$
and those of $\AAA^m_b$ are continuous maps:
$$ \Hv^{0,s}_{\loc}\left(\Xbar\right)\longrightarrow 
\Hv^{0,s-m}_{\comp}\left(\Xbar\right).$$
It is possible to ``micro-localize'' convergence properties in $\Hv^s$ and $\Hv^{0,s}$:
\begin{Def}
\label{defCV}
Let $\rho\in S^*X$. The sequence $(v^k)$ is said to be {\bf bounded
  (respectively converging to $0$) in $\Hv^s_{\rho}$ } when there exists
$A\in\AAA^s_i$, whose principal symbol is invertible near $\rho$ and such
that $(A v^k)$ is bounded in $\Hv^s$ (respectively converges to $0$ in
$\Hv^s$.\par
Let $\rho'\in S^*\partial X$. The sequence $(v^k)$ is said to be {\bf bounded
  (respectively converging to $0$) in $\Hv^{0,s}_{\rho}$ } when there exists
$A\in\AAA^s_i$, whose principal symbol is invertible near $\rho'$ and such
that $(A v^k)$ is bounded in $\Hv^{0,s}$ (respectively converges to $0$ in
$\Hv^{0,s}$).
\end{Def}
Note that, according to proposition
\ref{prop1.appendice} of the appendix, for a sequence of solutions of
(\ref{onde}), the convergence in $\Hv^{0,1}$ and $\Hv^1$ are
equivalent. The spaces $\Hv^{0,1}$ and the tangential operators are thus well
fitted for the description of the $\Hv^1$ convergence of $(u^k)$.
\subsubsection{Melrose's compressed cotangent bundle}
We shall now introduce a bundle which naturally contains as subbundles both
bundles $T^*X$ and $T^*\partial X$. For this purpose, set $T^*_{b}
X\egaldef \Xbar\times \RR^{n+1}$, endowed with its canonical topology and consider: 
\begin{eqnarray*}
T^*\Xbar &\overset{j}{\longrightarrow}& T^*_b X \\
(x,\xi',\xi_n) &\longmapsto & (x,\xi',r=x_n\xi_n).
\end{eqnarray*}   
The mapping $j$ restricts to a continuous map:
$$T^*X \longrightarrow T^*_b X\cap \{x_n > 0\},$$
 which identifies $T^*X$ to a subbundle of
dimension 
$2(n+1)$ of the interior of $T^*_b X$. Furthermore, the restriction of 
$j$ to $x_n=0$  defines a map from 
$T^*\Xbar\cap\{x_n=0\}$ to
$T^*_b X\cap\{x_n=0\}$, whose kernel is the set $\{\xi'=0\}$. This clearly identifies:
$$ T^*\partial X \approx (T^*\Xbar\cap\{x_n=0\}) / \RR_{\xi_n},$$
(quotient taken by identifiying all
the points  $(\tx',\txi',\xi_n)$, $\xi_n
\in \RR$) 
with a $2n$-dimensional subbundle of $T^*_b X$. The set of all sections of
$T^*_b X$, with the above identifications, may be seen as the dual bundle
of the bundle of all vector fields on $X$ tangent to $\partial X$. It is
called the compressed cotangent bundle.\par
We will also consider $S^*_b X$ the spherical bundle of $T^*_b X$, 
which naturally contains the spherical bundles $S^*X$ and $S^*\partial X$.

\subsubsection{Symbol of $P$ and related manifolds}
\label{var.car}
The equation (\ref{onde}) takes the following form in local coordinates:
\begin{equation}
\label{onde'}
P u^k=0,\quad
 P\egaldef-\detg^{-1/2}\partial_{x_n}\detg^{1/2} \partial_{x_n} 
+Q,
\end{equation}
where $Q$ is a scalar tangential differential operator of degree $2$. Let
$q(x,\xi')$ be the principal symbol of $Q$, and  
$p(x,\xi)=\xi_n^2+q(x,\xi')$ the principal symbol of $P$. They are both
scalar, homogeneous polynomials of degree $2$ with respect to $\xi$. 
Let:
$$ \CarP\egaldef\left\{(x,\xi)\in T^*\Xbar,\; p(x,\xi)=0 \right\}, \quad 
Z\egaldef j\left(\CarP\right),\quad 
\widehat{Z}\egaldef j\left(\CarP\right)\cup j\left(\{x_n=0\}\right),$$ 
and $\SCarP$, $SZ$, $S\widehat Z$ the corresponding spherical bundles.\par
Decompose $T^*(\partial X)$ 
(and $S^*(\partial X$)) into the disjoint union of the elliptic region $\EEE$, the
glancing region $\GGG$ and the hyperbolic region $\HHH$:
$$ \EEE\egaldef\{q_0>0\},\quad \GGG\egaldef\{q_0=0\},\quad \HHH\egaldef\{q_0<0\}.$$
\subsubsection{Global measure}
\label{par.recol}
The defect measure is at first constructed in each of the preceding local
coordinate systems. The objects obtained are then pieced together to
$M=\RR\times\Omega$. It is easy to define from local objects global Sobolev
spaces and bundles on $M$, such as Melrose's compressed cotangent bundle
$T^*_b M$. We shall use the same notations ($\CarP$, $Z$, 
$\widehat{Z}$, $\SCarP$, $SZ$, $S\widehat{Z}$,$\dots$) for the local and global objects.
The definition of global operators is less natural in our setting. 
The symbol $\AAA^m$ will denote the set of operators $A$ acting on
functions on $M$, which are of the form:
$$ A=\sum_{j=0}^J A_{(j)}.$$
where $A_{(0)}$ is a classical pseudo-differential operator of order $m$
with compact support in $M$ and each $A_{(j)}$ is an
operator of the sets $\AAA^m$ defined in each system of local
coordinates. The global space $\AAA^m$ depends of the coordinate
patches chosen, which shall not cause any problem in the remaining of the
article. For a totally intrisic construction, we could have used Melrose's
totally characteristic operators (see \cite[chap 18.3]{Hor}).  
\subsection{Existence of the measure}
\label{par.mes.existence}
\subsubsection{The existence theorem}
\label{mes.existence}
The next elementary proposition shows that for any $A\in \AAA^0$, the
behaviour of $Au^k$ in $\Hv^1$ only depends upon the restriction of its
principal symbol to $S\widehat{Z}$: 
\begin{prop}
\label{en.dehors.de.SZ}
Let $A_b\in \AAA_b^{-\eps}$. Then:
$$ A_bu^k \underset{k\rightarrow 0}\longrightarrow 0 \text{ in } 
\Hv^1.$$
Let $A_i\in \AAA_i^0$, whose principal symbol vanishes on 
$\CarP$. Then:
$$ A_i u^k \underset{k\rightarrow 0}\longrightarrow 0 \text{ in } 
\Hv^1. $$
\end{prop} 
According to proposition \ref{en.dehors.de.SZ}, it is sufficient to describe the
$\Hv^1$ convergence of $(u^k)$ near $S\widehat{Z}$, in the sense given by
definition \ref{defCV}. Let $\MMM$ be the set of matice-valued measures on
$S\widehat Z$, i.e. the dual space of:  
 $$\CCC\egaldef C^0_0\left(S\widehat{Z},M^N\left(\CC\right)\right),$$
and $\MMM^+$ the subset of all positive measures in  
$\MMM$, i.e. measures $\mu$ which satisfy:
\begin{equation*}
\forall z\in S\widehat Z,\quad b(z)\geq 0\implies \big< \mu,b \big> 
\geq 0.
\end{equation*} 
($M\geq 0$ means $M$ is positive hermitian).\par
Before coming to the main theorem of this paragraph, we shall introduce a
technical condition on $u^k$:
\begin{Def}
Let the sequence $(u^k)$ satisfies (\ref{CVfaible}) and (\ref{onde}). We
shall say that $(u^k)$ is {\bf regular on the boundary} when one the
following equivalent assumptions is satisfied:
\stepcounter{equation}
\newcounter{eq.traces.cpt}
\setcounter{eq.traces.cpt}{\value{equation}}
\begin{gather}
\tag{$\arabic{eq.traces.cpt}a$} 
u^k_{\rx} =o(1) \text{ in } \Hv^{1/2}_{\loc}(\partial X),\;k \longrightarrow +\infty\\
\tag{$\arabic{eq.traces.cpt}b$}
\partial_{x_n} u^k_{\rx} =o(1) \text{ in } 
\Hv^{-1/2}_{\loc}(\partial X),\;k\longrightarrow +\infty.
\end{gather} 
\end{Def}
Note that this is a very weak condition: the standard trace theorems imply
conditions (\arabic{eq.traces.cpt}) with $O(1)$ instead of $o(1)$. All the
sequences $(u^k)$ in this work shall satisfy this condition. For the proof
of the equivalence between (\arabic{eq.traces.cpt}a) and
(\arabic{eq.traces.cpt}b) see \cite[lemma 2.6]{BuLe99}.
\begin{theo}
\label{existence.mu}
Let $u^k$ be such that (\ref{CVfaible}), (\ref{onde}) and
(\arabic{eq.traces.cpt}) hold. Then there exists a subsequence of $(u^k)$,
still denoted by $(u^k)$, and a measure $\mu \in \MMM^+$, called
{\bf micro-local defect measure}, such that
$\mu(\EEE\cup\HHH=0)$ and: 
\begin{equation}
\label{mu.lim.2}
\forall A_j \in \AAA^j,\; j\in \{1,2\},\; \lim_{k\rightarrow +\infty} \left( A_2u^k+A_1 D_{x_n} 
u^k,u^k\right)=\big<\mu,\frac{a_2+\xi_n a_1}{\tau^2} \big>.
\end{equation}
\end{theo} 
In (\ref{mu.lim.2}), the notation $(.,.)$ stands for the $\Lv$-scalar
product on $M$ (in local coordinates, it is the scalar product on $X$ using the metric
$\detg^{1/2} dy \,dt$) and $\mu$ is considered as a measure on the subset
$\SCarP$ of $S^*\overline{X}$, using the canonical map $j$, which is an homeomorphism:
$$  \SCarP\overset{j}{\longrightarrow} SZ\backslash \HHH.$$ 
This is made possible by the fact $\mu(\EEE\cup \HHH)=0$. In the case where
$(u^k)$ is not regular on the boundary, it is still possible to define $\mu$,
but $\mu(\EEE)$ is non-null, which makes the statement of condition 
(\ref{mu.lim.2}) more intricate.
\begin{rem}
The measure $\tilde{\mu}=\mu\indic_{\{x_n>0\}}$ may be seen as the standard micro-local defect measure (cf \cite{PG91}) of the bounded
sequence $(u^k)$ of $\Hv^1_{\loc}(M)$. This interior measure describe the
compactness defect of $(u^k)$ in $\Hv^1_{\loc}(M)$ (in particular, it is
null when $(u^k)$ converges to $0$ in this space), but not in
$\Hv^1_{\loc}(\overline{M})$: $\tilde{\mu}$ vanishes when $(u^k)$ concentrates on
$\partial M$, even if it does not converge to $0$ in $\Hv^1_{\loc}(\overline{M}
)$. On the other hand:
$$ \varphi \in C^{\infty}_0(\overline{M})\Longrightarrow \int \varphi \left| \nabla_y u^k \right|^2 dx+\int \varphi \left| \partial_t u^k \right|^2 dx \tendk \left<\mu,\varphi\right>.$$
Thus, $\mu$ sees all the $\Hv^1_{\loc}(\overline{M})$ density of $(u^k)$ at
infinity. More precisely, it gives a micro-local description of this density: 
$$ \rho \in \supp \mu \iff u^k \underset{k\rightarrow 
+\infty}\longrightarrow +0 \text{ in } \Hv^1_{\rho}. $$ 
When  $\rho$ is a boundary point, one should replace $\Hv_{\rho}^1$ 
by $\Hv_{\rho}^{0,1}$.
\end{rem}
Theorem \ref{existence.mu} is a new formulation, using lemma 2.7 of 
\cite{BuLe99}, of proposition 2.5 of this article.\par
\subsubsection{A sufficient condition of nullity for $\mu$}
\label{CS.mu.0}
Let $\trho\in S\widehat{Z}$ an interior point and  $A\in 
\AAA^2_i$,  whose principal symbol is invertible at $\trho$. By elementary
symbolic calculus on classical operators, it is easy to show, with formula (\ref{mu.lim.2}):
\begin{equation}
\label{mu.0.int}
 A u^k\tendk 0 \dans \Hv^{-1}_{\trho} \Rightarrow \trho \notin 
\supp \mu.
\end{equation}
The same statement holds in $\GGG$:
\begin{prop}
\label{prop.mu.G.0}
Consider an operator of the form:
$$ A=A_0 D_{x_n}^2+A_{1}D_{x_n}+A_2,\; A_j \in \AAA^j, 
\;a_j:=\sigma(A_j),$$
such that:
\begin{equation}
\label{Auk.CV0}
A u^k \tendk 0 \dans \Hv_{\trho}^{0,-1}. 
\end{equation}
Assume that $(u^k)$ is regular on the boundary and let $\trho\in \GGG$ such
that $a_2(\trho)$ is invertible. Then:
$$ \trho \notin \supp \mu.$$
\end{prop}
\begin{rem}
Proposition \ref{prop.mu.G.0} is trivial when $A=A_2\in \AAA^2$ (it is
essentially the definition of $\Hv_{\trho}^{0,-1}$). 
\end{rem}
\begin{rem}
Near an hyperbolic point, it is more difficult to state proposition
\ref{prop.mu.G.0}. Indeed it is much more relevant to study $\mu$, in the set
$\{x_n\geq \eps_0>0\}$, near rays in and out of $\trho$.   
(see paragraph \ref{par.hyp}). 
\end{rem}
\begin{rem}
Note that according to the appendix, the convergence to $0$ of $(Au^k)$ in
the space $\Hv^{-1}(X)$ near $\tx$ would imply
(\ref{Auk.CV0}). Furthermore, the proof of the lemma will show that
assumption (\ref{Auk.CV0}) is equivalent to:
$$ (J A u^k,u^k) \tendk 0,\; \forall J\in \AAA^0,\text{ with support close
  enough to }\trho.$$
\end{rem}
\begin{proof}
Let $j=\sigma(J)$. According to (\ref{mu.0.int}), $\mu\indic_{\{x_n>0\}}$
is null, near $\trho$. The same property remains to be proved on 
$\mu\indic_{\{x_n=0\}}$. Let:
$$J\in \AAA^0,\quad \psi \in C^{\infty}_0(\RR)\text{ tel que } 
\psi(0)=1,\quad J_{\eps}\egaldef \psi\left(\frac{x_n}{\eps}\right) J.
$$ 
In view of (\ref{Auk.CV0}) and formula (\ref{mu.lim.2}): 
\begin{equation*}
\big< \mu,\psi\left(\frac{x_n}{\eps}\right)j \frac{a_0 
\xi_n^2+a_1\xi_n +a_2}{\tau^2}\big> =
\lim_{k\rightarrow +\infty} (J_{\eps}A u^k,u^k)
= 0.
\end{equation*}
Letting $\eps$ goes to $0$, the dominated convergence theorem and the fact
that $\xi_n$ is null on the support of $\mu\indic_{\{x_n=0\}}$ give:
\begin{equation}
\label{muG.0.a}
\big< \mu,\indic_{\GGG} \frac{j a_2}{\tau^2} \big>=0.
\end{equation}
Let $\psi\in S^0_b$ be scalar, positive, and compactly supported near
$\trho$ such that $a_2$ is invertible on the support of $\psi$, and choose
$J$ such that:
$$j(x,\xi')=\psi(x,\xi')a_2^{-1} \tau^2.$$
The equality (\ref{muG.0.a}) then shows that  
$\big<\mu,\indic_{\GGG}\psi\big>=0$,
which completes the proof using the positivity of $\mu$.
\end{proof}

\subsection{The propagation theorem}
\label{par.propa}
\subsubsection{The generalized bicharacteristic flow}
The characteristic curves of the hamiltonian flow of $p$:
$$ H_p=\partial_{\xi}p\partial_{x}-\partial_{x}p\partial_{\xi}$$
define a local flow on $T^*X$. The symbol $p$ is homogeneous of degree $2$
in $\xi$, so that the flow of $H_p$ does not yield a flow on the quotient
space $S^*_b X$. To get such a flow, we shall replace $p$ by $p/\tau$ which
is homogeneous of degee $1$. Note that on the support of $\mu$ (where
$\tau$ does not vanish), $p$ is null, so that  $\frac{1}{\tau}H_p$ and
$H_{p/\tau}$ are equal. Furtermore, the integral curves of
$\frac{1}{\tau}H_p$ and $H_p$ are the same.\par
Let $\Sigma$ be a small conic open subset of $Z=j(\CarP)$.
Set $q_0\egaldef q_{\rx}$, $q_1 \egaldef \partial_{x_n}q_{\rx}$ and
\begin{gather*}
\Sigma^0\egaldef\Sigma\cap\{x_n>0\}\\
\Sigma^1\egaldef\HHH=\Sigma\cap\{x_n=0,\;q_0<0\}\\
\Sigma^2\egaldef\Sigma\cap\{x_n=0,\;q_0=0\;q_1\neq 0\}\\
\Sigma^{k+3}\egaldef\Sigma\cap\{x_n=0,\; q_0=q_1=...=H_{q_0}^k q_1=0, \; 
H_{q_0}^{k+1}\,q_1\neq 0\}.
\end{gather*}
Assume that in $\Sigma$, there is no contact of infinite order between the
bicaracteristic curves of $P$ and the boundary, which means that for a
certain finite integer $J$:
\begin{equation}
\label{contact}
\exists J\in \NN,\quad \Sigma=\bigcup_{j \leq J} \Sigma^J. 
\end{equation}
Decompose $\Sigma^2$ in the disjoint union:
$$ \Sigma^2=\GGG^{2,+}\cup\GGG^{2,-}, \quad
\GGG^{2,+}\egaldef \Sigma^2\cap\{q_1< 0\},\; \GGG^{2,-}\egaldef \Sigma^2\cap\{q_1> 0\}.$$
The set $\GGG^{2,+}$ is the set of {\bf strictly diffractive} points  
and $\GGG^{2,-}$ the set of {\bf strictly gliding} points.
\begin{Def}
Let $\gamma$ be a map form a real interval $I$ to 
$\Sigma$ and:
$$\Gamma(s)=j^{-1}(\gamma(s)) \in \SCarP$$
which is defined as long as $\gamma(s) \notin \HHH$. Such a map
$\gamma(s)=(x(s),\xi(s))$ is called a {\bf ray}, or a 
{\bf general bicharacteristic curve} when 
$\gamma$ is continuous
from $I$ to $\Sigma$ and for all $s_0$ in $I$:
\begin{itemize}
\item if $x_n(s_0)>0$, $\Gamma$ is differentiable in $s_0$ and:
$$ \Gamma'(s_0)=\frac 1\tau H_p \Gamma(s_0); $$
\item if $\gamma(s_0) \in \HHH \cup \GGG^{2,+}$,
$$ \exists \eps>0, \quad \forall s \in 
]s_0-\eps,s_0[\cup]s_0,s_0+\eps[,
\quad x_n(s)>0;$$
\item if $\gamma(s_0) \in \GGG \backslash \GGG^{2,+}$, $\Gamma$ is well
  defined and differentiable near $s_0$ and:
$$ \Gamma'(s_0)=\frac{1}{\tau} H_{q_0} \Gamma(s_0). $$
(Thus, if $\gamma$ stays in this region, its spatial projection is a
geodesic of the boundary.) 
\end{itemize}
\end{Def}
Under the assumption (\ref{contact}), R.~Melrose and J.~Sj\"ostrand have
shown that for any $\rho\in \Sigma$, there exists an unique maximal ray
$\gamma$ taking values in $\Sigma$ such that $\gamma(0)=\rho$ (cf
\cite{MeSj78}, \cite[chap 24.3]{Hor}). In the sequel, we shall denote by
$\phi(s,\rho)$ the resulting flow (satisfying $\phi(0,\rho)=\rho$). The
function $p/\tau$ being homogeneous of degree $1$ in $\xi$, the flow
$\phi$ passes to the quotient and defines a flow on $\Sigma/\RR^*_+$. 
\subsubsection{The uniform Lopatinsky conditions}
\label{par.Lopa}
\begin{notations}
Let $S^m_{\partial}$ be the set of symbols $a(x',\xi')$ 
of
pseudo-differential operators on $\partial X$, with compact support in 
$x'$, with principal symbol homogenous of degree m in 
$\xi'$, 
and $\AAA_{\partial}^m$ the set of corresponding compactly supported
pseudo-differential operators (cf paragraph \ref{a.l.interieur} for precise
definitions). The Sobolev spaces on $\partial X$, defined as those on $X$,
shall be denoted by $\Hv^s_{\partial}$, $\Hv^s_{\loc,\partial}$,
$\Hv^s_{\trho,\partial}$.
\end{notations}
An approximate pseudo-differential equation on the traces of $u^k$ is said
to satisfy Lopatinsky conditions when it is independent of the equation
$Pu^k=0$. Precisely: 
\begin{Def}
\label{def.Lopa}
Under the assumptions
(\ref{CVfaible}) and (\ref{onde}), the sequence $(u^k)$ is said to 
satisfy {\bf uniform Lopatinsky} boundary conditions near $\trho\in
S^*\partial X$ when: 
\begin{itemize}
\item if $\trho \in \GGG$, $\exists B_{-1} \in \AAA^{-1}_{\partial}$ such that:
\begin{equation}
\label{LopaG}
\left\{ 
\begin{gathered}
u^k_{\rx}=B_{-1} \left( D_{x_n} u^k_{\rx}\right) +h^k\\
h^k \tendk 0 \text{ in } \Hv^1_{\trho,\partial};
\end{gathered}
\right.
\end{equation} 
\item if $\trho \in \HHH$, $\exists B_0 \in \AAA^0_{\partial}$ such that:
\begin{gather}
\tag{\ref{LopaG}'}
\label{LopaH}
\left\{ 
\begin{gathered}
D_{x_n} u^k_{\rx}-\Lambda u^k_{\rx}=B_{0} \left( D_{x_n} 
u^k_{\rx}+\Lambda u^k_{\rx}\right) +h^k\\
\sigma(B_0) \text{ invertible near } \trho\\
h^k \tendk 0 \text{ in } \Lv_{\trho,\partial}
\end{gathered}
\right.\\
\notag
\Lambda \in \AAA_{\partial}^1,\quad 
\sigma(\Lambda)=\sqrt{q_0(x',\xi')}=\sqrt{\nu^2\tau^2-\|\eta'\|^2}\text{
  near }\trho. 
\end{gather} 
( $\|\eta'\|^2=\transp{\eta'} g'{}^{-1} \eta'$ is the natural euclidian
norm in the local coordinate system). 
\end{itemize}
\end{Def}
\begin{exemples}
\begin{itemize}
\item The Dirichlet boundary condition, $u^k_{\rx}=0$, or more generally a
  pseudo-diffe\-ren\-tial boundary condition of the form:
\begin{equation}
\label{bord.iR}
\begin{gathered}
u^k_{\rx}=B_{-1} D_{x_n}u^k_{\rx} +h^k,\quad B_{-1} \in \AAA^{-1}\\ 
 h^k \tendk 0 \text{ in } \Hv^1_{\trho,\partial},
\end{gathered}
\end{equation}
where the eigenvalues of $\sigma(B_{-1})$ are all pure imaginary numbers
near $\trho$, is an uniform Lopatinsky boundary condition, whether $\trho$ is
glancing or hyperbolic. In the glancing case (\ref{bord.iR}) corresponds
exactly to the definition (\ref{LopaG}) and in the hyperbolic case, both
operators $\Id-\Lambda B_{-1}$ and $\Id+\Lambda B_{-1}$ are elliptic in
$\trho$, and it is easy to show (\ref{LopaH}), taking $B_0=(\Id-\Lambda
B_{-1})(\Id+\Lambda B_{-1})^{-1}$ (where $(...)^{-1}$ stands for a
parametrix near $\trho$).
\item Neumann condition:
$$ D_{x_n} u^k_{\rx} \tendk 0 \text{ in } \Lv_{\trho,\partial}$$
is an uniform Lopatinsky condition near any  
$\trho\notin \GGG$.
\item Boundary conditions which are not Lopatinsky conditions in the
  hyperbolic region are described in paragraph \ref{par.hyp}.
\end{itemize}
\end{exemples}
In the glancing case, a boundary condition of the form (\ref{LopaG}) locally
implies better estimates than the standard ones on the traces of $u^k$,
which we shall state in proposition \ref{Lopa.traces}. This shows in
particular that $(u^k)$ is regular near the boundary, and that
$\mu(\EEE)=0$. As a consequence, the bicharacteristic flow is defined
$\mu$-almost everywhere.\par
The set $\GGG^{2,+}$ is transverse to the bicharacteristic flow. The next
result is necessary to the propagation of $\mu$ by the flow, which is
treated in the next paragraph. 
\begin{lem}
\label{glissant.nul}
Suppose that on $\Sigma$, $u^k$ satisfies uniform Lopatinsky boundary
conditions. Then: 
$$ \mu\left(\GGG^{2,+}\cap \Sigma\right)=0.$$
\end{lem}
\subsubsection{The propagation theorem}
When an uniform Lopatinsky condition holds, $\mu$ propagates along the
integral curves of the bicharacteristic flow. In the hyperbolic region,
there is a jump (which depends upon the boundary condition). We shall only
state a propagation theorem for the support of $\mu$, without giving a
complete description of the propagation of $\mu$.
\begin{theo} 
\label{thm.propa}
Let $\trho\in \HHH\cap \GGG$ such that $(u^k)$ satisfies Lopatinsky
boundary conditions. Consider a small conic open neighbourhood $\Sigma$ of
$\trho$ in $S\hat{Z}$ such that on $\Sigma$, (\ref{LopaG}) (or
(\ref{LopaH})) holds. Then the support of $\mu$ is, in $\Sigma$ invariant by
the bicharacteristic flow.  
\end{theo}
(cf \cite[chap. 3.3, th.1]{BuLe99})\par
In other terms, if $\rho\in \Sigma$ is on the support of $\mu$, so is the
entire bicharacteristic passing through $\rho$ in $\Sigma$.
\begin{rem}
Inside $M$, theorem \ref{thm.propa} is an easy consequence of the
transport equation on $\mu$:
\begin{equation}
\label{eq.mu.int}
\big< \mu, \left\{p/\tau,a\right\} \big>,\quad a\in 
C^{\infty}_0\left(Z\cap\{x_n>0\}\right),
\end{equation}
which may be immediately derived, using symbolic calculus, from the
elementary property:
\begin{equation}
\label{commu.int}
\lim_{k\rightarrow +\infty} \left(A_1 P u^k-P A_1 u^k,u^k\right) 
=0,\quad 
A_1\in \AAA_i^1,
\end{equation}
obtained by integration by parts with the equation (\ref{onde'}).
Near a boundary point, property (\ref{commu.int}), with $A_1\in \AAA_b^1$,
still holds with an additional boundary term. Consequently,
(\ref{eq.mu.int}) holds only for a certain class of function $a\in
C^{\infty}_0(SZ)$, satisfying a particular boundary condition on
$\{x_n=0\}$ (condition chosen to kill the boundary terms when $k$ tends to
$\infty$). The proof of the propagation theorem, which is fairly technical,
uses (\ref{eq.mu.int}), and near strictly diffractive points, lemma
\ref{glissant.nul}. The boundary condition on $a$ gives the exact value of
the jump in the hyperbolic region. See \cite[par. 3]{BuLe99} for details.
\end{rem}
\subsection{Estimates on traces}
\label{mes.traces}
We now state precise properties of the traces of $u^k$ in the hyperbolic,
elliptic and glancing regions, which are one of the main tools of the
proofs of the following sections. Those results are fairly classical, and we
only shall give a proof (in the appendix) for the glancing case. See
\cite[appendix]{BuLe99} for proofs in the hyperbolic and elliptic cases.  
In this paragraph, we shall always assume $(u^k)$ satisfies (\ref{CVfaible})
and (\ref{onde}).
\subsubsection{Hyperbolic region}
\label{par.hyp}
Near an hyperbolic point, one gains without any boundary condition, half a
derivative in comparision with the standard traces theorem.  
\begin{prop}
\label{traces.hyp}
Let $\trho\in \HHH$. Then:
\begin{align*}
u^k_{\rx}&\underset{k\rightarrow \infty}{=}O(1) \dans 
\Hv^1_{\trho,\partial}\\
\partial_{x_n} u^k_{\rx}&\underset{k\rightarrow \infty}{=}O(1) \dans 
\Lv_{\trho,\partial}.
\end{align*}
\end{prop}
In view of the propagation theorem in the interior of $M$, the support of $\mu$ is, near $\trho$, the union of
incoming rays (integral curves of $H_{p/\tau}$ along which $\xi_n<0$) and
outgoing rays (integral curves of $H_{p/\tau}$ along which
$\xi_n>0$). When the sequence satisfies uniform Lopatinsky conditions,
theorem \ref{thm.propa} is equivalent to the fact that if an incoming
(respectively outgoing) ray is in the support of $\mu$, so is the outgoing
(respectively incoming) ray passing through the same hyperbolic point. In the
opposite case where the support of $\mu$ contains, locally, only incoming
(or only outgoing) rays, one gets a boundary condition which is in a
certain sense ``orthogonal'' to Lopatinsky uniform conditions: 
\begin{prop}
\label{traces.hyp.0}
Assume that near $\trho\in \HHH$, on the support of $\mu$, $\xi_n>0$. Then:
\begin{gather*}
D_{x_n} u^k_{\rx}+ \Lambda u^k_{\rx}=o(1) \dans \Lv_{\trho,\partial}\\
\notag
\Lambda \in \AAA^1_{\partial},\quad \sigma_1(\Lambda)=\sqrt{\nu^2 
\tau^2-\| \eta'\|^2}
\end{gather*}
On the other hand, if near $\trho\in \HHH$, on the support of $\mu$, 
$\xi_n<0$, then:
\begin{equation*}
D_{x_n} u^k_{\rx}-\Lambda u^k_{\rx}=o(1) \dans \Lv_{\trho,\partial}.
\end{equation*}
In particular, if $\mu$ is null near $\trho$, 
$$ D_{x_n} u^k_{\rx} \tendk 0 \dans \Lv_{\trho,\partial},\quad 
u^k_{\rx} \tendk 0 \dans \Hv^1_{\trho,\partial}.$$
\end{prop} 
\subsubsection{In the elliptic region}
In $\EEE$, the equation (\ref{onde'}) implies a pseudo-differential
traces equation on $u^k$:
\begin{prop}
\label{traces.ell}
Let $\trho\in \EEE$ and $(u^k)$ satisfy the assumptions of 
theorem \ref{existence.mu}. 
Let $M>0$. Then:
\begin{gather}
\label{eq.elliptique}
D_{x_n} u^k_{\rx} +\Xi u^k_{\rx} \underset{k\rightarrow 
+\infty}\longrightarrow 0 \text{ in }H^M_{\trho}\\
\notag
\Xi \in \AAA^1,\quad \sigma_1(\Xi)=i\sqrt{q_0} 
=i\sqrt{\|\eta'\|^2-\nu^2\tau^2} \text{ near } \trho.
\end{gather}
\end{prop}
In particular, if a boundary condition independent of (\ref{eq.elliptique})
holds on $u^k$ near $\rho_0$ (such a condtion is called as in the glancing
and hyperbolic cases an uniform Lopatinsky condition), the traces of $u^k$
converge to $0$ in appropriate Sobolev spaces $H^M_{\trho}$. Proposition
\ref{traces.ell} still holds in a much more general case, for example if
$P$ is replaced by a non-scalar operator $\mathbf{P}$. The principal symbol
of $\Xi$ depends again upon the principal symbol of $\mathbf{P}$. In the next
proposition, we only state a consequence of this fact when $(u^k)$
satisfies Dirichlet boundary conditions (which are of uniform Lopatinsky type).\begin{prop}
\label{traces.ell.M}
Let:
$$ \mathbf{P}\egaldef D_{x_n}^2 + \mathbf{Q}_1D_{x_n}+\mathbf{Q}_2,$$
where each $\mathbf{Q_j}$ is a matricial pseudo-differential operator
of degree $j$, with principal symbols $\mathbf{q}_j$. Let:
$$ \mathbf{p}(x,\xi)\egaldef \mathbf{q}_2+\xi_n \mathbf{q}_1 $$
be the principal symbol of $\mathbf{P}$, and
$\trho=(\tx',\txi')$ be a point of $S^* 
\partial X$ such that the matrix $\mathbf{p}(\tx',0,\txi',\xi_n)$ is
invertible for any real number $\xi_n$. Consider a sequence $(u^k)$, weakly converging to
$0$ in $\Hv^1_{\loc}(\Xbar)$ and satisfying:  
$$ \mathbf{P} u^k=0, \quad u^k_{\rx}=0.$$
Then for all $M$:
$$ \partial_{x_n} u^k_{\rx} \underset{k\rightarrow 
+\infty}\longrightarrow 0 \text{ in } \Hv^M_{\trho,\partial}.$$ 
\end{prop}
We shall later apply the preceding proposition in the elliptic zone
($\EEE_T\cap \EEE_L$ with notations of paragraph \ref{mesures.Lame}) of the Lam\'e operator  $\partial_t^2-\Delta_e$. 

\subsubsection{In the glancing region}
The strong results of the two preceding paragraphs do not hold in the
neighbourhood of a glancing point. In this case, one need boundary
conditions to get further estimation than the standard traces theorem with
loss of one half-derivative. In
the case of Lopatinsky boundary conditions, the results are similar to
those of the hyperbolic region.  
\begin{prop}
\label{Lopa.traces}
{\bf a)} Let $\trho\in \GGG$. Assume that $(u^k)$ satisfies Lopatinsky
uniform boundary conditions near $\trho$. Then:
\begin{equation}
\label{Lopa.tr.O}
u^k_{\rx}=O(1) \text{ in } \Hv^1_{\trho,\partial},\quad 
\partial_{x_n}u^k_{\rx}=O(1) \dans \Lv_{\trho,\partial}.
\end{equation}
Furthermore, if $\mu$ vanishes near $\trho$, then:
\begin{equation}
\label{Lopa.tr.o}
u^k_{\rx}=o(1) \text{ in } \Hv^1_{\trho,\partial},\quad 
\partial_{x_n}u^k_{\rx}=o(1) \dans \Lv_{\trho,\partial}.
\end{equation}
{\bf b)} Assume:
\begin{equation}
\label{Hyp.b.G}
u^k_{\rx}=o(1) \text{ in } \Hv^1_{\trho,\partial},\quad 
\partial_{x_n}u^k_{\rx}=o(1) \dans \Lv_{\trho,\partial},
\end{equation}
and that every $\trho\in \GGG$ is not diffractive, in the sense that at
least one off the two half-bicharacteristic passing through $\trho$ stays
in $\partial \Omega$ near $\trho$. Then $\mu=0$ near $\trho$. 
\end{prop}
\begin{rem}
As seen in propositions \ref{traces.hyp} and \ref{traces.hyp.0}, point {\bf a)} holds in the hyperbolic case, where no boundary condition is
required.
\end{rem}
The proof of proposition \ref{Lopa.traces} is given in the appendix.

\subsection{The Lam\'{e} system}
\label{par.Lame}
This subsection is devoted to the Lam\'e system with Dirichlet boundary
conditions on an open bounded subset $\Omega$ of $\RR^3$:
\begin{equation}
\label{Lame}
\left\{ 
\begin{gathered}
\partial_t^2 u-\Delta_e u=0,\; (t,y)\in \RR\times \Omega\\
u_{\restriction \partial \Omega}=0\\
u_{\restriction t=0}=u_0 \in \Hv^1_0,\; \partial_t u_{\restriction 
t=0}=u_1 \in \Lv.
\end{gathered}
\right.
\end{equation}
In paragraph \ref{ondes.T.L}, (\ref{Lame}) is decomposed into two 
wave equations. In paragraph \ref{mesures.Lame}, we shall introduce the defect
measures associated to these equations. Next paragraphs are devoted to a few
elementary properties of these measures.
\subsubsection{Transversal and longitudinal waves}
\label{ondes.T.L}
The natural energy:
\begin{equation*}
E(t)\egaldef \frac 12 \int_{\Omega}\left(|\partial_t u|^2+\mu |\nabla 
u|^2+(\lambda+\mu)|\div u|^2\right) \dvol y
\end{equation*}
is time-invariant. Let $E_0$ be its constant value. The nest classical
proposition is proved, for example, in \cite{BuLe99}.
\begin{prop}[Decomposition of the Lam\'e system]
\label{decompo}
There exists a constant $C>0$ such that for every solution $u$ of
(\ref{Lame}), there exists:
$$ u_T\in \Hv^1_{\loc}\left(\overline{M}\right),\; u_L \in 
\Hv^1_{\loc}\left(\overline{M}\right),$$
such that:
\begin{enumerate}
\item $u=u_T+u_L,\; \div u_T=0,\;\rot u_L=0.$
\item $(\partial_t^2-c_T^2\Delta)u_T=0, \quad$where $c_T^2\egaldef\mu.$
\item $(\partial_t^2-c_L^2\Delta)u_L=0, \quad$where $c_L^2\egaldef\lambda+2\mu.$
\item For every bounded interval $I$ of $\RR$, of length $|I|$:
$$\| u_L\|^2_{\Hv^1(I\times\Omega)}+\| 
u_T\|^2_{\Hv^1(I\times\Omega)}\leq C|I|E_0.$$ 
\item If $u$ is in the space vector generated by a finite number of
  eigenfunction of $\LLL$, then: 
$$u_T\in C^{\infty}(\overline{M}),\quad u_L\in C^{\infty}(\overline{M}).$$ 
\end{enumerate}
\end{prop}
\begin{Def}
The function $u_T$ is called {\bf transversal wave}, and the function 
$u_L$ {\bf longitudinal}.
\end{Def}
\begin{rem}
In the sequel we shall often reduce the longitudinal wave to a scalar
function, writing $u_L=\nabla \varphi$, with:
\begin{gather*}
\varphi \in \Hv^2_{\loc}\left(\overline{M}\right), \quad 
\|\varphi\|_{\Hv^2(I\times\Omega)} \leq C |I| E_0,\quad
(\partial_t^2-c_L^2\Delta)\varphi=0. 
\end{gather*} 
\end{rem}
\subsubsection{Measures}
\label{mesures.Lame}
Let $(u^k)$ be a sequence of solutions of the Lam\'e system with:
\begin{equation*}
(u_0^k,u_1^k) \tendkf 0 \dans \Hv_0^1(\Omega)\times \Lv(\Omega). 
\end{equation*}
In view of the continuity of the map introduce in proposition \ref{decompo}:
\begin{eqnarray*}
(u_0,u_1) &\longmapsto &(u_L,u_T)\\
\Hv_0^1\times \Lv &\longrightarrow & \left(\Hv^1(I\times \Omega)\right)^2,
\end{eqnarray*}
where $I$ is a bounded interval, both sequences  
$(u_T^k)$ and $(u_L^k)$ weakly converge to $0$ in 
$\Hv^1_{\loc}(\RR\times \overline{\Omega})$. Likewise, the sequence 
$(\partial_t\varphi^k)$ weakly converges to $0$ in $H^1_{\loc}(\RR 
\times \overline{\Omega})$ (it is more convenient to consider a derivative of
$\varphi$ in order to work in a $H^1$ space, for which the defect measures
introduced here
 are well fitted).
\begin{lem}
\label{traces.lame.0}
The sequences $(u_T^k)$, $(u_L^k)$, $(\partial_t \varphi^k)$ are regular on
the boundary.
\end{lem}
(cf \cite[lemme 4.2]{BuLe99})
\begin{notations}
Let:
\begin{itemize}
 \item $\mu_T$, $\mu_L$ and $\mu$ be the defect measures respectively
   associated (up to a subsequence) to  
$(u_T^k)$, $(u_L^k)$, and 
$(\partial_t \varphi^k)$ by theorem \ref{existence.mu};
\item
$ {\displaystyle \HHH_T,\;\HHH_L,\;\GGG_T,\;\GGG_L,\;\EEE_T,\;\EEE_L} $
the hyperbolic, glancing and elliptic region of  
the transversal and longitudinal waves.
\end{itemize}
All the calculation shall be carried out in one of the $J+1$ local
coordinate systems choosen in the beginning of this section. We shall make
a distinction between the spaces of scalar operators:
$\AAA^m,\;\AAA^m_{\partial}$ (defined in paragraphs  \ref{a.l.interieur}, \ref{au.bord}, \ref{par.recol},  
with $N=1$), and the spaces of matricial operators
$\AAAv^m,\;\AAAv^m_{\partial}$ (with $N=3$).\par
The notation $x=(y,t)$ always refers to local coordinates. When a
distinction is necessary, we shall write coordinates on $\Omega$ before the
change of variables $z=(z_1,z_2,z_3)$. This global system of
coordinates has been chosen so that the magnetic field is vertical: $ \Bbf=(B,0,0).$
\end{notations}
\begin{rem}
One may identify $\mu_T$ and $\mu_L$ to measures on 
$S^*_b\left(\RR\times\overline{\Omega}\right)$ with values endomorphism
of $\CC^3$.
\end{rem}
\begin{rem}
Condition $c_L^2\neq c_T^2$ means that the intersection of
$\GGG_T$ and $\GGG_L$ is empty. 
\end{rem}
\subsubsection{Link between $\mu$ and $\mu_L$}
Let $\chi$ be the local diffeo\-morphism from glo\-bal spa\-tial coordinates
$(z_1,z_2,z_3)$ to local coordinates:
$$ z=\chi(y),\quad \transp \chi'(y)\zeta=\eta.$$
\begin{prop}
\label{mu.muL}
$$ \forall a\in C^0_0\left(S^*_b(\RR\times\overline{\Omega}),M^3(\CC)\right),\; \big< \mu, \frac{\transp{\zeta} a \zeta}{\tau^2} \big>=\big< \mu_L,a\big>.$$
In particular measures $\mu$ and $\mu_L$ have the same supports.
\end{prop}
\begin{proof}
Let $ A_j \in\AAAv^j,\; j=1,2,\quad A\egaldef A_{-1} D_{x_n}+A_0$.\\
Set: $I_k\egaldef (A\partial_t u_L^k,\partial_t u^k_L)$.\\
On one hand:
$$ I_k=-(\partial_t A\partial_t u_L^k,u^k_L)\tendk \left< \mu_L,a \right>.$$
On the other hand:
$$ I_k=-(\div A\nabla \partial_t \varphi^k,\partial_t \varphi^k)+o(1)\tendk
\left< \mu_L,\frac{\transp{\zeta} a\zeta}{\tau^2}\right>.$$
The boundary terms of this preceding integration by parts converge to $0$
according to lemma \ref{traces.lame.0}. This implies proposition
\ref{mu.muL} when $a$ is of the form $\xi_n a^{-1}+a_0$, and then by a density
argument for any $a$.
\end{proof}
\subsubsection{Polarization of $\mu_T$ and $\mu_L$}
Let $\pi$ bet the orthogonal projection in $\CC^3$ on the line generated by
$\zeta$, and $\pi_{\bot}$ the orthogonal projection on the plane normal to $\zeta$.  
\begin{equation}
\label{def.pi}
\pi(V)\egaldef |\zeta|^{-2}(\transp \zeta.V) \zeta,\quad \pi_{\bot}\egaldef \Id_{\CC^3}-\pi.
\end{equation}
Projectors $\pi$ and $\pi_{\bot}$ are defined by 
formulas (\ref{def.pi}), on $S^*\overline X$. 
\begin{prop}
The measure $\mu_L$ is polarized along the direction of
propagation, and $\mu_T$ orthogonally to this direction:
$$ \mu_L=\pi\mu_L\pi,\quad \mu_T=\pi_{\bot} \mu_T \pi_{\bot}.$$
\end{prop}
\begin{proof}
The statement on $\mu_L$ is an immediate consequence  of proposition
\ref{mu.muL}. To show the statement on $\mu_T$, take $A_0\in \AAA_0$. The
nullity of $\div  u_T^k$ implies: 
$$ 0=(A_0\nabla \div u_T^k,u_T^k)\tendk \big< \mu_T, \frac{a_0\zeta \transp{\zeta}}{\tau^2}\big>.$$
Thus:
$$ <\mu_T,a\pi>=0,\quad <\mu_T,a\pi_{\bot}>=<\mu_T,a>.$$
The symmetry of $\mu_T$ completes the proof.
\end{proof}
\begin{rem}
To get more intrisic formulations of the preceding results, i.e. statements
where the two
coordinate system do not mix, one should have considered $u_L^k$ and $u_T^k$ as
section of the tangent space $T\Omega$, and defined measures with values
endomorphism of $T\Omega$ (instead of endomorphism of $\CC^3$).   
\end{rem}
\subsubsection{A decoupling lemma}
\label{decouplage}
The next result, converting an approximate differential equation on $u^k$
into two equations on $u^k_L$ and $u^k_T$, is of crucial importance in the
sequel. As before, $(.,.)$ stand for the $\Lv$ scalar product on $\RR\times\Omega$.  
\begin{lem}[decoupling lemma]
\label{lem.decouplage}
Let $A$ be a pseudo-differential operator of order $2$ of the following form:
\begin{gather}
\label{forme.A}
A=\sum_{j=-M}^2 A_j \partial_{x_n}^{2-j}+A_i\\
\notag
A_i\in \AAAv^2_i,\quad A_j \in \AAAv^j.
\end{gather}
Then:
\begin{equation}
\label{decouplage.a}
\lim_{k\rightarrow +\infty} (Au_T^k,u_L^k)= 
\lim_{k\rightarrow +\infty}(Au_L^k,u_T^k)=0. 
\end{equation}
Let $A$ be a $(1,3)$ matrix of pseudo-differential operators, with
coefficients of the form (\ref{forme.A}), but with {\bf scalar} operators. Then:
\begin{equation}
\label{decouplage.b}
\lim_{k\rightarrow +\infty} (Au_T^k,\partial_t \varphi^k)=0.
\end{equation}
\end{lem}
\begin{proof}
We shall only prove the convergence to $0$ of $(Au_T^k,u_L^k)$. The proof
of rest of the lemma is very much the same. We may obviously assume that
the oprators $A_j$ have compact support in one of  the local coordinate
system introduced in paragraph \ref{par.recol}. In view of the equations:  
\begin{equation*}
-\detg^{-1/2} \partial_{x_n} \detg^{1/2} \partial_{x_n} u^{k}_{T} 
+Q_{T}u^{k}_{T}=0,\quad -\detg^{-1/2} \partial_{x_n} \detg^{1/2} 
\partial_{x_n} u^{k}_{L} +Q_{L}u^{k}_{L}=0,
\end{equation*}
where $Q_T$ and $Q_L$ are tangential, it is sufficient  to prove the lemma
in the cases $j=1,2$ and in the interior case.\par
 {\bf First case:} $A\in \AAA^i$.\par 
In view of: $\nu_L\neq\nu_T$ it is easy to construct two operators:
\begin{gather*}
\Psi_T,\Psi_L\in \AAA_i^0,\quad (\Psi_T+\Psi_L)_{\restriction U}= 
{\rm Id}\\
\supp \sigma_0(\Psi_T)\cap \ZZZ_L=\supp \sigma_0(\Psi_L) \cap 
\ZZZ_T=\emptyset,
\end{gather*} 
where $U$ is a small open subset of $\Omega$ such that there exists a function
$\varphi\in C^{\infty}(U)$ satisfying $\varphi A \varphi=A$. Writing:
$$A=A\varphi=A\Psi_L\varphi+A\Psi_T\varphi,$$ 
we may assume that the principal symbol of $A$ does not intersect $\ZZZ_T$
(or does not intersect $\ZZZ_L$). For such operators, (\ref{decouplage.a})
holds trivially. For example, in the first case we have: 
$$ Au^k_T=O(1) \dans \Lv.$$\par
{\bf Second case}: $A\in \AAA^2_{\partial}$.\par
We know that the support of $\mu_T\indic_{\{x_n=0\}}$ is included in
$\GGG_T$ and the support of $\mu_L\indic_{\{x_n=0\}}$ in $\GGG_L$.\\
As a consequence, we may write $A=A(\Theta_T+\Theta_L)$ where $\Theta_{T}$ and 
$\Theta_L$ are tangential operators of degree $0$ such that:
$$\supp \sigma_0(\Theta_T) \cap \GGG_L=\supp \sigma(\Theta_L) \cap 
\GGG_T=\emptyset. $$
We may thus assume that the support of the principal symbol $a_2$ of $A$ is
 disjoint from one of the two glancing sets, say $\GGG_T$. We have:
$$ \big{(}A u_T^{k},u_L^{k}\big{)}=\big{(}\chi(x_n/\eps)A 
u_T^{k},u_L^{k}\big{)} 
+\big{(}(1-\chi(x_n/\eps))Au^k_T,u_L^{k}\big{)}.$$
where $\chi$ is a compactly supported function in $\RR$ equal to $1$ near
the origin. We first fix $\eps$ and let $k$ tend to $\infty$. The second
term of the sum tends to $0$ in view of  the preceding case. As for the
first term, we have, $u_L^k$ being bounded in $\Hv^1$: 
\begin{align*}
|\big{(}\chi(x_n/\eps)A
u_T^{k},u_L^{k}\big{)}|\leq &
C\|\chi(x_n/\eps)Au_T^{k}\|_{L^2(0,l,\Hv^{-1}(X'))}\| 
u_L^k\|_{L^2(0,l,\Hv^{1}(X'))}\\ 
\leq & C\| 
\Lambda'_{-1} \chi(x_n/\eps) A u_T^{k}\|_{\Lv} +o(1),\; 
k\rightarrow +\infty,
\end{align*}
where $\Lambda'_{-1}\in \AAA_{\partial}^{-1}$, with principal symbol equal
to $\|\xi'\|^{-1}$ near the support of $A$. We have:
\begin{gather*}
\|\Lambda'_{-1}\chi(x_n/\eps)A\,u_T^{k}\|_{\Lv}^2
\tendk \Big{<}
\mu_T,\frac{(\chi(x_n/\eps))^2\,|a_2|^2}{\tau^2\|\xi'\|^2}\Big{>}\\
\limsup_{k\rightarrow+\infty}\,\big{|}\big{(}A 
u_T^{k},u_L^{k}\big{)}\big{|}\leq
\Big< \mu_T,\frac{(\chi(x_n/\eps))^2\,|a_2|^2}{\tau^2\|\xi'\|^2}\Big{>}.
\end{gather*}
When $\eps$ goes to $0$, the right side of this inequality converges (by
the dominated convergence theorem) to:  
$$ \big<\mu_T, \indic_{\{x_n=0\}} |a_2|^{2}\tau^{-2}\|\xi'\|^{-2} \big>,$$
which is null, because $\GGG_T$ and the support of $a_2$ are disjoint.
\par
{\bf Third case}: $A=A_1 D_{x_n}$, $A_1$ tangential.\par
As in the preceding case, we  may assume that $\sigma_1(A_1)$ is disjoint
with one of the two glancing sets, say $\GGG_T$. Then:
$$ \big{(}A u_T^{k},u_L^{k}\big{)}=\big{(}\chi(x_n/\eps)A
u_T^{k},u_L^{k}\big{)}+\big{(}(1-\chi(x_n/\eps))A
u_T^{k},u_L^{k}\big{)}$$
The second term converges to zero when $k\CVF \infty$ for the same reasons as
in the case $A\in\AAA_i^2$. The first term may be written:
\begin{align*}
 \left(\chi(x_n/\eps)A_1 D_{x_n} u_T^{k},u_L^{k}\right)&=
\left(D_{x_n}(\chi(x_n/\eps)A_1
u_T^{k}),u_L^{k}\right)+\left( R_\eps
u_T^{k},u_L^{k}\right),\; R_{\eps} \in \AAAv^0\\
&=\left((\chi(x_n/\eps)A_1
u_T^{k}),D_{x_n}u_L^{k}\right)+\left( R_\eps
u_T^{k},u_L^{k}\right)+\{boundary\;terms\}
\end{align*}
The boundary terms tend to zero when $k$ tends to infinity because
$(u^k_T)$ and $(u^k_L)$ are regular on the boundary. The proof may be completed as in the preceding case,
letting $k$ go to infinity then $\eps$ go to zero. 
\end{proof}

\section{Sufficient condition}
\label{chap.CS}
If $\gamma$ is a $\Bbf$-resistant ray defined on a real interval
$]a,b[$, we shall call life-length the positive quantity $|t(b)-t(a)|$.
In this section we use tools of the preceding section to prove the following:
\begin{prop}
\label{propCS}
Let $T>0$. Assume that every $\Bbf$-resistant ray in $\Omega$ is of life-length strictly less than $T$. Then there exists $C>0$ such that for evey solution
of the Lam\'e system (\ref{Lame}):
\begin{equation}
\label{inegCS2}
\|u_0\|^2_{\Hv^1_0}+\|u_1\|^2_{\Lv}\leq C \left( \|\rot(\partial_t 
u\wedge 
\Bbf)\|^2_{\Hv^{-1}((0,T)\times \Omega)} 
+\|u_0\|^2_{\Lv}+\|u_1\|^2_{\Hv^{-1}}\right)
\end{equation}
\end{prop}

Inequality (\ref{inegCS2}) is the sufficient condition (\ref{inegCS}) for uniform decay
stated in point a) of proposition \ref{CNSineg}. Proposition \ref{propCS}
thus completes the proof of the sufficient condition of theorem \ref{th1},
namely that the non-existence of arbitrarily large  $\Bbf$-resistant rays on
$\Omega$ implies the uniform decay of the energy for solutions of the system
of magnetoelasticity.\par
To show (\ref{inegCS2}), we shall argue by contradiction, considering the
defect measures $\mu_{T,L}$ of subsection  \ref{mesures.Lame} associated to
a sequence $(u^k)$ which contredicts (\ref{inegCS2}) (cf subsection
\ref{CS.absurde}). The bound on $\rot(\partial_t u^{k}\wedge \Bbf)$ given
by the negation of (\ref{inegCS2}) implies a strong condition on the supports
 of these measures (see subsection \ref{CS.support}). In subsection
 \ref{CS.bord}, we make use of this condition, together with propagation
 arguments near the boundary of $\Omega$. Subsection \ref{CS.conclu}
 completes the proof, using the assumption of non-existence of
 $\Bbf$-resistant rays of life-length larger than T.

\subsection{Introduction of measures}
\label{CS.absurde}
Assume that (\ref{inegCS2}) does not hold. Then there exists a sequence
$(u^k)$ of solutions of the Lam\'e system such that:
\begin{equation}
\label{CS.a.ineg}
1=\|u_0^k\|^2_{\Hv^1_0}+\|u_1^k\|^2_{\Lv}> k \left( \|\rot(\partial_t 
u^k\wedge 
\Bbf)\|^2_{\Hv^{-1}((0,T)\times \Omega)} 
+\|u^k_0\|^2_{\Lv}+\|u_1^k\|^2_{\Hv^{-1}}\right).
\end{equation} 
Up to the extraction of a subsequence, one may assume that $(u_0^k,u_1^k)$
weakly converges in $\Hv^1_0\times \Lv$. Inequality (\ref{CS.a.ineg})
implies that its weak limit is $0$. We may thus introduce the defect measures
$\mu_T$, $\mu_L$ and $\mu$ of paragraph \ref{mesures.Lame}, associated to
the sequences $(u^k_T)$, $(u^k_L)$ and $(\partial_t \varphi^k)$. To
contredict (\ref{CS.a.ineg}), we need to show that these measures are
null. Note that (\ref{CS.a.ineg}) implies:
\begin{equation}
\label{rot.CVF.0}
\rot(\partial_t u^k\wedge \Bbf) \tendk 0 \dans H^{-1}((0,T)\times 
\Omega).
\end{equation}
\begin{rem}
By a density argument, it suffices to show 
(\ref{inegCS2}) with $(u_0,u_1)$ generated by a finite number of
eigenfunctions of $\LLL$. We may thus assume, that $u_L^k$ and
$u_T^k$ are $C^{\infty}$.
\end{rem}
\begin{rem}
\label{rem.CS}
We will indeed show a more precise statement than proposition \ref{propCS}
namely that if $(u^k)$ is a sequence of solutions of the Lam\'e system
converging weakly to $0$ in the energy space and satisfying
(\ref{rot.CVF.0}) then the set $(\supp \mu_T\cup \supp \mu_L)\cap
\{t\in (0,T)\}$ a an union of $\Bbf$-resistant rays of length $T$.
\end{rem}
\subsection{Condition on the supports}
\label{CS.support}
We may see $\Bbf$ as a vector field on $\Omega$,  i.e. a section of $T\Omega$.
To avoid confusions, the magnetic field considered as a vector field shall
be refered as $\Bv$. In a local coordinate system, if:
 $$ \chi:U\subset \overline{\Omega} \longrightarrow \overline{\RR}_+^{n}$$
is the change of coordinates, and $\chi'$ its differential, $\Bv$ is equal
to $\chi' \Bbf$. Notation $\Bbf$ shall always refer to the vector of
$\RR^3$ of coordinates $(B,0,0)$. As before, $(z_1,z_2,z_3)$ refers to the
global spatial coordinates on $\Omega$, before the change of variable.
\begin{lem}
\label{lem.supp.mu}
Assume (\ref{rot.CVF.0}). Then, on the interval $(0,T)$, the transversal
measure charges set of all points
whose direction of propagation is orthogonal to $\Bbf$ and the longitudinal
measure charges set of all points whose direction of propagation is
parallel to $\Bbf$.
\begin{align}
\label{supp.muT}
\mu_T \indic_{(0,T)}=\mu_T \indic_{(0,T)} \indic_{\Bv^{\bot}},\quad
 \Bv^{\bot}:= 
\{(t,y,\tau,\eta);\;\transp{\Bv}\eta=0\}\\
\label{supp.muL}
\mu_L \indic_{(0,T)} =\mu_L 
\indic_{(0,T)}\indic_{\Bv^{\sslash}},\quad \Bv^{\sslash}:= 
\{(t,y,\tau,\eta);\; \eta 
\in
\text{vect}\,(g\Bv)\}.
\end{align}
\end{lem}
\begin{proof}
Set:
\begin{equation}
\label{Ru=rot}
Ru\egaldef \rot(\partial_t u\wedge \Bv)= \partial_t \left( 
\begin{array}{c}
-\partial_{z_2} u_2-\partial_{z_3} u_3\\
\partial_{z_1} u_2\\
\partial_{z_1} u_3
\end{array}\right).
\end{equation}
\underline{Transversal measure}.
The measure $\mu_T$ does not charge neither $\HHH_T$ nor $\EEE_T$. Thus, it
suffices to check:
$$ \supp \mu_T\indic_{(0,T)}\indic_{\{x_n>0\}\cup \GGG_T}\subset 
\Bv^{\bot}.$$
Near the boundary, by proposition \ref{prop1.appendice} of the appendix, 
$$ Ru^k \tendk 0 \dans L^2([0,l[,\Hv_{\loc}^{-1}(X')).$$
Thus, according to the decoupling lemma, whether $A_0$ has support in the interior or near the boundary:
\begin{equation}
\forall A_0 \in \AAA^0,\; \supp A_0\subset\{t\in(0,T)\},\; (A_0 R 
u_T^k,u_T^k) \tendk 0.
\end{equation} 
Because $\div u_T^k=0$, formula (\ref{Ru=rot}) may be written: $ Ru_T^k=\partial_t \partial_{z_1} u_T^k.$
This implies (by paragraph \ref{CS.mu.0}):
$$ \mu_T\Big( \{t\in (0,T)\} \cap \left\{ \sigma_2(\partial_t \partial_{z_1}) \neq 0 \right\}\Big)=0.$$
This clearly shows the announced result near an interior point. When
$\trho\in \GGG_T$ one may write:  $\partial_{z_1}=f_0
\partial_{x_n} +F_1$, where $f_0$ is a function and $F_1$ a first order
tangential differential operator, and the following basic fact completes
the proof: 
\begin{equation*}
\sigma(F_1)(\trho)=0 \iff \teta' \bot \Bv.
\end{equation*}
 \underline{Longitudinal measure}.
The first coordinate of $Ru_L^k$ is:
$$ -\partial_t(\partial_{z_2}^2+\partial_{z_3}^2) \varphi^k.$$
Its scalar product with $\partial_t\varphi^k$ gives:
$$\forall A_0 \in \AAA^0,\; \supp(A_0)\subset \{t\in (0,T)\} 
\Rightarrow \lim_{k\rightarrow 
+\infty}\left(A_0(\partial_{z_2}^2+\partial_{z_3}^2)\partial_t 
\varphi^k,\partial_t \varphi^k\right)=0.$$
This implies (again by paragraph \ref{CS.mu.0}), that $\mu$ (thus $\mu_L$)
vanishes, in $(0,T)$, on the set of all $\rho$ such that:
$$ \rho\in \{x_n>0\}\cup \GGG_T,\quad 
\sigma_2(\partial_{z_2}^2+\partial_{z_3}^2)(\rho)\neq 0.$$
Hence (\ref{supp.muL}).
\end{proof}
\subsection{Support of the measure near the boundary}
\label{CS.bord}
\begin{figure}
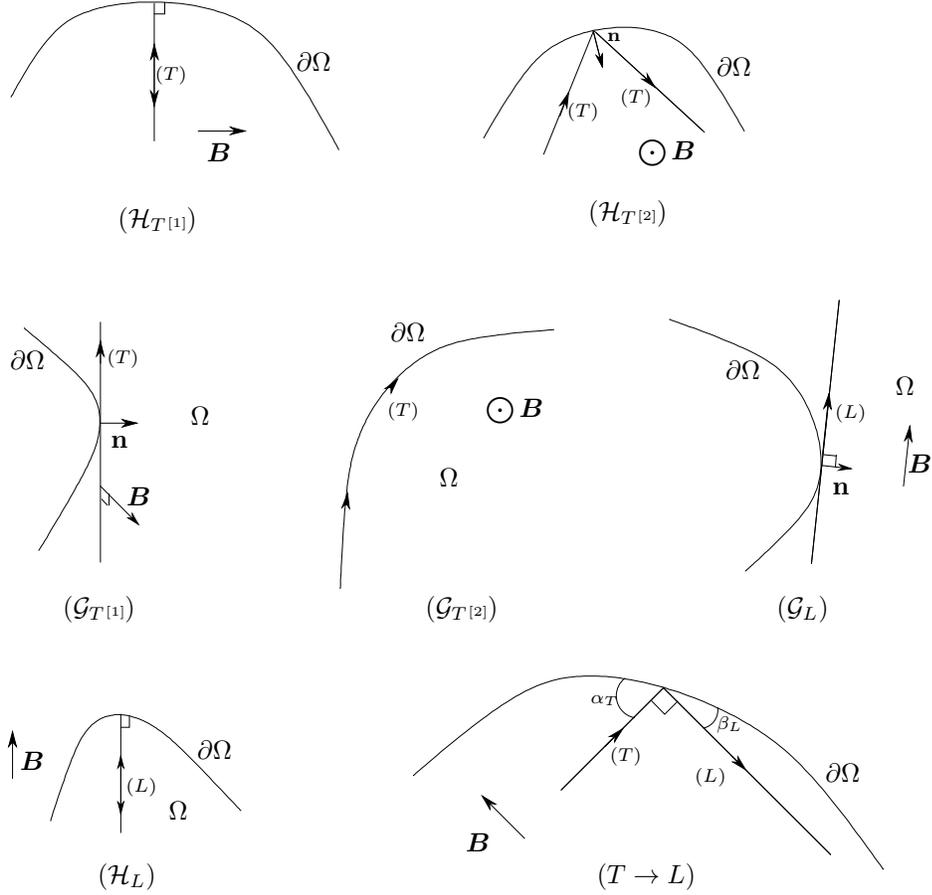

\begin{center}
\dessCS
\end{center}
\caption{Cases arising in lemma \ref{lem.CS.bord}}
\label{fig.CS}
\end{figure}

For any symbol $q_0$ with support in $\{x_n>0\}$, we have:
$$ \big<\mu_T,H_{p_T/\tau} q_0\big>=0 \qquad \big<\mu_L,H_{p_L/\tau}
q_0\big>=0, $$ 
which shows, in the interior of $\Omega$, the invariance of each measure by
the appropriate hamiltonian flow. Unfortunately, the condition:
\begin{equation}
\label{Cbord}
u_{T\restriction \partial 
\Omega}^{k}+u_{L\restriction\partial\Omega}^{k}=0
\end{equation}
is not always sufficient to extend such a property in the neighbourhood of
$\partial \Omega$. Indeed, without any additional assumption,  $\mu_L$ and
$\mu_T$ are not deterministic: the value of the two measures for time
$t< t_0$ is not uniquely determined by their value for time $t\leq
t_0$. In our case, this convenient property holds thanks to the strong
conditions on the support of $\mu_T$ and $\mu_L$. As announced before, we
shall only describe the propagation of the supports of the measure. 
\begin{lem}
\label{lem.CS.bord}
Let $\mu_T$ and $\mu_L$ be the defect measures associated to a sequence  of
 solutions of the Lam\'e system satisfying (\ref{rot.CVF.0}). Let:
$$ \trho=(\tx',\txi')=(\tilde{t},\ty',\ttau,\teta')\in S^*\partial X$$
and $\bfn$ the unitary exterior normal vector to $\partial \Omega$ at
$\ty'$. Them $\mu_T$ and $\mu_L$ both vanish near $\trho$ except possibly in
the following cases (cf figure \ref{fig.CS}): 
\begin{enumerate}
\item \label{bord.mu.T} $\mu_L$ is null. 
The support of $\mu_T$ propagates along the transversal flow and:
\begin{itemize}
\item \underline{$(\HHH_T\un)$ case}: $\trho\in \HHH_T$, $\teta'=0$ and 
$\Bv$ is
orthogonal to $\bfn$;
\item \underline{$(\HHH_T\deux)$ case}: $\trho\in \HHH_T$, $\teta'\neq 
0$ and $\Bv$ is 
normal to the reflection plane;
\item \underline{$(\GGG_T\un)$ case}: $\trho$ is diffractive for the
  transversal wave (i.e. $\trho\in \GGG_T$ and the bicharacteristic ray
  passing through $\trho$ only intersect the boundary at $\trho$), and
  $\Bv$ is orthogonal to $\teta'$; 
\item \underline{$(\GGG_T\deux)$ case}: $\trho\in \GGG_T$ is not
  diffractive for the transversal wave, and $\Bv$ is normal to the
  reflection plane.
\end {itemize}
\item \label{bord.mu.L} $\mu_T$ is null, the support of $\mu_L$ propagates
  along the longitudinal flow and:
\begin{itemize}
\item \underline{$(\HHH_L)$ case} $\trho\in \HHH_L$, $\teta'=0$ and 
$\Bv$ 
is parallel to $\bfn$;
\item \underline{$(\GGG_L)$ case} $\trho$ is a diffractive point for the
  longitudinal wave and $\Bv$ is parallel to $\teta'$.
\end{itemize}
\item \label{bord.mu.T.L} Both measures $\mu_T$ and $\mu_L$ are non null,   
$\trho\in \HHH_T\cap\HHH_L$ and:
\begin{itemize}
\item \underline{$(T\rightarrow L)$ case}: $\Bv$ is orthogonal to the
  transversal ray coming in, and parallel to the longitudinal ray going
  out of $\trho$. The support of $\mu_T$ is an union of incoming transversal
  rays. The support of $\mu_L$ is the union of all outgoing longitudinal
  rays going out of points of $\HHH_T\cap \HHH_L$ where the transversal
  rays of the support of $\mu_T$ come in; 
 \item \underline{$(L\rightarrow T)$ case}: $\Bv$ is orthogonal to the
  transversal ray going out of and parallel to the longitudinal ray coming in $\trho$. The support of $\mu_L$ is an union of incoming longitudinal
  rays. The support of $\mu_T$ is the union of all outgoing transversal
  rays going out of points of $\HHH_T\cap \HHH_L$ where the longitudinal
  rays of the support of $\mu_L$ come in. 
\end{itemize}
\end{enumerate}
\end{lem}
All the assertions of lemma \ref{lem.CS.bord} should be understood in a
neighbourhood of $\trho$. The reflection plane at a boundary point
$\rho=(t,y',\tau,\eta')$, defined as long as
$\eta'\neq 0$, is the plane passing through $y'$ and generated by $\bfn$
and $\eta'$, thus containing the bicharacteristic ray passing through
$\rho$. The statement ``$\Bv$ is parallel to $\eta'$'' must be understood
as ``the vector $g\Bv$ of the cotangent bundle of $\Omega$ is parallel to $\eta'$''.
\begin{notation}
Let $\rho$ be an hyperbolic point for the transversal (respectively
longitudinal) wave. We shall denote by $\xi^-_T$, 
$\xi^+_T$ (respectively $\xi^-_L$, $\xi^+_L$) the incoming and outgoing
vectors through $\rho$:
\begin{equation*}
\xi^+_T\egaldef \left(
\begin{array}{c}
\xi' \\
\xi_{nT}=\sqrt{\nu_T^2 \tau^2-\|\eta'\|^2}
\end{array}
\right) 
\quad
\xi^-_T=\left(
\begin{array}{c}
\xi' \\
-\xi_{nT}
\end{array}
\right) 
\end{equation*} 
(respectively with ``$L$'' instead of ``$T$'').\par 
We shall write $\eta^{\pm}_{T}$,  $\eta^{\pm}_{L}$ the spatial
components of this vectors. For example:
$$\eta^+_T=\left(
\begin{array}{c}
\eta' \\
\sqrt{\nu_T^2 \tau^2-\|\eta'\|^2}
\end{array}
\right).$$
\end{notation}
Let's postpone the proof of lemma \ref{CS.bord} to show, as stated in the
introduction of this article, that in the ($T\rightarrow L$) and ($L\rightarrow T$) cases,
the angles of refraction and incidence have a fixed value, determined by
the quotient $c_T/c_L$. Consider for example the ($T\rightarrow
L$) case. Let $\alpha_T$ be the angle of incidence of the transversal wave,
$\beta_L$ the angle of refraction of the longitudinal wave, $a_T$ and $a_L$
the following numbers:
$$ a_T\egaldef \tan \alpha_T=\frac{\txi_{nT}}{\|\teta'\|},\quad a_L\egaldef \tan \beta_L=\frac{\txi_{nL}}{\|\teta'\|}.$$
The incident and refracted waves are orthogonal, so that:
\begin{equation}
\label{eq.aT.aL.1}
\|\teta'\|^2-\txi_{nT}\txi_{nL}=0,\;\text{i.e.}\; a_T a_L=1.
\end{equation}
Furthermore, the definition of $\txi_{nL}$ and $\txi_{nT}$ yields
$$ c_T^2\|\teta^+_T\|^2-\tau^2=0,\quad c_L^2\|\teta_L^-\|^2-\tau^2=0,$$
which gives the equation:
\begin{equation}
\label{eq.aT.aL.2}
c_T^2(1+a_T^2)=c_L^2(1+a_L^2).
\end{equation}
Equations (\ref{eq.aT.aL.1}) and (\ref{eq.aT.aL.2}) imply the formula
announced in the introduction:
\begin{equation*}
 \alpha_T=\arctan \frac{c_L}{c_T},\; \beta_L=\arctan\frac{c_T}{c_L}.
\end{equation*}
By a similar calculation, one gets, in the ($T\rightarrow L$) case: 
\begin{equation*}
 \alpha_L=\arctan \frac{c_T}{c_L},\; \beta_T=\arctan\frac{c_L}{c_T}.
\end{equation*}
There are thus very strong constraints for the possible transfer of energy
from one wave equation to the other.
\begin{proof}[Proof of lemma \ref{lem.CS.bord}]
\underline{Case (\ref{bord.mu.T.L})}: $\mu_T\neq 0$, $\mu_L\neq 0$ 
near $\trho$.\par 
In this  case, $\trho \notin \EEE_T\cup \EEE_L$. It is also easy to show
that $\trho \notin \GGG_T\cup \GGG_L$. Indeed, if $\trho\in \GGG_T$ then it
also belongs to $\HHH_L$ (it cannot be a point of $\EEE_L$, and $\GGG_L$ and
$\GGG_T$ are disjoint). But $\mu_T$ being non-null near $\trho$, $\teta$ is
orthogonal to $\Bv$ (by lemma \ref{lem.supp.mu}) so neither  $\teta_L^+$
nor $\teta_L^-$ are parallel to $\Bv$, which implies (again by lemma
\ref{lem.supp.mu}) that $\mu_L=0$ near $\trho$, contradicting our
assumptions. Likewise, if $\trho \in \GGG_L$, $\teta'$ must be parallel to
$\Bv$ and $\mu_T$ null near $\trho$. 
Thus  $\trho \in \HHH_T \cap \HHH_L$. The support of measures $\mu_T$ and
$\mu_L$ is, near $\trho$, an union of incoming and outgoing maximal
rays.\par
Let's first assume that the support of $\mu_L$ contains the
ray going out of $\trho$. Then:
$$\teta^+_L \sslash \Bv$$ 
so that $\teta^+_T$ is not orthogonal to $\Bv$. As a consequence, the
support of $\mu_T$ is only made of incoming rays, and the fact that
$\mu_T\neq 0$ implies:
$$ \teta^-_T \bot \Bv.$$
Thus $\teta^-_L$ is not parallel to $\Bv$. This is the $(T
\rightarrow L)$ case, and it remains to show the statement of the lemma
about the transfer from transversal incoming waves to longitudinal outgoing
waves, which may be formulated as follow: for any $\dot{\rho}\in \HHH_T
\cap \HHH_L$ near $\trho$ the following equivalence holds:
\begin{equation}
\label{dot.rho}
\dot{\rho}\in \supp \mu_T \iff \dot{\rho} \in \supp\mu_L. 
\end{equation}
Let's assume for example $\dot{\rho} \notin \supp \mu_T$. Then $\mu_T$ is 
null near rays coming in and going out of $\dot{\rho}$ and by the
hyperbolic theory (see proposition \ref{traces.hyp.0}):  
$$ u^k_{L\rx}=-u_{T\rx}^k \tendk 0 \dans 
\Lv_{\dot{\rho},\partial}.$$  
By the propagation theorem, the support of $\mu_L$ propagates near
$\dot{\rho}$. But this support is an union of outgoing rays. Consequently,
it is empty near $\dot{\rho}$ and $\mu_L$ is null near $\dot{\rho}$.  The implication
$ \dot{\rho} \notin \supp \mu_L\Rightarrow \dot{\rho} \notin \supp \mu_T $
may be shown in the same manner.\par
If the support of $\mu_L$ contains no outgoing ray near $\trho$,
it must contain incoming rays. This corresponds to the $(L\rightarrow T)$ case, which may be treated
as the $(T\rightarrow L)$ case.\par 
The study of all other cases relies on a transfer argument on boundary
conditions, stated in the following technical lemma \ref{transfert.trace}:
rougly, a boundary condition on the longitudinal wave implies one on the
transversal wave and vice versa. 
\begin{notation}
Let $v_{T,L}^k$ be the functions $u_{T,L}^k$ 
considered as vector fields on $\Omega$. In local coordinates, if $\chi$
denotes the change of coordinates, we have:
$$ v^k_T=\chi'(y) u^k_T=\left(\begin{array}{c} v^k_{T1}\\v^k_{T2}\\v^k_{Tn}\end{array}\right)  \quad v^k_L=\chi'(y) u^k_L=\left(\begin{array}{c} v^k_{L1}\\v^k_{L2}\\v^k_{Ln}\end{array}\right) .$$ 
\end{notation}
\begin{lem}
\label{transfert.trace}
Let $\trho\in S^*\partial X$, and $(u^k)$ be any sequence of solutions of the
Lam\'e system weakly converging to $0$ in $\Hv^1_{\loc}(\RR\times\overline \Omega)$.\par
\begin{itemize}
\item
Assume the following aproximate equation for some $A_1\in \AAA^1$:
\begin{equation}
\label{t.t.1}
\partial_{x_n} v^k_{Tn\rx}=A_1 v^k_{Tn\rx}+o(1) \dans L^2_{\trho,\partial}.
\end{equation}
Then:
\begin{equation}
\label{t.t.2}
\tag{\ref{t.t.1}'}
\Delta_{y'} \varphi^k_{\rx}=-A_1 \partial_{x_n} \varphi^k_{\rx} +o(1) \dans 
H^2_{\trho,\partial}.
\end{equation}
\item Conversely, if, for some $A_{-1}\in \AAA^{-1}$ the following equation
  holds:
 \begin{equation}
\label{t.t.3}
\varphi^k_{\rx}=A_{-1}\partial_{x_n}\varphi^k_{\rx}+o(1) \dans 
H^2_{\trho,\partial}.
\end{equation}
Then:
\begin{equation}
\label{t.t.4}
\tag{\ref{t.t.3}'}
\partial_{x_n} v_{Tn\rx}^k=-\Delta_{y'} A_{-1}v_{Tn\rx}^k+o(1) \dans 
L^2_{\trho,\partial}.
\end{equation}
\item Moreover if, in addition to (\ref{t.t.3}), $\teta'\neq 0$ and 
$\sigma(A_{-1})(\trho)\neq 0$, then:
\begin{gather}
\label{t.t.5}
\tag{\ref{t.t.3}''}
v^k_{T\rx}=\mathbf{Z}_{-1}\partial_{x_n} v^k_{T\rx}+o(1) \dans 
\Hv^1_{\trho,\partial}\\
\notag
\mathbf{Z}_{-1} \in \AAAv^{-1},\; \sigma_{-1}(\mathbf{Z}_{-1})=
\|\eta'\|^{-2}\left(
\begin{array}{ccc}
\begin{array}{c} 0\\0 \end{array} & \begin{array}{c} 0\\0 \end{array} 
&
i{g'^{-1}} 
\left[\begin{array}{c} \eta_1\\ \eta_2 \end{array} \right] \\
0 & 0 & \sigma(A_{-1})^{-1} 
\end{array}
\right). 
\end{gather}
\end{itemize}
\end{lem}
\begin{proof}
First note that the Dirichlet condition on $u^k$ implies:
\begin{equation}
\label{vTn..phi}
\partial_{x_n} v_{Tn\rx}^k= \Delta_{y'}\varphi^k_{\rx}+O(1) \dans 
H^{1/2}_{\loc}(\partial X),
\end{equation}
Indeed, the equation $\div u_T^k=0$ implies, in local coordinates:
$$ \partial_{x_n}v_{Tn}^k+\partial_{y_1} 
v_{T1}^k+\partial_{y_2}v_{T2}^k=O(1) \dans H^1_{\loc}(X).$$
Which may be written, using $v^k_{n\rx}=0$:
$$\partial_{x_n}v^k_{Tn\rx}=\partial_{y_1}v_{L1\rx}^k+
\partial_{y_2}v_{L2\rx}^k+O(1) \dans H^{1/2}_{\loc}(\partial X),$$
yielding (\ref{vTn..phi}) by the definition of 
$\varphi$.\par
Assume (\ref{t.t.1}). By (\ref{vTn..phi}) and the nullity of
$v^k_{n\rx}$:
$$ \Delta_{y'} \varphi^k_{\rx}=-A_1 
v_{Ln\rx}^k+o(1)=-A_1\partial_{x_n}\varphi_{\rx}^k+o(1) \dans \Lv_{\trho,\partial}.$$
Now assume (\ref{t.t.3}). Hence:
$$\Delta_{y'} 
\varphi^k_{\rx}=\Delta_{y'}A_{-1}\partial_{x_n}\varphi^k_{\rx}+o(1) 
\dans L^2_{\trho,\partial}.$$
Which implies
(\ref{t.t.4}) using (\ref{vTn..phi}) on the left side of the equation, and
the Dirichlet condition on $v_m^k$ on its right side.\par
If, in addition to the assumption (\ref{t.t.3}), $\teta'$ and 
$\sigma(A_{-1})(\trho')$ are non zero, 
both operators $A_{-1}$ and
$\Delta_{y'}$ are elliptic at 
$\trho$, and equations (\ref{t.t.4}) and (\ref{vTn..phi}) may be rewritten:
\begin{gather}
\label{transfert.1}
v_{Tn\rx}^k=Y_{-1}\partial_{x_n}v_{Tn\rx}^k+o(1) \dans 
L^2_{\trho,\partial} \quad 
\sigma(Y_{-1})=\|\teta'|\|^{-2}\sigma(A_{-1})^{-1},\\
\notag
\varphi^k_{\rx}=E_{-2}\partial_{x_n}v^k_{Tn\rx}+O(1) \dans 
H^{3/2}_{\trho,\partial},\quad \sigma(E_{-2})=-\|\eta'\|^{-2} \text{ 
near } \trho\\
\label{transfert.2}
\left( \begin{array}{c}
v_{T1\rx}^k\\
v_{T2\rx}^k
\end{array}\right)=
g'^{-1} \left( \begin{array}{c}
\partial_{y_1}\varphi_{\rx}^k\\
\partial_{y_2}\varphi_{\rx}^k
\end{array}\right)=
Z_{-1} \partial_{x_n} v_{Tn\rx}^k +O(1) \dans H^{1/2}_{\trho,x_n=0}\\
\notag
\sigma(Z_{-1})=-i\|\eta'\|^{-2}g'^{-1} 
\left( 
\begin{array}{c}
\eta_1\\
\eta_2
\end{array}
\right)
\text{ near } \trho.
\end{gather}
Assertion (\ref{t.t.5}) is an easy consequence of (\ref{transfert.1}) and (\ref{transfert.2}).
\end{proof}
We may now study cases (\ref{bord.mu.T}) and
(\ref{bord.mu.L}) of lemma \ref{lem.CS.bord}.\par
\underline{case (\ref{bord.mu.T})}: assume $\mu_L=0$, $\mu_T\neq 0$
near $\trho$. There are three possibilities:
\begin{itemize}
\item If $\trho \in \HHH_L$, the nullity of $\mu_L$ implies, by standard
  hyperbolic theory:
$$ u^k_{T\rx}=-u^k_{L\rx}\tendk 0 \dans \Hv^1_{\trho,\partial}.$$
So near $\trho$, the support of $\mu_T$ propagates. It is easy to see that
condition (\ref{supp.muT}) on the support of $\mu_T$ implies, if $\mu_T$
does not vanish, that this is one of the four cases described in the
(\ref{bord.mu.T}) of lemma \ref{lem.CS.bord}.
\item Assume $\trho\in \EEE_L$. The standard elliptic theory (proposition
  \ref{traces.ell}) implies:
$$ \partial_{x_n} \varphi_{\rx}^k +\Xi \varphi^k_{\rx}\tendk 0 \dans 
H^1_{\trho,\rx}.$$
With lemma \ref{transfert.trace}, this yields the following equation on the
traces of $u^k_T$:
\begin{equation}
\label{eq.uT.El}
u_{T\rx}^k=i\widetilde{\Zbf}_{-1} D_{x_n} u_{T\rx}^k+o(1) \dans 
\Hv^1_{\trho,\partial}, \quad \widetilde{\Zbf}_{-1}=\chi'^{-1} \Zbf_{-1} \chi'
\end{equation}
where the principal symbol of the operator $\Zbf_{-1}\in\AAAv^{-1}$ is given by
(\ref{t.t.5}). Notice that the eigenvalues of $\sigma_{-1}(i\Zbf_{-1})$,
thus those of $\sigma_{-1} (i\widetilde{\Zbf}_{-1})$ are pure
imaginary numbers. As a consequence, the boundary condition
(\ref{eq.uT.El}) is an uniform Lopatinsky boundary condition near $\trho$
(see the example following definition \ref{def.Lopa}), which shows again
the propagation of $\mu_T$. As in the case where $\trho\in \HHH_L$, it is easy
to see that this is one of the four cases of lemma \ref{lem.CS.bord}.

\item The case $\trho\in \GGG_L$ is the most difficult. When $\mu_T$ is
  non-null $\trho$, must be in $\HHH_T$. We use a contradiction argument to
  prove the propagation of the support of $\mu_T$. Let
  $\dot{\rho}\in\HHH_T$ such that the ray coming in $\dot{\rho}$ is in the
  support of $\mu_T$, but not the ray going out of $\dot{\rho}$. According
  to the standard hyperbolic theory (proposition \ref{traces.hyp.0}): 
$$ D_{x_n} u^k_{T\rx} -\Lambda_T u^k_{T\rx} \tendk 0 \dans 
\Lv_{\trho,\partial},\quad \sigma_1(\Lambda_T)=\sqrt{\nu_T\tau^2-\|\eta'\|^2}. 
$$ 
This implies, by lemma \ref{transfert.trace}, a boundary equation on
$\partial_t \varphi^k$, of the following form:
$$ \partial_t \varphi_{\rx}^k=Y_{-1} D_{x_n} \partial_t\varphi^k_{\rx}+o(1) \dans 
H^1_{\dot{\rho},\partial},$$   
which is an uniform Lopatinsky condition near $\dot{\rho}$ because
$\trho\in \GGG_L$. In view of proposition \ref{Lopa.traces} on traces in
the glancing region, such an equation implies, with the nullity of $\mu$
near $\dot{\rho}$ the following conditions: 
$$ \partial_t \varphi^k_{\rx} \rightarrow 0 \dans H^1_{\drho,\partial},\quad\partial_n \partial_t \varphi^k_{\rx} \rightarrow 0 \dans L^2_{\drho,\partial} .$$
The operator $\partial_t$ being elliptic at $\drho$, this shows that
$u^k_L$ tends to $0$ in $\Hv^1_{\drho,\partial}$, and thus:
$$ u^k_{T\rx} \tendk 0 \dans \Hv^1_{\drho,\partial}.$$ 
Hence the propagation of the support of $\mu_T$ near $\drho$, which
contredicts the assumption on rays coming in and going out of $\drho$.\par
Similar arguments show that if the ray going out of $\drho$ is
in the support of $\mu_T$, so is the ray coming in $\drho$. This proves that
the support of $\mu_T$ propagates near $\trho$. Notice that this is
necessarily the $(\HHH_T\deux)$ case.
\end{itemize}
\begin{rem}
Case (\ref{bord.mu.T}), which appears in the study of linear
thermoelasticity, was precisely described in \cite{BuLe99}. The authors
show a result of propagation of $\mu_T$, determining all the
characteristic elements of this propagation, which gives in particular the
polarization properties of $\mu_T$. In the case of the system of
magnetoelasticity, the polarization causes no problem by and 
it suffices to show the
propagation of the support of $\mu_T$ (or that of $\mu_L$ in case (2)). As mentionned
in the introduction, one may consider that the only component of $\mu_T$
and $\mu_L$ which is resistant to the dissipation is the component parallel to
$\Bbf$, cancelling the quantity $u\wedge \Bbf$.
\end{rem}
\underline{Case \ref{bord.mu.L}:} we assume now that $\mu_T=0$ and $\mu_L
\neq 0$. We argue in a similar way, considering three possibilities:
\begin{itemize}
\item If $\trho\in \HHH_T$, the standard hyperbolic theory gives an
  approximate boundary equation on $u^k_T$, which implies:  
$$ u^k_{L\rx} \tendk 0 \dans \Hv^1_{\trho,\partial}.$$
As a consequence, the support of $\mu_L$ propagates. On this support,
$\eta$ is parallel to $\Bv$, which shows, as stated in lemma
\ref{lem.CS.bord}, that $\tilde{\eta'}=0$ and
$\trho\in \HHH_L$, or $\tilde{\eta'} \sslash
 \Bv$ and  $\trho\in
\GGG_L$. 
\item If $\trho\in \EEE_T$, we write (as in the similar situation when  
 $\trho \in \EEE_L$), the boundary equation of the elliptic region:
$$ D_{x_n} u_{T\rx}^k+\Xi_T u_{T\rx}^{k} \tendk 0 \dans 
\Lv_{\trho,\partial}.$$
This implies in view of lemma \ref{transfert.trace} an uniform Lopatinsky
boundary equation on $\partial_t \varphi^k$, thus the propagation of the
support of $\mu$, which is the same as that of $\mu_L$. The fact that
$\eta'\neq 0$ shows that $\trho$ cannot be hyperbolic for the longitudinal
wave (in this case outgoing and incoming directions are not parallel, thus
at least one is not parallel to $\Bv$). Consequently, $\trho\in
\GGG_L$. More precisely, it is a diffractive point: the bicharacteristic
passing through $\trho$ must stay parallel to $\Bv$, thus its direction is
constant which is not possible for gliding rays because $\Omega$ has no
contact of infinite order with its tangents. We are in the 
$(\GGG_L)$ case of lemma \ref{lem.CS.bord}.   
\item If $\trho \in \GGG_T$, then $\trho\in \HHH_L$. The fact that $\trho\in
\GGG_T$ implies that $\teta'\neq 0$ , so directions $\teta_L^+$ and 
$\teta_L^-$ cannot be both parallel to $\Bv$. Consequently, the support of
$\mu_L$ is an union of only ingoing rays (or only outgoing rays). This
gives a boundary equation of the following form:
$$ \partial_{x_n} \partial_t \varphi^k_{\rx} + \iota \Lambda_L 
\partial_t \varphi^k_{\rx} \tendk 0 \dans \Lv_{\trho,\partial},$$
where $\iota\in\{+1,-1\}$. 
Notice that $\partial_t$ is elliptic at $\trho$, so that we may
rewrite this last property taking out all the $\partial_t$ and with $\Hv^1$
instead of $\Lv$. This yields, in view of lemma \ref{transfert.trace}, a
unform Lopatinsky boundary condition on $u^k_T$. The nullity of $\mu_T$
gives as before (by proposition \ref{Lopa.traces}):

$$ u^k_L=-u^k_T \tendk 0 \dans \Hv^1_{\trho},$$
so that $\mu_L$ propagates, and in view of the particular form of its
support, vanishes near $\trho$. This shows that this particular situation
($\trho\in \GGG_T$ and $\mu_L\neq 0$) is impossible, and completes the
proof of lemma \ref{lem.CS.bord}.  
\end{itemize}
\end{proof} 
\begin{Def}
We shall call $\Bbf$-{\bf admissible} points the point of the boundary
of $\Melrose$ which are of one of the eight types described in lemma  \ref{lem.CS.bord}.
\end{Def}
\subsection{Conclusion of the proof}
\label{CS.conclu}
Let  $S\egaldef \supp \mu_L\cup\supp\mu_T$, and $\MR$ the subset of
$\Melrose$, of all points $\rho$ satisfying one of the following properties:
\begin{itemize}
\item $x_n>0$, $\rho\in S\widehat{Z}_T$ and $\eta \sslash \Bv$;
\item $x_n>0$, $\rho\in S\widehat{Z}_L$ and
 $\eta \bot \Bv$;
\item $x_n=0$ and $\rho$ is $\Bbf$-admissible. 
\end{itemize}
Let $\Phi_T(\rho,s)$ and $\Phi_L(\rho,s)$ be the bicharacteristic flows for
the transversal and longitudinal waves. We shall define a local continuous
flow on $\MR$, denoted by:
$$\Phi(\rho,s)=\left(\Phi_{x'},\Phi_{x_n},\Phi_{\xi'},\Phi_{\xi_n}\right),$$
in the following way:
\begin{itemize}
\item if $\Phi_{x_n}(\rho,s)>0$ and $\Phi(\rho,s)\in S\widehat{Z}_T$, or if
  $\Phi(\rho,s)$ is a $\Bbf$-admissible boundary point of the form
(\ref{bord.mu.T}) of lemma \ref{lem.CS.bord}, $\Phi$ is near $(\rho,s)$ the
restriction to $\MR$ of the transversal bicharacteristic flow;

\item if $\Phi_{x_n}(\rho,s)>0$ and $\Phi(\rho,s)\in S\widehat{Z}_L$, or if
  $\Phi(\rho,s)$ is a $\Bbf$-admissible boundary point of the form
(\ref{bord.mu.L}) of lemma \ref{lem.CS.bord}, $\Phi$ is near $(\rho,s)$ the
restriction to $\MR$ of the longitudinal bicharacteristic flow;
\item if $\Phi(\rho,s)$ is $\Bbf$-admissible of type $\TL$, then:
\begin{align*}
\Phi(\rho,r)&=\Phi_T(\rho,r),\text{ if } r<s\\
\Phi(\rho,r)&=\Phi_L(\rho,r),\text{ if } r>s;
\end{align*}
\item if $\Phi(\rho,s)$ is $\Bbf$-admissible of type $\LT$, then:
\begin{align*}
\Phi(\rho,r)&=\Phi_L(\rho,r),\text{ if } r<s\\
\Phi(\rho,r)&=\Phi_T(\rho,r),\text{ if } r>s.
\end{align*}
\end{itemize} 
In view of lemma \ref{CS.bord}, $S$ (which a subset of $\MR$) is
stable under
the flow $\Phi$ on $(0,T)$. Furthemore, if for some $\rho\in \MR$,
$$ \Phi(\rho,s)\underset{s\rightarrow \tilde{s}}{\longrightarrow} \trho \notin \MR,\; \tilde{s}\in (0,T),$$
(thus $\trho$ is a boundary point which is not $\Bbf$-admissible),
then $\rho$ is not in $S$. Consequently, $S$ is an union of
$\Bbf$-resistant rays of life-length $T$. The assumption of non-existence of such
rays made in proposition \ref{propCS} shows that $S$ is empty, which
completes the proof.
\section{Necessary condition}
\label{chap.CN}
\begin{prop}
\label{PCN}
Assume that for all $T>0$, there exists a
$\Bbf$-resistant ray of life-length $T$.
Then for all $T>0$, there exists a sequence $(u^{k})$ of solutions of the
Lam\'e system such that:
\begin{align}
\label{condS1}
\|\partial_t u^{k}_{\restriction t=0}\|^2_{\Lv}+\|u^{k}_{\restriction t=0}\|^2_{\Hv^1_0} \tendk 1\\
\label{condS2}
\|u^{k} \wedge \Bbf\|_{\Hv^1((0,T)\times\Omega)} \tendk 0
\end{align}
\end{prop}
\begin{corol}
\label{corolCN}
Under the assumption of proposition \ref{PCN}, the energy of the solutions
of the magnetoelasticity equations does not decay uniformly. In other
terms, the necessary condition of theorem \ref{th1} holds. 
\end{corol}
Corollary \ref{corolCN} is a direct consequence of point b) of proposition \ref{CNSineg}.
\begin{proof}[Proof of proposition \ref{PCN}]
This proof is very much inspired by that of theorem 4 of \cite{BuLe99}.
Denote by $Z$ one of the indices $T$ or $L$ and set
$P_Z\egaldef \Delta-\nu_Z^2\partial_t^2$. We start by an elementary remark:
\begin{rem}
\label{remCN}
If $E$ is a vector subspace of $\CC^3$ and $\pi_E$  the orthogonal
projection on $E$, the defect measure of  $\pi_E u^{k}_Z$ is $\pi_E^*
\mu_Z \pi_E$. Furthermore, $P_Z \pi_E u^{k}_Z=0$, so that theorem
of propagation \ref{thm.propa} holds for the measure $\pi_E^*
\mu_Z \pi_E$ if an uniform
Lopatinsky
boundary condition holds on $\pi_E u^k_Z$. Notice that any {\bf scalar} uniform Lopatinsky condition on $u^k_Z $
yields such a condition on $\pi_E u^k_Z$. If $\pi_E^* \mu_Z \pi_E=\mu_Z$,
the measure $\mu_Z$ will be said to be {\bf polarized along $E$}. If $E$ is
the line generated by a vector $H$ of $\CC^3$ we shall also use the phrase
``polarized along $H$''. If both measures $\mu_T\indic_{]-\eps,T+\eps[}$ and
$\mu_L \indic_{]-\eps,T+\eps[}$ are polarized along $\Bbf$, then condition 
(\ref{condS2}) is fullfilled.
\end{rem}
Let $T'>T$. Consider a $\Bbf$-resistant ray defined on an open interval $I$
of length $T'$:
$$ \gamma(s)=\left(t_\gamma(s),y_\gamma(s),\tau_\gamma(s),\eta_\gamma(s)\right)=(x_\gamma(s),\xi_\gamma(s)) $$
 If $T'$ is large enough, then one of the two following assertions holds:
\begin{enumerate}
\item[a)]$\gamma(I)$ contains an interior point;
\item[b)]$y_\gamma(I)=\Gamma\subset \partial \Omega$ 
where $\Gamma$ is a closed curved, contained in a plane $P$ which is normal
to $\Bbf$, boundary of a convex subset of $P$, and such that on $\Gamma$,
$\bfn$ is orthogonal to $\Bbf$. 
\end{enumerate}
Case b) occures when there exists an infinite boundary $\Bbf$-resistant
ray. This case reduces to case a), choosing a transversal ray contained in
$P$ which only meets the boundary at hyperbolic points.\par
Thus, we may assume that $\gamma(I)$ has an interior point. We may also
assume, possibly moving the origin of coordinates, that this interior point
is $\gamma(0)$, and that $(t_\gamma(0),y_{\gamma}(0))=(0,0)$.
Recall that the magnetic field is vertical: $\Bbf=(B,0,0)$. We shall denote
by $-T^-$ and $T^+$ the extremal points of $I$: $I=(-T^-,T^+)$.\par
If $\eta_\gamma(0)$ is parallel to $\Bbf$ (i.e. if $\gamma(0)$ is in the
longitudinal characteristic set), choose a non-zero function $\varphi\in
C^\infty_0(\Omega)$, and set:
\begin{align*}
\varphi^{k}(y)=K^{-1} k^{-5/4}e^{iky_1}\varphi(\sqrt{k}y)\\
u_0^{k}=\nabla \varphi^{k}, \quad u_1^{k}=i k c_L u_0^{k}.
\end{align*}
Where $u^{k}$ is the solution of the Lam\'e system with initial data:
$$(u^{k},\partial_t u^{k})_{\restriction t=0}=(u_0^{k},u_1^{k}).$$
Then:
$$ \|u_0^k\|_{\Hv^1} \tendk K^{-1} \|\varphi\|_{\Lv},\quad
\|u_1^k\|_{\Lv}\tendk K^{-1}\|\varphi\|_{\Lv}.$$
Thus, condition (\ref{condS1}) is fullfilled with an appropriate choice of $K$.\par
For small $t$, by finite speed of propagation for the wave equation,
$u_L^{k}$ has compact support in $\Omega$. 
Thus
$u_T^{k}=0$ and $u_L^{k}=u^{k}$. As a consequence, for small $t$:
\begin{enumerate}
\item $\mu_T=0$;  
\item $\mu_L$ is polarized along $\Bbf$;
\item the projection of the support of $\mu_L$ on $\RR_t\times
  \overline{\Omega}$ is contained in $x_{\gamma}(I)$.
\end{enumerate} 
If $\eta_{\gamma}(0)$ is orthogonal to $\Bbf$, we construct 
a sequence of solution of the Lam\'e system, with the following initial
data: (cf \cite{BuLe99})
\begin{align*}
\psi^{k}=K^{-1} k^{-5/4} e^{iky_1} & \psi(\sqrt{k}y)\\
u_0^{k}=\rot(0,0,-\psi^{k}), & \quad u_1^{k}=ikc_T u_0^{k}. 
\end{align*}
In this case, condition $(\ref{condS1})$ is fullfilled for an appropriate
$K$ and the defect measures satisfy the following properties for small $t$: 
\begin{enumerate}
\item $\mu_L=0$;
\item $\mu_T$ is polarized along $\Bbf$;
\item  the projection of the support of $\mu_T$ on $\RR_t\times
  \overline{\Omega}$ is contained in $x_{\gamma}(I)$.
\end{enumerate} 
To show (\ref{condS2}), we shall prove that both measures $\mu_T$ and
$\mu_L$, are, for $t\in I$, polarized along $\Bbf$. For $t>0$, we shall denote
by $\PPP(t)$ the following property:
{\bf in a neighbourhood of $[0,t]$, both measures $\mu_L$ and $\mu_T$
  are polarized along $\Bbf$ and the projections of their support on
  $\RR_t\times \overline{\Omega}$ are contained in  $x_{\gamma}(I)$}.\\
Let $\TTT$  be the set of $t$ in $[0,T^+)$ such that $\PPP(t)$
holds. By its definition, $\TTT$ is an open subset of $[0,T^+)$. We have just
shown that $0$ in $\TTT$. We shall now prove that $\TTT$ is closed. We
shall use the next lemma:
\begin{lem}
\label{propanul}
Let $\trho=(\tilde{t},\ty,\tilde{\tau},
\teta)\in S_b^*M$. If $\mu_T \indic_{t<\tilde{t}}$ and $\mu_L
\indic_{t<\tilde{t}}$ vanish in a neighbourhood of $\trho$, so do both measures
$\mu_T$ and $\mu_L$.
\end{lem}
This is a trivial assertion in the interior of $\Omega$ by the propagation of
both measures. Near a point of the boundary of $\Omega$, one may show lemma
\ref{propanul} using the Dirichlet boundary condition on $u^k$
and the theorem of propagation $\ref{thm.propa}$, together with the same
type of arguments as in lemma \ref{lem.CS.bord}.\par
Let $s_0>0$ such that $\PPP(s_0)$ holds for $s<s_0$. We must check that
$\PPP(s_0)$ holds. Three cases arise, depending on the nature of
$\rho\egaldef \gamma(s_0)$.\par
{\bf i)} $\rho$ is an interior point.\par
 $\PPP(s_0)$ is obvious in view of the propagation of both measures in the
 interior of $\Omega$.\par
{\bf ii)} $\rho$ is of the type (\ref{bord.mu.T}) of lemma \ref{lem.CS.bord}.\par
This case, where $\mu_L$ vanishes for time $t<t_{\gamma}(s_0)$ near
$t_\gamma(s_0)$, were studied in \cite{BuLe99}. The authors show that
$\mu_L$ remains null for times greater than $t_{\gamma}(s_0)$ and that
$\mu_T$ propagates near $\gamma(s_0)$, in such a way that in our case, its
polarization along $\Bbf$ is preserved. In particular property $\PPP(s_0)$ holds. 

{\bf iii)} $\rho$ is of the type $(\HHH_L)$ of lemma \ref{lem.CS.bord}: $y'_\gamma(s_0) \in \partial \Omega$, $\eta'_{\gamma}(s_0)=0$ and ${\bf n} \sslash \Bbf$.\par
In view of lemma \ref{propanul}, the support of the measure $\mu_L$ is
contained, near $\rho$ in the union of the longitudinal ray coming in
$\rho$ and the ray going out of $\rho$. The support of $\mu_T$, if not
empty, if the transversal ray going out of $\rho$. Let $E$ be the
plane orthogonal to $\Bbf$ in $\CC^3$. The polarization of $\mu_L$ along
$\Bbf$ shows that  $\pi_E^*\mu_L\pi_E=0$ and thus, by remark \ref{remCN}
and the standard hyperbolic theory of proposition \ref{traces.hyp}:
$$\pi_E u^{k}_{T\restriction x_n=0}=-\pi_E u^{k}_{L\restriction x_n=0} \CVF 0 \text{ in }\Hv^1_{\rho,\partial}$$
which implies, using again remark  \ref{remCN} that
$\pi^*_E \mu_T \pi_E$ propagates along the transversal flow near
$\gamma(s_0)$. Thus $\pi^*_E\mu_T\pi_E$ vanishes near $\rho$. But $\mu_T$
is polarized orthogonally to its direction of propagation which
is exactly $\Bbf$ on the support of $\mu_T$ near $\rho$. This show that
$\mu_T$ vanishes near $\rho$, completing the proof of $\PPP(s_0)$.\par
{\bf iv)} $\rho$ is of the type ($\GGG_L$) of lemma \ref{lem.CS.bord}: $\rho$ is a
diffractive point for the longitudinal wave, and $\eta'_\gamma(s_0)$ is
parallel to $\Bbf$.\par
Then $\rho\in \HHH_T\cup \EEE_T$. Furthermore, in the longitudinal
hyperbolic case,  $\mu_T \indic_{t<t_{\gamma}(s_0)}$ vanishes near
$\rho$. Thus, according to standard elliptic or hyperbolic theory (cf
propositions \ref{traces.hyp.0} and \ref{traces.ell}), $u_T^k$ satisfies a
boundary condition of the following form: 
\begin{gather*}
D_{x_n} u^{k}_{T\restriction
  x_n=0}=A u_{T\restriction x_n=0}^{k}+o(1) \dans {\Lv_{\rho,\partial}},\quad A\in \AAAv^1_{\partial}\\
\sigma_{1}(A)=-i\sqrt{\|\eta'\|^2-\nu_T^2\tau^2} \text{ in the elliptic case,}\\
 \sigma_1(A)=-\sqrt{\nu_T^2 \tau^2-\|\eta'\|^2} \text{ in the hyperbolic case.}\end{gather*}
Each of this equation yields, in view of lemma \ref{transfert.trace}, an
uniform Lopatinsky boundary equation on $\varphi^k$:
\begin{equation}
\label{tr iv)}
\varphi^{k}_{\rx}=B_{-1} D_{x_n} \varphi_{\rx}^{k}+o(1) \dans \Hv^1_{\rho,\partial}.
\end{equation}
As a consequence, the support of $\mu$ (and that of $\mu_L$) propagates
near $\rho$. The polarization of $\mu_L$ along $\Bbf$ is immediate. The
nullity of $\mu_T$ near $\rho$ remains to be checked. This is a general 
property in the elliptic case $\rho \in \EEE_T$. In the hyperbolic case,
first note that the propagation theorem of Burq and Lebeau
\cite[th. 1]{BuLe99} implies with boundary condition (\ref{tr iv)})
that $\mu$ is invariant by the longitudinal flow near diffractive
points. So the total mass of $\mu_L$ is preserved by time, for $t$ close
enough to $t_{\gamma}(s_0)$. The next lemma, which is a measure
version of the conservation of energy for the Lam\'e system,
completes the proof of $\PPP(s_0)$:
\begin{lem}
\label{conservation}
Let $\varphi\in C_0^{\infty}(\RR)$. Then:
$$ \big<\mu_T+\mu_L,\varphi'(t)\big>=0.$$
In other terms, the total mass of the measure $(\mu_T+\mu_L)_{\restriction
  t=s}$ is well defined, and does not depend on $s$.
\end{lem}  
\begin{proof}
\begin{gather*}
(\partial_t^2-\Delta_e)u^k=0\\
\Re \int \partial_t^2 u^k \overline{u}^k \varphi(t) dx-\Re \int \Delta_e u^k \overline{u}^k \varphi(t) dx =0.
\end{gather*}
Set:   $\nabla_e u\egaldef (\mu \nabla u,(\lambda+\mu)\div u)\in \CC^4$. A
simple integration by parts yields:
\begin{equation*}
\int \varphi'(t) \left|\partial_t u^k\right|^2 dx+\int \varphi'(t) \left| \nabla_e u^k \right|^2 dx. 
\end{equation*}
Using another integration by parts, and then the decoupling lemma \ref{lem.decouplage}: 
\begin{gather*}
-\int \partial_t^2 u^k \overline{u}^k \varphi'(t) dx-\int \Delta_e u^k \overline{u}^k \varphi'(t) dx=o(1) \text{ as } k\rightarrow +\infty\\
\int(\partial_t^2+c_T^2\Delta)u_T^k \overline{u}^k_T \varphi'(t) dx+\int(\partial_t^2+c_T^2\Delta)u_L^k \overline{u}^k_L \varphi'(t) dx=0.
\end{gather*}
When $k$ tends to $\infty$, we get, by the definition of $\mu_T$ and $\mu_L$:
$$
\big<\mu_T, \frac{\tau^2+c_T^2 \|\eta\|^2}{2\tau^2}\varphi'(t) \big>+\big<\mu_L, \frac{\tau^2+c_L^2 \|\eta\|^2}{2\tau^2}\varphi'(t) \big>=0.
$$
This completes
 the proof, noting that on the support of $\mu_T$
(respectively $\mu_L$), $c_T \|\eta\|$ (respectively $c_L \|\eta\|$) is
equal to $\tau$.
\end{proof}
Lemma \ref{conservation} and the mass conservation of $\mu_L$ as time
goes by imply that the mass of $\mu_T$ is also preserved near $\rho$, which
shows that $\mu_T$ vanishes in a neighbourhood of $\rho$.\par
{\bf v)} $\rho$ is of the type $(L\rightarrow T)$ of lemma \ref{lem.CS.bord}.\par
In view of lemma \ref{propanul} and of the assumption $P(s)$ for
$s<s_0$, the support of $\mu_L$ is contained in the two longitudinal half-rays
passing through $\rho$, and that of $\mu_T$ is only contained in the ray
going out of $\rho$. To prove $P(s_0)$, it remains to show that $\mu_L=0$
along the longitudinal ray going out of $\rho$. We shall do so by a simple
polarization argument. Let $H$ (respectively $J$) be an unitary vector of
$\CC^3$ parallel to the direction of the transversal (respectively
longitudinal) ray going out of $\rho$. The
polarization of $\mu_T$ shows that $\pi_H^* \mu_T \pi_H$ is null, so that 
$\pi_H^* \mu_L \pi_H$ propagates near $\rho$. Now, $H$ is orthogonal to the
direction of the longitudinal ray coming in $\rho$, so that $\pi_H^* \mu_L
\pi_H$ vanishes, along incoming rays but also, in view of the propagation,
along outgoing rays. Furthermore $\mu_L\indic_{t>t_{\gamma}(s_0)}$ is
polarized along $J$. It is easy to see that this last measure
vanishes. Indeed, the polarization of $\mu_L$ implies:
$$\indic_{t>t_{\gamma}(s_0)}\pi_J^* \pi_H^* \pi_J^* \mu_L  \pi_J \pi_H
\pi_J=\indic_{t>t_{\gamma}(s_0)}\pi_J^* \pi_H^* \mu_L \pi_H \pi_J=0$$
But: 
$$\pi_J \pi_H \pi_J=\frac{<H,J>^2}{|H|^2 |J|^2} \pi_J$$
Noting that $H$ and $J$ are not orthogonal this yields the nullity of 
$\mu_L\indic_{t>t_{\gamma}(s_0)}$ in a neighbourhood of $\rho$. 
{\bf vi)} $\rho$ is of the type  $(T\rightarrow L)$ of lemma \ref{lem.CS.bord} .\par
One may argue as before, showing that for every vector $K$ orthogonal to
$J$, $\pi_K^* \mu_T \pi_K=0$ near $\rho$, which implies the nullity of
$\mu_T \indic_{t>t_{\gamma(s_0)}}$ near $\rho$.\par
The proof is completed by reversing time, which yields $\PPP(s)$ for $-T^-<t<0$.
\begin{figure}
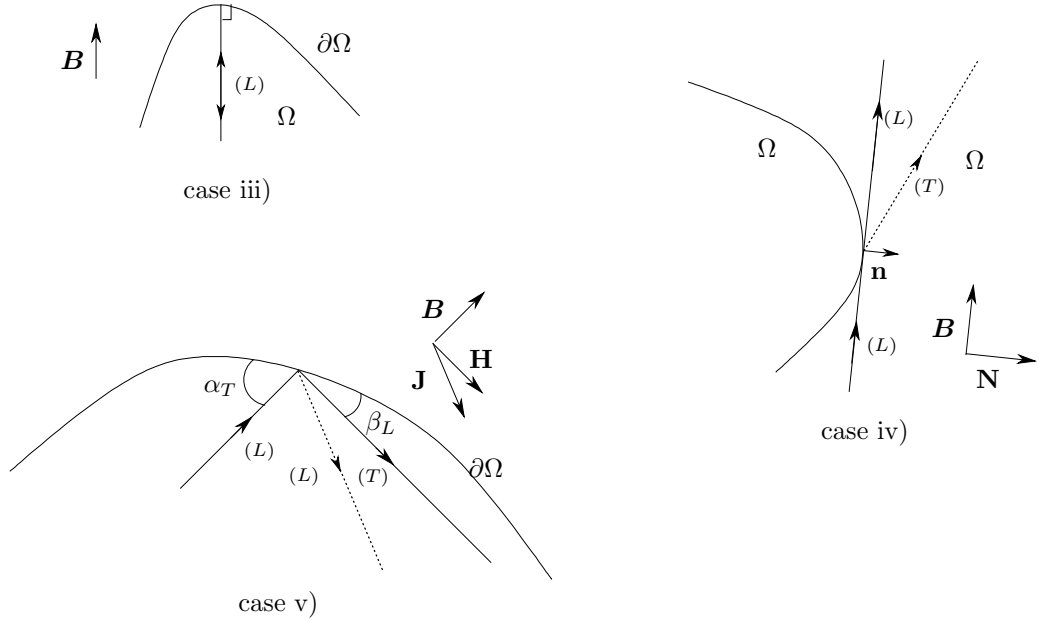

\begin{center}
\dessCN
\end{center}
\caption{Cases iii), iv) and v)}
\label{figureCN}
\end{figure}
\end{proof}
\section{Polynomial decay}
\label{chap.poly}
We shall prove in this section the following proposition:
\begin{prop}
\label{prop.inegpolyB}
There exist $T,C>0$ such that for every solution of the Lam\'e system with
initial data:
$$(u_0,u_1)\in D(\LLL^N),$$
we have:
\begin{equation}
\label{inegpolyB}
\|u_0\|^2_{\Hv_0^1}+\|u_1\|^2_{\Lv}\leq C\left(\QQQ^N_T(u)+\|u_0\|^2_{\Lv}+\|u_1\|^2_{\Hv^{-1}} \right),
\end{equation}
where:
\begin{itemize}
\item $N=1$ if there does not exist on $\Omega$ any {\bf boundary}
  $\Bbf$-resistant ray of infinite life-length; 
\item $N=K$ elsewhere.
\end{itemize}
\end{prop}
The integer $K\geq 2$ was introduced in theorem
\ref{thdecpoly}:
\begin{equation*}
K\egaldef \sup_{j=1..m} \aleph(\Gamma_j),
\end{equation*}
where $\Gamma_1$,...,$\Gamma_m$ are the spatial images of infinite
boundary $\Bbf$-resistant rays and $\aleph(\Gamma_j)$ is the minimal order of
contact of $\partial \Omega$ with its tangents
parallel to $\Bbf$ at points of $\Gamma_j$. The quadratic form $\QQQ_N^T$ is defined in section \ref{chap.A.F} by:
$$\QQQ^N_T(u)\egaldef \sum_{l=0}^N \|\rot(\partial_t^{l+1} u\wedge
\Bbf)\|^2_{\Hv^{-1}((0,T)\times\Omega)}.$$ 
Inequality (\ref{inegpolyB}) is precisely the sufficient condition of
polynomial decay given by proposition \ref{decpoly}, which completes the
proof of theorem \ref{thdecpoly}. As in section \ref{chap.CS}, we shall
argue by contradiction, using the defect measures of section \ref{chap.mesures}.
\subsection{Introduction of measures}
\label{poly.intro}
First note that by a density argument, it suffices to show
(\ref{inegpolyB}) for initial generated by a finite number of eigenfunctions
of $\LLL$. Assume that (\ref{inegpolyB}) does not hold. This yields a sequence
$u^k$ of smooth solutions of the Lam\'e system such that:  
\begin{equation}
\label{inegabsurdepoly}
 1=\|u_0^k\|^2_{\Hv^1_0}+\|u_1^k\|^2_{\Lv}>k\left( \QQQ_T^N(u^k)+\|u_0^k\|_{\Lv}^2+\|u_1^k\|_{\Hv^{-1}}^2\right) 
\end{equation}
As in section \ref{chap.CS}, we may assume that $u^k$ converges weakly to
$0$ in $\Hv^1_{\loc}(\RR\times \overline{\Omega})$ and introduce the defect
mesures of subsection \ref{par.Lame}. Note that (\ref{inegabsurdepoly}) implies:
\begin{equation}
\label{rot.CVF.0.poly}
\QQQ_T^N(u^k) \tendk 0.
\end{equation}
As a consequence, condition (\ref{rot.CVF.0}) of section \ref{chap.CS} is
fulfilled. In this section, we proved (cf remark \ref{rem.CS}) that this
condition implies that in the interval $(0,T)$, the supports of $\mu_T$ and
$\mu_L$ are unions of $\Bbf$-resistant rays.\par
We would like to show, as in section \ref{chap.CS}, that both measures
$\mu_L$ and $\mu_T$ are null, which would contredict
(\ref{inegabsurdepoly}). We shall first prove by invariance arguments that,
as long as $N\geq 1$, both measures are null in the interior of $\Omega$
(subsection \ref{poly.interieur}), which yields proposition
\ref{prop.inegpolyB} in the favorable case of non-existence of boundary
$\Bbf$-resistant rays of infinite life length. Subsection \ref{poly.bord}
studies those boundary rays, using traces theorem to restrict
(\ref{rot.CVF.0.poly}) to the
boundary of $\Omega$.\par 
In all this section, $z=(z_1,z_2,z_3)$ denotes a global spatial orthonormal
coordinate system of $\RR^3$, in which $\Bbf=(B,0,0)$. As before, notations
$y$ and $x=(t,y)$ shall only be used for local coordinates.   
\subsection{Nullity of $\mu_T$ and $\mu_L$ in the interior of $\Omega$}
\label{poly.interieur}
\begin{lem}
\label{lem.poly.int}
Let $(u^k)$ be the sequence of solutions of Lam\'e system introduced in the
preceding subsection, with $N\geq 1$. Then:
\begin{itemize}
\item $\mu_L\indic_{(0,T)}$ does not depend, in the interior of $\Omega$, on
  the variables $z_2$ and $z_3$;
\item $\mu_T\indic_{(0,T)}$ does not depend, in the interior of $\Omega$, on
  the variable $z_1$.
\end{itemize}
\end{lem}
\begin{corol}
\label{corol.poly.int}
Under the assumptions of lemma \ref{lem.poly.int} and if $T$ is large enough:
\begin{gather}
\label{muL.poly.0} \mu_L\big(\{ t\in]0,T/2[\}\big)=0\\
\label{muT.int.poly.0} \mu_T\big(\{t\in]0,T/4[,\;z\notin \partial \Omega\}\big)=0. 
\end{gather}
In particular, if there does not exist any boundary $\Bbf$-resistant ray of infinite
life-length, (\ref{inegpolyB}) holds, with $N=1$ and large enough $T$.
\end{corol}
\begin{proof}[Proof of corollary \ref{corol.poly.int}]
The second part of corollary  \ref{corol.poly.int} is an immediate
consequence of (\ref{muL.poly.0}) and (\ref{muT.int.poly.0}).\par
We shall prove those two conditions by contradiction. Assume that
(\ref{muL.poly.0}) does not hold. Then there is, in the support of $\mu_L$,
a point $\rho$ such that:
 $$ \rho=(t,z,\tau,\zeta),\quad z \in \Omega,\quad t \in ]0,T/2[.$$
Consider $D$, the half-line of origin $z$ with direction $\zeta$. Let $z_0$
be the first point where $D$ intersects the boundary, $\rho_0=(t_0,z_0,\tau_0,\zeta_0)$ the point of
the longitudinal ray coming from $\rho$ whose spatial projection is $z_0$,
and $\alpha$ the angle between $D$ and the exterior normal vector to the
boundary $\bfn$ in $z_0$ (see figure \ref{fig.muL.0}). If $T$ is large
enough (namely if $C_LT/2$ is greater than the diameter of $\Omega$),
$t_0$ is in $(0,T)$. In view of lemma \ref{lem.CS.bord}, whose assumptions
are fulfilled because
(\ref{rot.CVF.0.poly}) implies (\ref{rot.CVF.0}), one of the three
following holds:
\begin{itemize}
\item $\alpha=0$ (corresponding to the $(\HHH_L)$ case of lemma \ref{lem.CS.bord});
\item $\alpha=\alpha_L=\arctan(c_T/c_L)$ (corresponding to the
  $(L\rightarrow T)$ case of lemma \ref{lem.CS.bord}). 
\item $\alpha=\frac{\pi}{2}$ ($D$ is tangent to the boundary at $z_0$,
  which corresponds to  the $(\GGG_L)$ case of lemma \ref{lem.CS.bord}).
\end{itemize}
\begin{figure}
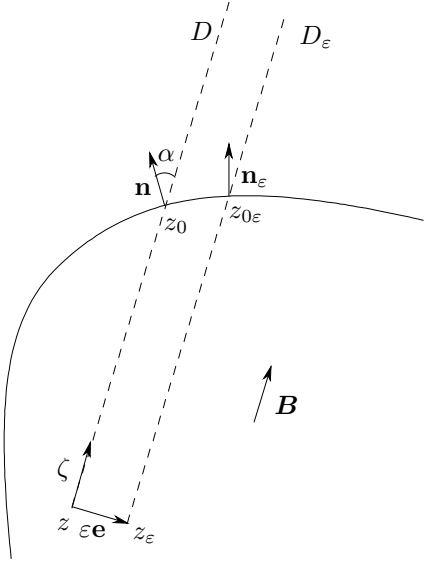

\begin{center}
\dessmuLpoly
\end{center}
\caption{Nullity of $\mu_L$.}
\label{fig.muL.0}
\end{figure}  
Let $\bf{e}$ be an arbitrary vector, orthogonal to $\Bbf$. 
By lemma \ref{lem.poly.int}, the following points of $S\widehat{Z}_L$ are
in the support of $\mu_L$: 
$$ \rho_\eps:=(t,z_{\eps},\tau,\zeta),\quad z_{\eps}\egaldef z+\eps
\bf{e},\; \eps>0 \text{ small. }$$
Consider the half-line $D_{\eps}$ going from $z_{\eps}$ in the direction
$\Bbf$, and denote by $z_{0\eps}$ its first intersection point with
$\partial \Omega$, and by $\alpha_{\eps}$ the angle between $n_{\eps}$ (the
exterior normal to the boundary in $z_{0\eps}$), and $D_{\eps}$. Then:
\begin{itemize}
\item if the intersection of $D$ with $\partial \Omega$ in $\rho$ is a
  transverse intersection (i.e. if $\alpha \neq \frac{\pi}{2}$), then:
$$ \rho_{0\eps} \underset{\eps \rightarrow 0}{\longrightarrow} \rho_0,\quad
\alpha_{\eps} \underset{\eps \rightarrow 0}{\longrightarrow} \alpha,$$
and the assumption that $\Omega$ does not have any contact of infinite
order with its tangents implies that for small, non-zero
$\eps$, $\alpha_{\eps}$ is close to, but distinct from $\alpha$, which yields: 
\begin{equation}
\label{alpha.not.in}
\alpha_{\eps}\notin \{0,\alpha_L,\pi/2\}.
\end{equation}
As a consequence, $\rho_{\eps}$ cannot be a $\Bbf$-admissible point, and
$D_{\eps}$ is not the spatial projection of a $\Bbf$-resistant ray. 
For small enough $\eps$, $t_{0\eps}$ is still in the interval $(0,T)$ and thus:
\begin{align*}
\rho_{\eps} \notin& \supp \mu_L,\\
\rho_{\eps} \notin& \supp \mu_L \text{ (by propagation along the
  longitudinal flow)},\\
\rho \notin& \supp \mu_L \text{ (by lemma \ref{lem.poly.int})},
\end{align*}
thus contradicting the definition of $\rho$;
\item if $D$ is tangent at $z$ to $\partial \Omega$, we choose $\bf{e}=\bfn$
  (which is orthogonal to $\Bbf$). For small $\eps$, the intersection of
  $D_{\eps}$ and $\partial\Omega$ is a transverse
  intersection and is close to $\rho$. Thus, if $\eps$ is small but strictly
  positive, $\alpha_{\eps}$ is close to, but
  distinct from $\pi/2$, which again shows (\ref{alpha.not.in}), and as
  before that $\rho$ is not in the support of $\mu_L$. 
\end{itemize}
\begin{figure}
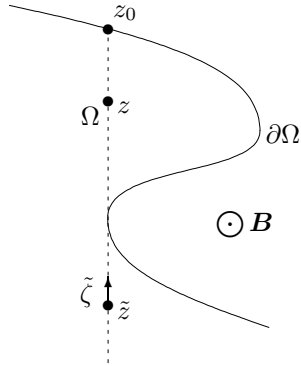

\begin{center}
\dessmuTpolyb
\end{center}
\caption{Choice of $\rho$}
\label{choix.de.z}
\end{figure}
Assume now:
$$ \mu_T\big( \{t\in ]0,T/4[, \; z\in \Omega \} \big) \neq 0.$$
Then there exists a point $\trho$ such that:
\begin{equation}
\label{poly.muT.neq.0}
 \trho=(\tz,\tilde{t},\tzeta,\tilde{\tau}) \in \supp \mu_T,\; \tilde{t}\in ]0,T/4[,\; \tz \in \Omega. 
\end{equation}
Consider a transversal ray passing  through $\trho$. It meets the boundary
at a non-diffractive point $\rho_0$, at a certain time $t_0>\tilde{t}$, after
possibly passing through diffractive points. Choose an interior point
$\rho$ which is located on this transversal ray, after all the diffractive
points, but before $\rho_0$ (see figure \ref{choix.de.z}). Thus:
$$ \rho \in \supp \mu_T,\; x \in \Omega,\; t \in ]0,T/2[.$$
The condition on $t$ holds for $T$ large enough (if suffices to take $T$ so
that the length $c_T T/4$ of a transversal ray of life length $T/4$ is greater
than the diameter of $\Omega$).\par
We may choose coordinates $(z_2,z_3)$ so that $\tzeta$ is parallel to $(0,1,0)$.
Let $\PPP$ be the plane passing through $z$ and generated by the two
orthogonal vectors $\Bbf$ and $\tzeta$. Consider  $U\egaldef \Omega \cap
\PPP$, which is an open subset of $\PPP$,
and the following family of points:
$$ \rho_{\eps}\egaldef(t,z_{\eps}=z+\eps\Bbf,\tau,\tzeta), \quad|\eps|< \eps_0.$$
For small enough $\eps_0$, $z_{\eps}$ stays in the interior of $\Omega$, so
that, in view of lemma \ref{lem.poly.int}, $\rho_{\eps}$ is in the support
of $\mu_T$. Let $\rho_{0\eps}$ be the point of
$S^*_b(\RR\times\overline{\Omega})$ where the transversal ray coming from
$\rho$ hits the boundary. Denote by $\bfn(\dot{z})$
the exterior unitary normal to $\partial \Omega$ at $\dot{z}\in
\partial\Omega$ and, if $\dot{z}\in \PPP$, by  $\bfn'(\dot{z})$ the
exterior unitary normal to $\partial U$ at $\dot{z}$. If  $\bfn(\dot{z})$
is not normal to the vector plane $\Pbf$ generated by $\dot{\zeta}$ and
$\Bbf$, the (non null) orthogonal projection of $\bfn(\dot{z})$ on $\Pbf$
is parallel to $\bfn'(\dot{z})$.\par
\begin{figure}
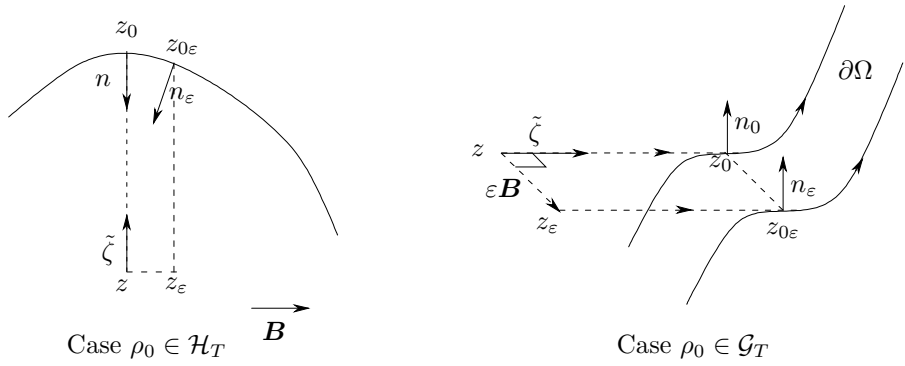

\desspolyc
\caption{Nullity of $\mu_T$ in the interior of $\Omega$}
\end{figure}
We shall note $\bfn_0,\bfn'_0,\bfn_{\eps},\bfn'_{\eps}$ instead of
$\bfn(z_0),\bfn'(z_0),\bfn(z_{0\eps}),\bfn'(z_{0\eps})$. Two cases arise:
\begin{itemize}
\item \underline{$\rho_0\in\HHH_T$}. 
The point $\rho_0$ is $\Bbf$-admissible, and ($t_0$ being in $(0,T/2)$),
$\mu_L$ vanishes near $\rho_0$ by (\ref{muL.poly.0}), which implies, if
$\mu_T$ is non null, that $\rho_0$ is of the type (\ref{bord.mu.T}) of lemma \ref{CS.bord}.
Consequently, $\bfn_0$ is orthogonal to $\Bbf$ but not to $\tzeta$ (or else
$\rho_0$ should be in $\GGG_T$). Thus the (non-null) orthogonal projection of
$\bfn_0$ on $\Pbf$ is ortogonal to $\Bbf$. Hence:
$$\bfn'_0\bot \Bbf.$$
For small enough $\eps$, $z_{0\eps}$ is close to $z_0$. Furthermore, the
preceding argument is still valid when replacing $\rho_{0}$ by
$\rho_{0\eps}$, so that $\bfn'_{\eps}$ stays orthogonal to
$\Bbf$. Consequently, $\partial U$ is, in a neighbourhood of $z_0$, a line
segment parallel to $\Bbf$, contradicting the assumption that $\Omega$ has
no contact of infinite order with its tangents.  
\item \underline{$\rho_0 \in \GGG_T$}.
By the choice of $\rho$, $\rho_0$ is not strictly diffractive and is of
the $\GGG_T\deux$ type of lemma \ref{lem.CS.bord}, which implies that
$\bfn_0$ is normal to the plane $\Pbf$. On the transversal ray coming from
$\rho_0$, whose spatial projection is a geodesic curve of $\partial
\Omega$, $\bfn$ is still orthogonal to $\Bbf$ (at least for time less
than $T$). Consider points $\rho_{0\eps}$, which are close to $\rho_0$ for
small $\eps$. If for one $\eps$ such that $|\eps|<\eps_0$,
$\rho_{0\eps}$ belongs to $\HHH_T$, we may reduce to the preceding case with
$\rho_{\eps}$ instead of $\rho$. Thus, we may assume:  
$$ \forall \eps,\; |\eps|<\eps_0 \Longrightarrow \rho_{0\eps}\in \GGG_T.$$
The same argument as before shows that along rays coming from
$\rho_{0\eps}$, the normal to the boundary $\bfn$ stays orthogonal to
$\Bbf$. This yields a small opens subset of $\partial \Omega$ in which
$\bfn$ is orthogonal to $\Bbf$, which shows that $\Bbf$ is a tangent of
infinite order to $\partial \Omega$, contradicting the assumptions on $\Omega$.
 \end{itemize}  
\end{proof}
\begin{proof}[Proof of lemma \ref{lem.poly.int}]
We may rewrite condition (\ref{rot.CVF.0.poly}):
\begin{equation}
\label{R.explicite} 
\forall l=0..N,\quad \partial_t^{l+1} \left[ 
\begin{array}{c}
-\partial_{z_2} u_2^k -\partial_{z_3} u_3^k \\
\partial_{z_1} u_2^k\\
\partial_{z_1} u_3^k
\end{array}\right]
\tendk 0 \dans \Hv^{-1}((0,T)\times \Omega).
\end{equation}
When $N\geq 1$, (\ref{R.explicite}) still holds with $l=1$. The
characteristic manifolds of $P_T$ and $P_L$ are disjoint, which shows that
(\ref{R.explicite}) holds, in the interior of $\Omega$, if one replaces
$u^k$ by $u^k_T$ or $u^k_L$. Furthermore, on each of these characteristic
manifolds, $\tau$ does not vanish so that the operator $\partial_t^2$ is
elliptic. Hence: 
\begin{gather}
\label{CV.poly.int.T}
\partial_{z_1} u_T^k \tendk 0 \dans \Hv^1_{\loc}((0,T)\times \overline{\Omega}) \\
\label{CV.poly.int.L}
\partial_{z_2} u_{L1}^k,\; \partial_{z_3} u_{L1}^k \tendk 0 \dans \Hv^1_{\loc}((0,T)\times \overline{\Omega}).
\end{gather}
We have used the nullity of $\div u_T^k$ to get (\ref{CV.poly.int.T}), and
the nullity of $\rot u_L^k$ on the two last lines of (\ref{R.explicite}) to
get (\ref{CV.poly.int.L}). We shall use properties \eqref{CV.poly.int.T},
\eqref{CV.poly.int.L} of convergence to $0$ with gain of one derivative in
a classical way (see the proof of the propagation of defect measures),
calculating the commutator of the appropriate operators  ($\partial_{z_1}$
for $u_T^k$, $\partial_{z_2}$ and $\partial_{z_3}$ for $u_L^k$) with a
``test'' pseudo-differential operator of order $2$. For example,
(\ref{CV.poly.int.T}) implies, by integration by parts:
\begin{align*}
\forall A\in \AAA^2_i,\;&([A,\partial_{z_1}] u_T^k,u_T^k) \tendk 0\\
&\big<\mu_T,\partial_{z_1}\frac{a_2}{2\nu_T^2\tau^2} \big>=0.
\end{align*} 
(We used that on the support of $\mu_T$, $\nu_T^2\tau^2=\|\eta\|^2$.)
This shows that $\mu_T$ does not depend upon $z_1$. To prove that $\mu_L$
does not depend upon $z_2$ nor $z_3$, it suffices to prove the same
property for the defect measure $\mu_{L1}$ of $u_{L1}^k$ ($\mu$ is  polarized along $\Bbf$). To do so, we
use the above argument on (\ref{CV.poly.int.L}). 
\end{proof}
\subsection{Nullity of $\mu_T$ on the boundary of $\Omega$.}
\label{poly.bord}
We now assume that there exists a boundary $\Bbf$-resistant ray, $\gamma$,
of infinite life length. It is transversal and its spatial projection lives
on a plane curve $\Gamma$, contained in the intersection of $\partial
\Omega$ with a plane $\PPP$ normal to $\Bbf$. Furthermore, all points of
$\gamma$ being gliding points, $\Gamma$ is the boundary of a convex set of
$\PPP$. We choose $T$ large enough so that in the interval $(0,T)$, the
spatial image of $\gamma$ is the entire curve $\Gamma$. We shall work near
a point $\tilde{z}$ of $\Gamma$, such that $\Bbf$ is tangent at the
order $\aleph_0=\aleph(\Gamma)$ at $\tz$ to $\partial \Omega$. This choice
is made possible by the definition of $\aleph_0$. Recall that by
assumption, $N\geq \aleph_0$. Let $\trho$ be a point in the image of
$\gamma$ whose spatial projection is $\tz$.
\begin{notations}
We choose $\tilde{z}$ as the origin of the orthonormal frame
$(0,\ebf_1,\ebf_2,\ebf_3)$, and
assume $\ebf_2$ to be tangent, in $0$, to $\Gamma$ and $\ebf_3$ equal to
the exterior unitary normal at $0$ to $\partial \Omega$. Consider the local
coordinates $(s_1,s_2)$ on $\partial \Omega$ defined by:
$$ s_1=z_1,\; s_2=z_2.$$
The assumption on $\tz$ implies that near $0$:
\begin{equation}
\label{z3.s1}
 z_3(s_1,s_2)= s_1^{\aleph_0}T_0(s_1,s_2),\quad \frac{\partial{z_3}}{\partial_{s_1}} (s_1,s_2)= s_1^{\aleph_0-1}T_1(s_1,s_2)\quad T_j(0,0)\neq 0. 
\end{equation} 
\end{notations}
The strategy of this last part of the proof is a simple one. We shall
restrict condition (\ref{rot.CVF.0.poly}) to the boundary of $\Omega$ by an
appropriate trace theorem. Such a theorem does not exist, in general, in
spaces $H^s$, $s\leq 1/2$, but in our particular case, $u^k_T$ and $u^k_L$
being  solutions of a differential
equation which is transverse to the boundary, it is 
possible to take the
trace of (\ref{rot.CVF.0.poly}) on the boundary, losing as in standard
trace theorems only half a derivative (cf proposition
\ref{traces.appendice} in the appendix). The boundary equations thus
obtained will yield the nullity of $\mu_T$, by a simple lemma giving bounds
of $L^2$ norms with loss of derivatives (lemma \ref{pertes}), and the usual
boundary equations given by standard hyperbolic, elliptic and glancing theory.
\subsubsection*{First step: restriction to the boundary.}
We shall first prove the following:
\begin{gather}
\label{traces.poly.T}
\rot(u_T^k\wedge \Bbf)_{\rbord}=B \left( 
\begin{array}{c}
-\partial_{z_2} u_{T2}^k-\partial_{z_3} u_{T3}^k\\
\partial_{z_1} u_{T2}^k\\
\partial_{z_1} u_{T3}^k
\end{array}
\right){}_{\rbord}  \tendk 0 \dans \Hv^{N-1/2}_{\loc}((0,T)\times\partial\Omega) \\
\label{traces.poly.L}
\rot(u_L^k\wedge \Bbf)_{\rbord}=B \left( 
\begin{array}{c}
-\partial_{z_2} u_{L2}^k-\partial_{z_3} u_{L3}^k\\
\partial_{z_1} u_{L2}^k\\
\partial_{z_1} u_{L3}^k
\end{array}
\right){}_{\rbord}
\tendk 0 \dans \Hv^{N-1/2}_{\loc}((0,T)\times\partial\Omega).
\end{gather}
First of all we need to decouple condition (\ref{rot.CVF.0.poly}) into one
condition on the longitudinal wave and one condition on the transversal
wave. Set:
\begin{gather*}
w^k\egaldef\rot(\partial^{N+1}_t u^k \wedge \Bbf)=w^k_T+w^k_L\\
w_{T,L}^k\egaldef \rot(\partial^{N+1}_t u_{T,L}^k \wedge \Bbf).
\end{gather*}
Then:
\begin{equation*}
\nu_T^2 \partial_t^2 w^k_T - \Delta w^k_T =0,\quad
\nu_L^2 \partial_t^2 w^k_L-\Delta w^k_L =0
\end{equation*}
Adding these two equations and using (\ref{rot.CVF.0.poly}) we get:
\begin{equation}
\label{poly.somme.1}
\nu_T^2 \partial_t^2 w_T^k+\nu_L^2 \partial_t^2 w_L^k \tendk 0 \dans \Hv^{-3}((0,T)\times\Omega).
\end{equation}
Furthermore, condition (\ref{rot.CVF.0.poly}) twice differentiated with
respect to time yields: 
\begin{equation}
\label{poly.somme.2}
 \partial_t^2 w_T^k+ \partial_t^2 w_L^k \tendk 0 \dans \Hv^{-3}((0,T)\times\Omega).
\end{equation}
We deduce from (\ref{poly.somme.1}) and (\ref{poly.somme.2}), $\nu_T$ and $\nu_L$ being distincts:
$$ \partial_t^2 w_{T,L}^k \tendk 0 \dans
\Hv^{-3}((0,T)\times\overline{\Omega}).$$
Both functions $w^k_T$ and $w^k_L$ being solutions of wave equations,
proposition \ref{traces.appendice} of the appendix implies:
\begin{equation*}
\rot\left( \partial_t^{N+3} u_{T,L}^k \wedge \Bbf\right)_{\rbord} \tendk 0 \dans
\Hv_{\loc}^{-7/2}((0,T)\times\partial \Omega)
\end{equation*}
Thus, if $\trho=(\tilde{t},\tilde{y}',\ttau,\teta')$ is a boundary point,
and if $\tilde{\tau}$ is non null, the ellipticity of $\partial_t^{N+3}$ yields:
\begin{equation*}
\rot\left(u_{T,L}^k \wedge \Bbf\right)_{\rbord} \tendk 0 \dans \Hv^{N-1/2}_{\trho,\partial}.
\end{equation*}
It remains to prove the same property when $\ttau=0$. In this case
$\trho\in \EEE_T\cap\EEE_L$. Let $M>0$. The standard elliptic theory
(propositions \ref{traces.ell.M} for the first line and 
\ref{traces.ell} for the following) implies:
\begin{gather}
\label{ell.poly.1}
\partial_{x_n} u^k_{\rx} \tendk 0 \dans \Hv^{M}_{\trho,\partial}\\
\label{ell.poly.2}
D_{x_n} u^k_{T\rx} +\Xi_{T} u^k_{T\rx} \underset{k\rightarrow 
+\infty}\longrightarrow 0, \dans \Hv^M_{\trho,\partial} \\
\label{ell.poly.3}
D_{x_n} u^k_{L\rx} +\Xi_{L} u^k_{L\rx} \underset{k\rightarrow 
+\infty}\longrightarrow 0
\dans \Hv^M_{\trho,\partial}\\
\notag
\Xi_{T,L} \in \AAA^1,\quad \sigma_1(\Xi_{T,L}) 
=i\sqrt{\|\eta'\|^2-\nu_{T,L}\tau^2} \text{ near } \trho.
\end{gather} 
Adding (\ref{ell.poly.2}) and (\ref{ell.poly.3}), and then using
(\ref{ell.poly.1}), we get:
$$ \Xi_T u_T^k+\Xi_L u_L^k \tendk 0 \dans \Hv^{M}_{\trho,\partial}.$$
This shows, using the Dirichlet boundary condition on $u^k$ and the ellipticity
of the operator $\Xi_T-\Xi_L$ that:
$$ u_{T,L}^k \tendk 0 \dans \Hv^{M}_{\trho,\partial}.$$ 
In view of (\ref{ell.poly.2}) and (\ref{ell.poly.3}), the same property
holds in the space $\Hv^{M-1}$ on $\partial_{x_n} u_{T,L}^k$, which
completes the proof of (\ref{traces.poly.T}) and (\ref{traces.poly.L}).
\subsubsection*{Second step}
We shall now prove that under the assumptions of subsection
\ref{poly.intro}, $\mu_T$ is polarized along $\Bbf$:
$$ \pi_B \mu_T\pi_B\indic_{t\in (0,T)}=\mu_T\indic_{t\in(0,T)}.$$
In other terms, denoting by $u_j$ the $jth$ component of $u$ in the basis $(\ebf_1,\ebf_2,\ebf_3)$: 
$$ u_{T2}^k \underset{k\rightarrow +\infty}{\longrightarrow} 0,\quad u_{T3}^k \underset{k\rightarrow +\infty}{\longrightarrow} 0.$$
First of all, we prove the following:
\begin{equation}
\label{dn.u2.u3}
 \partial_n u_{j\rbord}^k \tendk 0 \dans L^2_{\trho,\partial},\; j=2,3. 
\end{equation}
The sum of the last line of relations (\ref{traces.poly.T}) and (\ref{traces.poly.L}) yields:
\begin{equation}
\label{poly.u2.1}
\partial_{z_1} u_2^k \tendk 0 \dans H^{N-1/2}_{\loc}\left((0,T)\times \partial \Omega\right)
\end{equation}
On the other hand:
$$ \partial_{s_1}u_{2\rbord}^k=\partial_{z_1} u_{2\rbord}^k+\frac{\partial z_3}{\partial s_1} \partial_{z_3} u_{2\rbord}^k.$$ 
The Dirichlet boundary condition on $u^k$ allows us to take out, in the
preceding equality, all tangential derivatives. In view of (\ref{z3.s1}),
and noting the $\frac{\partial}{\partial n}$-component
of $\frac{\partial}{\partial z_3}$ does not vanish at $0$, we deduce from
(\ref{poly.u2.1}) that in a neighbourhood $U_0$ of $0$ in $\partial\Omega$:
\begin{equation}
\label{s_1^L.dnu2}
s_1^{\aleph_0-1} \partial_n u_{2\rbord}^k \tendk 0 \dans H^{N-1/2}_{\loc}(U_0).\end{equation}
Lemma \ref{pertes} below allows us to take out the factor
$s_1^{\aleph_0-1}$ in \ref{s_1^L.dnu2}, in return for the loss of
$\aleph_0-1$ derivatives:
\begin{lem}
\label{pertes}
Let $r> -1/2$, $p>0$, $d\in \NN^*$ and $f\in H^{r+p}(\RR^d)$, with compact
support. Let $y=(y_1,..,y_d)$ be the canonical coordinates on $\RR^d$. Then: 
$$ \|f\|_{H^r} \leq C\|y_1^p f\|_{H^{r+p}}. $$
\end{lem} 
\begin{proof}
It is sufficient to prove the inequality with $p=1$ and $f\in
C^{\infty}_0(\RR^d)$. Set:
\begin{equation*}
N(f)\egaldef -\text{Re}\, \int (1+|\eta|^2)^r \eta_1 \frac{\partial \hat{f}}{\partial \eta_1} \Bar{\Hat{f}} \,d\eta
\end{equation*}
Cauchy-Schwarz inequality implies:
\begin{equation}
\label{Nf.1}
|N(f)|\leq \| f \|_{H^{r}} \| y_1 f\|_{H^{r+1}},
\end{equation}
A simple integration by parts yields:
\begin{align}
\notag
N(f)=&-\frac 12 \int \frac{\partial}{\partial_{\eta_1}} |\hat f|^2 (1+|\eta|^2)^r \eta_1\,d\eta \\
\notag
=&\int |\hat f|^2(1+|\eta|^2)^{r-1} \left\{\frac 12 +\frac 12 |\eta|^2+r\eta_1^2\right\}\,d\eta\\
\label{Nf.2}
N(f) \geq & c_r \| f\| _{H^r}^2,\quad c_r>0.
\end{align}
To get inequality (\ref{Nf.2}), we used the assumption
$r+1/2>0$. Inequalities (\ref{Nf.1}) and (\ref{Nf.2}) yield the announced result.
\end{proof}
From (\ref{s_1^L.dnu2}), lemma \ref{pertes} and the assumption
$N-\aleph_0\geq 0$ we get (\ref{dn.u2.u3}) with $j=2$. 
A similar argument yields the same result on $\partial_n u_3^j$.\par
If $\trho \in \HHH_L$ (i.e. when $c_T>c_L$) we have:
$$u_{L\rbord}^k \CVF 0 \text{ in }\Hv^1_{\trho,\partial\Omega}$$
which implies by boundary Dirichlet condition on $u^k$ the same property on
$u_T^k$.\par
In the elliptic case, standard elliptic theory (lemma
\ref{traces.ell}) yields the following boundary condition:
$$ \forall j\in \{2,3\},\;D_n u^k_{Lj}=\Xi_1 u_{Lj}^k +o(1) \text{ in } L^2_{\trho}. $$
This condition is still valid with $T$ instead of $L$, by
(\ref{dn.u2.u3}) and the Dirichlet boundary condition on $u^k$.\par
Let $\mu'_T$ be the defect measure of the sequence
$(u_{T2}^k,u_{T3}^k)$. Denoting by $\pi_{23}$ the orthogonal projection of
$\CC^3$ on the plane $({\mathbf e}_2,{\mathbf e}_3)$, one may identify
$\mu_T'$ with $ \pi_{23} \mu_T \pi_{23} $. The support of $\mu_T'$ is, in a
neighbourhood of $\trho$, contained in $\GGG^T$ and its spatial projection
is contained in $\Gamma$. According to the propagation theorem of N.~Burq
and G.~Lebeau \cite[theorem 1]{BuLe99}, there exists a function $M$,
continuous (except possibly at hyperbolic points) and inversible on the
support of $\mu_T'$, such that $M^* \mu_T' M$ propagates along the
transversal flow $\Phi_T$. The {\bf scalar} boundary conditions above show
that one may also apply the propagation theorem on each component
$u_{T2}^k$ and $u_{T3}^k$, so that $M$ is of the form $m {\rm Id}_{\CC^2}$, where
$m$ is a complex valued function. Furthermore, $\mu_T$ is polarized orthogonally to
the direction of propagation, which is, along $\gamma$, the direction
tangential to $\Gamma$. As a consequence, denoting by $\pi_{\rho}$ the
orthogonal projection of $\CC^2$ on the line generated by
$(\zeta_2,\zeta_3)$, the following equality yields near $\Gamma$:
$$ \pi_{\rho} \mu_T' =\mu_T', \quad t\in(0,T)$$
by scalar propagation, we get, for  $t\in (0,T)$:
$$ \forall s,\; \pi_{\Phi_T(s,\rho)} \mu_T'=\mu_T'.$$
This is impossible unless $\mu_T'=0$. Hence:
$$ u_{Tj}^k \CVF 0 \text{ in }H^1_{\tx}((0,T)\times\Omega), \quad i\in\{2,3\}.$$
\subsubsection*{Third step} We are now able to conclude to the nullity of $\mu_T$.
The two first lines of (\ref{traces.poly.L}) may be rewritten, using the
nullity of $\rot u_L^k$:
$$ (\partial_{z_2} u_{L1}^k,\partial_{z_3} u_{L1}^k) \CVF 0 \text{ in }
H_{loc}^{N-1/2}((0,T)\times\partial \Omega ).$$
Noting that $\partial_{s_2} u=\partial_{z_2}u +\frac{\partial
  z_3}{\partial_{s_2}} \partial_{z_3}u$, we get:
\begin{equation}
\label{partial.s2}
\partial_{s_2} u_{T1}^k=-\partial_{s_2} u_{1L}^k \tendk 0 \dans
H_{\loc}^{N-1/2}((0,T)\times \partial \Omega),  
\end{equation}
which yields:
\begin{equation}
\label{u1.trace.0}
u_{T1\rx}^k \tendk 0 \dans H^1_{\trho, \partial}
\end{equation}
Consider now the first line of
(\ref{traces.poly.T}), :
\begin{equation}
\label{arbre}
\partial_{z_1} u_{T1\restriction \partial \Omega}^k \tendk 0 \text{ in } H_{loc}^{N-1/2}((0,T)\times \partial \Omega ).
\end{equation}
We have:
\begin{equation*}
\frac{\partial}{\partial s_2}\left(
  \partial_{z_1} u_{T1\rbord}^k-\frac{\partial z_3}{\partial s_1} \partial_{z_3} u_{T1\rbord}^k\right)=\partial_{s_1}\partial_{s_2} u_{T1\rbord}^k.
\end{equation*}
Together with (\ref{partial.s2}) and (\ref{arbre}) this yields:
$$ \frac{\partial}{\partial s_2}\left( \frac{\partial z_3}{\partial s_1} \partial_{z_3} u_{T1\rbord}^k\right)\tendk 0 \dans H^{N-3/2}_{\loc}((0,T)\times \partial \Omega).$$
As a consequence, in view of the expression (\ref{z3.s1}) of $z_3$, there
exists a neighbourhood $U_0$ of $0$ in $\partial \Omega$ such that:
$$s_1^{\aleph_0-1} \frac{\partial}{\partial s_2} \left( T_1(s_1,s_2) \partial_{n} u^k_{T1\rbord} \right) \tendk 0 \dans H^{N-3/2}_{\loc}((0,T)\times U_0).$$
With lemma \ref{pertes} one gets:
$$ \frac{\partial}{\partial s_2} \left( g(s_1,s_2) \partial_{n} u^k_{T1\rbord} \right) \tendk 0 \dans H^{N-\aleph_0-1/2}_{\loc}((0,T)\times U_0),$$
which yields, using the ellipticity of the operator $\partial_{s_2}$ at $\trho$:
$$ \partial_{n} u_{T1\rbord} \tendk 0 \dans H^{N-\aleph_0+1/2}_{\trho,\partial}.$$
Thus, since $N\geq \aleph_0$:
\begin{equation}
\label{dnu1.trace.0}
 \partial_n u_{T1\rbord}^k \tendk 0 \dans \Lv_{\trho,\partial}.
\end{equation}
Conditions (\ref{u1.trace.0}) and (\ref{dnu1.trace.0}) imply, by point b)
of lemma \ref{Lopa.traces}, the nullity of $\mu_{T1}$, component of $\mu_T$
along $\Bbf$. This complete the proof of the nullity of $\mu_T$ in view of our
second step above.
\section{Appendix}
\subsection{Two useful results on boundary value problems}
We shall state here two lemmas concerning solutions of a partial
diffential equation which is transverse to the boundary of an open set. The
first one says that for such functions, the control of
derivatives which are tangential to the boundary suffices to control all
the derivatives. The
second one states a trace theorem, with loss of one-half derivative, in all
$H^s$ spaces, even if $s\leq 1/2$. As in section \ref{chap.mesures}, we
shall work in an open subset $\overline{X}$ of $\overline{\RR}^{n+1}_+$, of
the form $X'\times [0,l[$, where $X'$ is an open subset of $\RR^n$.
Let $P$ be a diffential operator of degree $r$ on $\Xbar$, of the following
form: 
$$ P=\sum_{j=0..r} Q_{r-j}\partial_{x_n}^j,$$
where the $Q_j$'s are $N\times N$ matrices of tangential differential
operators, with $C^{\infty}\left(\Xbar\right)$  coefficients, and $Q_0$ is
the identity of $\CC^N$. To simplify the following statements, we suppose
$u\in C^{\infty}(\Xbar)$. 
\begin{prop}  
\label{prop1.appendice}
Let $s \geq 0$, $j\in\NN$, and suppose that
$Pu=0$ on $\Xbar$. Then:
$$ \forall \varphi \in C_0^{\infty}\left(\overline{X}\right),\quad\exists \tilde{\varphi} \in C_0^{\infty}\left(\overline{X}\right), \quad \|\varphi\partial_{x_n}^j u\|_{\Hv^{s-j}} \leq C \|\tilde{\varphi} u\|_{L^2\left(0,l;\Hv^s\left(X'\right)\right)},$$
where $C$ does not depend on $u$.
Likewise, if  again $Pu=0$ then:
$$ \forall \varphi \in C_0^{\infty}\left(\overline{X}\right),\quad\exists
\tilde{\varphi} \in C_0^{\infty}\left(\overline{X}\right),\quad \|\varphi
u\|_{L^2\left(0,l;\Hv^{-s}(X')\right)} \leq C \|\tilde{\varphi} u\|_{\Hv^{-s}(X)}.$$
\end{prop}
\begin{prop}
\label{traces.appendice}
Let $s\in\RR$, $j\in\NN$ and suppose that $Pu=0$ on
$\Xbar$. Then:
$$\forall \varphi \in C_0^{\infty}\left(\overline{X}\right),\quad\exists \tilde{\varphi} \in C_0^{\infty}\left(\overline{X}\right),\quad \|\varphi\partial_{x_n}^{j} u_{\rx}\|_{\Hv_{\loc}^{s-j-1/2}} \leq C\|\tilde{\varphi}u\|_{\Hv^s(X)},$$
where $C$ does not depend on $u$.
\end{prop}

\subsection{Proof of proposition \ref{Lopa.traces}}
As mentionned in the introduction of section \ref{chap.mesures}, whe shall
assume that each $u^k$ is smooth enough, so that all the quantities appearing
in the following calculation are well defined and finite. The general
result may be obtained with a technical smoothing argument (cf \cite[lemma 2.8]{BuLe99}).
Let:
$$ A=A_0 D_{x_n},\; A_0\in\AAA^0,\; 
\CCC^k:=\left([P,A]u^k,u^k\right).$$
Take the support of $A_0$ in a small enough neighbourhood of $\trho$. The
operator $P$ is formally self-adjoint. A simple integration by parts yields:
$$ 
\CCC^k=\underbrace{-(APu^k,u^k)+(Au^k,Pu^k)}_{0}-(Au^k_{\rx},iD_{x_n}u^k_{\rx})_{\partial}+
(iD_{x_n}Au^k_{\rx}, u^k_{\rx} )_{\partial},$$
where $(.,.)_{\partial}$ is the $\Lv$ scalar product on $\{x_n=0\}$, with
respect to the measure $\sqrt{\detg_{\rx}} dx'dt$. We have:
\begin{align*}
D_{x_n}A_0D_{x_n}u^k=& A_0 D_{x_n}^2u^k+[D_{x_n},A_0] D_{x_n} u^k\\
=& -A_0 Q u^k+R_0 D_{x_n}u^k,\quad R_0\in \AAA^0\\
\CCC^k=&i(A_0 
D_{x_n}u^k,D_{x_n}u^k)_{\partial}-i(A_0Qu^k,u^k)_{\partial}+(R_0D_{x_n}u^k,u^k)_{\partial}.
\end{align*}
So, using condition (\ref{LopaG}),
\begin{multline*}
\CCC^k=i(A_0 D_{x_n}u^k,D_{x_n}u^k)_{\partial} - i(A_0 Q B_{-1} 
D_{x_n}u^k,B_{-1}D_{x_n}u^k)_{\partial} -i (A_0 
Qh^k,B_{-1}D_{x_n}u^k)_{\partial}\\- 
i(A_0QB_{-1}D_{x_n}u^k,h^k)_{\partial}+
(R_0D_{x_n}u^k,B_{-1} D_{x_n}u^k)_{\partial}+o(1), \; 
k\rightarrow +\infty
\end{multline*}
Take $A_0$ of the form $T_0^* T_0$, where $T_0\in \AAA_0$ is scalar,
elliptic at $\trho$ and has support in a small neighbourhood of $\trho$. Then: \begin{gather*}
\CCC^k=i(E_0 D_{x_n}u^k,D_{x_n}u^k)_{\partial}-i(T_0 Q h^k,B_{-1} 
T_0D_{x_n} u^k) -i(QB_{-1}T_0 D_{x_n}u^k,T_0 h^k)+o(1),\; 
k\rightarrow +\infty\\
E_0 \in \AAA^0,\; E_0=T_0^* T_0 -B_{-1}^*T_0^* T_0 Q B_{-1}+B_{-1}^* 
R_0.
\end{gather*}
Denoting by $t_0$ the principal symbol (which is scalar) of $T_0$, we have:
$$ \sigma_0(E_0)=|t_0|^2\left(1-\transp{b}_{-1}^* 
b_{-1}q_2\right).$$
Since $q_2$ vanishes at $\trho$, we may choose $t_0$ with support in a
small enough neighbourhood of $\trho$ such that:
$$ \sigma(E_0)\geq 1/2 |t_0|^2$$
(in the sens of quadratic positive hermitian forms).
The weak G\aa rding inequality, applied to the operator $E_0-1/2 T_0^*T_0$,
thus yields:
$$ \liminf_{k\rightarrow +\infty} \Re (E_0 D_{x_n} u^k,D_{x_n}u^k)-\frac{1}{2}\| 
T_0 D_{x_n}u^k\|^2_{\Lv_{\partial}} \geq 0.$$
This implies, using the convergence to $0$ of $h^k$ in $\Hv^1_{\trho}$:
$$\liminf_{k\rightarrow +\infty} \Im \CCC^k \geq \frac 14 \|T_0 
D_{x_n}u^k\|^2_{\Lv_{\partial}}. $$
Thus $D_{x_n}u^k_{\rx}$ is bounded in 
$\Lv_{\trho,\partial}$, 
which yields, with the boundary condition, that  $u^k_{\rx}$ is bounded in
$\Hv^1_{\trho,\partial}$. The proof of (\ref{Lopa.tr.O}) is complete.\par
When $\mu$ is null near 
$\trho$, we have:
\begin{equation*}
\lim_{k\rightarrow +\infty} 
\CCC^k=\big<\mu,\frac{\{p,a_1\xi_n\}}{\tau^2} \big>=0,
\end{equation*}
which yields (\ref{Lopa.tr.o}).\par
Point b) of proposition \ref{Lopa.traces}  may be seen as a consequence of
the propagation theorem of \cite{BuLe99}. The assumptions (\ref{Hyp.b.G}) imply
that {\bf any} uniform Lopatinsky boundary condition  holds on the traces
of $u^k$, which shows that the measure $\mu$ propagates near $\trho$ with
any smooth multiplicative factor, which is impossible unless $\mu$ is null.

\end{document}